\documentclass[10pt,reqno]{amsart}
\usepackage[margin=1in]{geometry}
\usepackage[dvipsnames]{xcolor}
\usepackage{amsmath,amssymb,amsthm,mathtools,booktabs,bm,graphicx,float,array,enumitem,multirow,multicol,rotfloat,caption,subcaption,adjustbox,tikz,tikz-cd,pdflscape,csquotes,lscape,rotating,empheq,mathrsfs,longtable}
\usepackage[para]{threeparttable}
\usepackage[all]{xy}
\usepackage[numbers,sort&compress]{natbib}
\usepackage[colorlinks=true,linkcolor=blue!45!black,citecolor=blue!45!black,urlcolor=blue!45!black]{hyperref}
\hypersetup{pdftitle={Critical schemes, jet syzygies, and twisted periods},pdfauthor={Monica Jinwoo Kang}}
\numberwithin{equation}{section}
\theoremstyle{plain}
\newtheorem{theorem}{Theorem}[section]
\newtheorem{proposition}[theorem]{Proposition}
\newtheorem{lemma}[theorem]{Lemma}
\newtheorem{corollary}[theorem]{Corollary}
\theoremstyle{definition}
\newtheorem{definition}[theorem]{Definition}
\newtheorem{remark}[theorem]{Remark}

\setlist{itemsep=2pt,topsep=4pt}
\makeatletter
\def\l@section{\@tocline{1}{4pt plus 1pt}{0pt}{1.8em}{\bfseries}}
\def\l@subsection{\@tocline{2}{1pt}{1.8em}{2.4em}{\small}}
\renewcommand{\tocsection}[3]{\indentlabel{\@ifnotempty{#2}{\ignorespaces#2.\quad}}#3}

\makeatother

\newcommand{\Sym}{\mathfrak{S}}
\DeclareMathOperator{\tr}{tr}
\DeclareMathOperator{\Spec}{Spec}
\DeclareMathOperator{\Crit}{Crit}
\DeclareMathOperator{\Pf}{Pf}
\DeclareMathOperator{\ad}{ad}
\DeclareMathOperator{\Ad}{Ad}
\DeclareMathOperator{\Res}{Res}
\DeclareMathOperator{\gr}{gr}
\DeclareMathOperator{\Tor}{Tor}
\DeclareMathOperator{\Hess}{Hess}
\DeclareMathOperator{\rank}{rank}
\DeclareMathOperator{\diag}{diag}
\DeclareMathOperator{\wt}{wt}
\DeclareMathOperator{\HS}{HS}
\newcommand{\red}{\mathrm{red}}
\newcommand{\Eu}{\mathrm{Eu}}
\newcommand{\Kos}{\mathrm{Kos}}
\newcommand{\Fl}{\mathrm{Fl}}
\newcommand{\TP}{\mathrm{TP}}
\newcommand{\sT}{\mathsf{T}}
\newcommand{\sU}{\mathsf{U}}
\newcommand{\sV}{\mathsf{V}}
\newcommand{\ip}[2]{\langle #1, #2\rangle}
\newcommand{\Mac}{\mathrm{Mac}}

\title[Critical schemes, jet syzygies, and twisted periods]{Critical schemes, jet syzygies, and twisted periods}

\author{Monica Jinwoo Kang}
\address{Mitchell Institute for Fundamental Physics and Astronomy, Texas A\&M University, College Station, TX 77843, U.S.A.}
\email{monicak@tamu.edu}

\keywords{Critical schemes, Jacobian ideals, jet schemes, arc spaces, syzygies, polynomial twisted periods, Dunkl operators, hypergeometric polynomials, nilpotent orbits, Slodowy slices, Argyres--Douglas theories}

\begin{document}

\begin{abstract}
We study nonreduced critical schemes with prescribed support on nilpotent orbit closures and on the Kleinian singularity $\mathbb C^2/\mathbb Z_n$, and show that one specified minimal invariant syzygy among the first jet equations fixes their scheme structure. For $\mathfrak g=\mathfrak{sl}_N$ or $\mathfrak{so}_{2N}$ with Coxeter number $h$ and $p\ge2$ coprime to $h$, the selected polynomial is the coefficient of $s^{p+1}$ in $\det(1-sM)^{p/h}$, with reduced critical locus $\{M^p=0\}$, and the selecting syzygy is equivalent to the twisted period equation of Dubrovin and Feigin--Silantyev. The minimal invariant syzygies of that degree correspond to the lowest nonquadratic twisted periods, which, for $\mathfrak{so}_8$, form the standard representation of triality; in type D the Betti number alone does not determine the scheme. On $\mathbb C^2/\mathbb Z_n$, the selected polynomial is given by a terminating hypergeometric series. Restriction to Slodowy slices commutes with forming critical schemes and their nilradicals. These results determine, uniformly in rank and degree, the rings and nilradical filtrations through which the Schur and Macdonald indices of the corresponding Argyres--Douglas theories are proposed to arise as Hilbert series of arc spaces.
\end{abstract}

\maketitle

\setcounter{tocdepth}{2}
\tableofcontents
\setlength{\parskip}{4pt plus 1pt minus 1pt}

\section{Introduction}\label{sec:intro}

The critical locus of a polynomial carries more information as a scheme than as a variety. Its coordinate ring is the quotient of the polynomial ring by the first derivatives, the Jacobian ring of singularity theory \cite{AGZV}, and it may contain nonzero nilpotent elements. Passing to the reduced critical locus discards these elements and the relations they satisfy. The question we consider is how to recover this scheme structure from the reduced geometry, symmetry, and a specified relation among the first jet equations. For the families studied here, these data determine every coefficient of the polynomial and lead to closed-form formulas uniform in rank and degree. The question comes from a proposed construction of protected operator rings in four-dimensional quantum field theory, which is recalled at the end of this introduction; here it is treated as a question about critical schemes.

Let $P$ be a weighted homogeneous polynomial in $S=\mathbb{C}[u_1,\dots,u_r]$, with positive degrees $\deg u_a=d_a$, and write
\begin{equation*}
	I_P=(\partial_{u_1}P,\dots,\partial_{u_r}P),\qquad R_P=S/I_P,\qquad \Crit(P)=\Spec R_P .
\end{equation*}
Following the physics literature, we call $P$ the potential. The degrees and symmetry representations of the variables restrict the monomials of $P$, and prescribing the reduced critical locus restricts it further. In general undetermined coefficients survive. Since distinct Jacobian ideals can share a radical, these coefficients change the scheme structure, and hence the derivative relations in the arc space. The additional condition introduced in \cite{Kang:2026nge,KLSnonpoint} concerns the first jet scheme, obtained by substituting $u_a=u_{a,0}+u_{a,1}t$ with $t^2=0$ into the critical equations. At a distinguished bidegree, specified in Section~\ref{sec:universal}, the first jet equations are required to acquire a new minimal syzygy invariant under the symmetry group, that is, an identity among the defining equations that is not a combination, with polynomial coefficients, of identities of lower degree \cite{Eisenbud}. We call this condition \emph{syzygy maximization}. Its strength comes from the first jet equations themselves: they contain the Hessian of $P$, so the required identity becomes a differential equation for the potential. Solving this equation determines coefficients at arbitrarily high degree. The numerical condition suffices in the four-variable and type A families under the hypotheses below. In type D, one must also specify the form of the identity; a numerical count alone can fail to distinguish different Jacobian ideals.

We carry this out for three families. The first consists of $U(1)$-invariant polynomials in four variables, and its coefficients are those of a terminating hypergeometric series \cite{AAR}.

\begin{theorem}[Theorem~\ref{thm:A1} and Proposition~\ref{prop:A1crit}]\label{thm:mainA}
Let $n\ge4$, give $x,y,z,w$ intrinsic degrees $(n/2,n/2,1,2)$, let $x,y$ carry $U(1)$ charges $\pm1$, and let $z,w$ be neutral. Among $U(1)$-invariant potentials of degree $n+2$ whose reduced critical locus is the Kleinian surface 
\begin{equation*}
    \{w=0,\ xy+z^n=0\}\cong\mathbb{C}^2/\mathbb{Z}_n, 
\end{equation*}
the first jet scheme has a minimal $U(1)$-invariant syzygy at bidegree $(n+3,1)$ if and only if, up to normalization and a rescaling of $w$,
\begin{equation*}
	P_n=(xy+z^n)w+\sum_{k=1}^{\lfloor n/2\rfloor}C_{n,k}\,z^{n-2k}w^{k+1},\qquad
	C_{n,k}=\frac{2^k\,n!}{(n-2k)!\,[n(n-1)]^k\,k!\,(k+1)!}.
\end{equation*}
In that case, $\beta^{U(1)}_{2,(n+3,1)}(J_1(R_{P_n}))=1$. Equivalently, 
\begin{equation*}
    P_n=xyw+z^nw\,{}_2F_1\left(-\frac n2,\frac{1-n}2;2;\frac{8w}{n(n-1)z^2}\right).
\end{equation*}
In $R_{P_n}$, the nilradical is generated by $w$.
\end{theorem}

The same coefficient formula extends to $n=2,3$, whose minimal presentations require separate treatment in Section~\ref{sec:A1-min}. The apparent negative powers of $z$ in the hypergeometric expression disappear when its terminating series is expanded.

The other two families consist of invariant polynomials on the Lie algebras $\mathfrak{sl}_N$ and $\mathfrak{so}_{2N}$. Let
\begin{equation}\label{eq:introgh}
	(\mathfrak{g},h)=\big(\mathfrak{sl}_N(\mathbb{C}),\,N\big),\ N\ge3,\qquad\text{or}\qquad (\mathfrak{g},h)=\big(\mathfrak{so}_{2N}(\mathbb{C}),\,2N-2\big),\ N\ge4,
\end{equation}
so that $h$ is the Coxeter number \cite{Humphreys}, which equals the dual Coxeter number for these algebras. For $p\ge2$ with $\gcd(p,h)=1$, define
\begin{equation}\label{eq:introP}
	P_{p,\mathfrak{g}}(M)=[s^{p+1}]\,\det(1-sM)^{p/h},\qquad M\in\mathfrak{g}.
\end{equation}
Here,the determinant is taken in the defining representation, and $[s^m]$ denotes the coefficient of $s^m$ in the formal expansion at $s=0$, where the fractional power has constant term one, so no choice of branch is involved. We write $P^A_{p,N}$ and $P^D_{p,N}$ for the two families. Let $c_\mathfrak{g}$ be the smallest fundamental degree of $\mathfrak{g}$ greater than two, so that $c_\mathfrak{g}=3$ for $\mathfrak{sl}_N$ and $c_\mathfrak{g}=4$ for $\mathfrak{so}_{2N}$, and fix the invariant
\begin{equation}\label{eq:introQ}
	Q_{\mathfrak{sl}_N}=C_3,\qquad Q_{\mathfrak{so}_{2N}}=Q_4=\frac12C_4-\frac{3}{4h}C_2^2,\qquad C_k=\frac1k\tr(M^k).
\end{equation}

\begin{theorem}[Theorems~\ref{thm:typeA} and~\ref{thm:typeD}]\label{thm:mainB}
$P_{p,\mathfrak{g}}$ is a nonzero invariant polynomial of degree $p+1$, and
\begin{equation}
	\Crit(P_{p,\mathfrak{g}})_\red=\{M\in\mathfrak{g} : M^p=0\}.
\end{equation}
This is the closure of the nilpotent orbit with partition $[p^a,r]$ for $\mathfrak{sl}_N$, where $N=ap+r$, and for $\mathfrak{so}_{2N}$, where $2N=ap+r$, the closure of the orbit whose partition is the $D$-collapse of $[p^a,r]$.\footnote{A very even partition, in which all parts are even and each occurs with even multiplicity, labels two nilpotent orbits in $\mathfrak{so}_{2N}$ \cite{CM}. Since $p$ is odd, the $D$-collapse of $[p^a,r]$ always has an odd part, so it labels a single orbit.} In both cases, $0\le r<p$, and a zero part is omitted.
\end{theorem}

Here, partitions describe the Jordan type in the defining representation \cite{CM}. The $D$-collapse is the largest partition below $[p^a,r]$ in dominance order in which each even part has even multiplicity; it is computed in Section~\ref{sec:typeD}. The variety in Theorem~\ref{thm:mainB} is also the associated variety of the simple affine vertex algebra at the boundary admissible level $k=-h+h/p$ \cite{KW88,KWb,AvEM}.

We take gradients and Hessians with respect to $\ip AB=\tr(AB)$ on $\mathfrak{sl}_N$ and $\ip AB=\tfrac12\tr(AB)$ on $\mathfrak{so}_{2N}$; these conventions are recalled in equation \eqref{eq:pairings}. To state the syzygy results, let $\varepsilon_0,\varepsilon_1$ be the adjoint-valued generators of the presentation module of $J_1(R_P)$, mapped to the first jet equations 
\begin{equation*}
    F_0=\nabla P(M),\qquad F_1=H_P(M)Y,
\end{equation*}
where $Y$ is the first jet coordinate and $H_P$ the Hessian. For an invariant $Q$ of degree $q+1$ put
\begin{equation}\label{eq:introSigmaQ}
	\Sigma_Q=q\,\ip{\nabla Q(M)}{\varepsilon_1}-p\,\ip{H_Q(M)Y}{\varepsilon_0},
\end{equation}
which is a syzygy exactly when
\begin{equation*}
    q\,H_P\nabla Q=p\,H_Q\nabla P .
\end{equation*}
Let $K$ be the first syzygy module, $\mathfrak{m}$ the ideal of jet coordinates, and $G=SL_N$ or $SO(2N)$. We write $\TP_d(\mathfrak g)$ for the degree-$d$ invariant polynomials whose Cartan restrictions satisfy the twisted period equation \eqref{eq:TP}. This equation and the minimal quotient $K/\mathfrak mK$ are explained in Section~\ref{sec:prelim}.
\newpage
\begin{theorem}[Propositions~\ref{prop:identity-TP} and~\ref{prop:basicTP}, Theorems~\ref{thm:typeA},~\ref{thm:typeD}, and~\ref{thm:TPform}]\label{thm:mainC}
Let $c=c_\mathfrak{g}$.
\begin{enumerate}[label=(\roman*)]
\item An invariant homogeneous polynomial $P$ of degree $p+1$ satisfies
\begin{equation}\label{eq:introidentity}
	(c-1)\,H_P\nabla Q_\mathfrak{g}=p\,H_{Q_\mathfrak{g}}\nabla P
\end{equation}
if and only if it is a polynomial twisted period. The space of such $P$ is spanned by $P_{p,\mathfrak{g}}$. Moreover $Q_\mathfrak{g}$ is proportional to $[s^{c}]\det(1-sM)^{(c-1)/h}$.
\item For $P=P_{p,\mathfrak{g}}$ the map $Q\mapsto[\Sigma_Q]$ is an isomorphism from the space $\TP_c(\mathfrak{g})$ of degree-$c$ invariant twisted periods onto $(K/\mathfrak{m} K)^G_{(p+c-1,1)}$. One has $\TP_c(\mathfrak{g})=\mathbb{C} Q_\mathfrak{g}$, except that $\TP_4(\mathfrak{so}_8)=\mathrm{span}\{Q_4,\Pf\}$ is the standard representation of the triality group $\Sym_3$. Hence
\begin{equation}
	\beta^{G}_{2,(p+c-1,1)}\big(J_1(R_{P_{p,\mathfrak{g}}})\big)=
	\begin{cases}
		2,&\mathfrak{g}=\mathfrak{so}_8,\\ 
		1,&\text{otherwise.}
	\end{cases}
\end{equation}
\item For $\mathfrak{sl}_N$, a nonzero invariant homogeneous $P$ of degree $p+1$ whose reduced critical locus is contained in the nilpotent cone satisfies $\beta^{SL_N}_{2,(p+2,1)}\big(J_1(R_P)\big)\neq0$ if and only if it is proportional to $P^A_{p,N}$. For $\mathfrak{so}_{2N}$, the identity in equation \eqref{eq:introidentity} cannot in general be replaced by the invariant Betti number, together with the reduced critical locus: on $\mathfrak{so}_{10}$ at $p=3$ there is a one-parameter family of potentials with the same reduced critical locus and the same invariant Betti number but pairwise distinct Jacobian ideals.
\end{enumerate}
\end{theorem}
It is important to note that coprimality enters in two places. It excludes the higher-dimensional twisted period solution spaces in the uniqueness argument, and it forces critical matrices to be nilpotent in the characteristic polynomial argument.

Slodowy slices \cite{Slodowy} are transverse slices to nilpotent orbits. Under quantum Drinfeld--Sokolov reduction, associated varieties of vertex algebras are intersected with Slodowy slices \cite{ArakawaW}, and the following result shows that the critical schemes behave in the same way.

\begin{theorem}[Proposition~\ref{prop:slice} and Corollary~\ref{cor:nilradical}]\label{thm:mainD}
Let $P$ be an adjoint-invariant polynomial on a complex semisimple Lie algebra $\mathfrak{g}$, let $\mathcal{S}_e$ be the Slodowy slice through a nilpotent $e$, and let $\iota:\mathcal{S}_e\hookrightarrow\mathfrak{g}$ be the inclusion. Then, $I_{P|_{\mathcal{S}_e}}=\iota^*I_P\cdot\mathbb{C}[\mathcal{S}_e]$, so 
\begin{equation*}
    \Crit(P|_{\mathcal{S}_e})=\Crit(P)\times_\mathfrak{g}\mathcal{S}_e
\end{equation*}
as schemes, and the radical of the restricted ideal is the restriction of the radical. Consequently, the nilradical of the restricted Jacobian ring and all its powers are the images of their ambient counterparts. At every finite jet order and on the arc space, the critical scheme identity induces the fiber product identity, which is given in equation \eqref{eq:jets-fiber}.
\end{theorem}

The proofs of Theorems~\ref{thm:mainB} and~\ref{thm:mainC} rest on the theory of polynomial twisted periods. The Frobenius structure on the orbit space of a Coxeter group is built from Saito's flat coordinates \cite{Saito,Dubrovin:1993nt,Dubrovin:1994hc}. Twisted periods of this Frobenius manifold are discussed in Dubrovin's work on almost duality \cite{DubrovinAlmost}, and Feigin and Silantyev \cite{FS} classified their invariant polynomial solutions in terms of singular polynomials for Dunkl operators \cite{Dunkl,DJO}. The Cartan restrictions of our determinant formulas are the classical
polynomial twisted periods, which admit residue expressions
\cite[Section~5]{FS}. Here,we characterize their invariant extensions
through first jet identities, determine their reduced critical loci,
and compute the minimal invariant syzygies at the distinguished
bidegree. On a Cartan subalgebra, the first jet identities become differential equations for the gradient of the potential, and Weyl symmetry turns them into these Dunkl equations. Lemma~\ref{lem:pairing} shows that any two twisted periods $P$ and $Q$ of degrees $p+1$ and $q+1$ respectively satisfy 
\begin{equation*}
    qH_P\nabla Q=pH_Q\nabla P ,
\end{equation*}
so that $\Sigma_Q$ is a syzygy. For $Q=C_2$, this is the Euler identity; for $Q=P$, it is the Koszul identity; for $Q=Q_\mathfrak{g}$, it is the distinguished syzygy. Since $Q_\mathfrak{g}$ is itself the degree-$c$ member of the determinant family, the distinguished syzygy pairs each potential with the lowest non-quadratic member of its own family. The exceptional behavior of $\mathfrak{so}_8$ has the same origin. The Pfaffian is a twisted period in every rank, but its degree $N$ equals $c_\mathfrak{g}=4$ only when $N=4$, and the parameter $3/h=1/2$ is then precisely the value at which the classification of \cite{FS} admits a two-dimensional space of solutions.

The three families above originate in four-dimensional $\mathcal N=2$ superconformal field theories. Argyres--Douglas theories \cite{Argyres:1995jj,Argyres:1995xn} are the prototypical examples without a known Lagrangian description, and the families considered here, the $(A_1,A_{2n-1})$ theories \cite{Cecotti:2010fi,Xie:2012hs} and the $\mathcal{D}_p(G)$ theories \cite{Cecotti:2012jx,Cecotti:2013lda}, have been studied extensively through their protected operator sectors. The $\mathcal{D}_p(G)$ theories also serve as building blocks: gauging diagonal subgroups of their flavor symmetries produces infinite families of $\mathcal N=2$ and $\mathcal N=1$ superconformal theories, including theories with $a=c$ \cite{Kang:2021lic,Kang:2021ccs}.\footnote{For subsequent work on four-dimensional superconformal theories with $a=c$ and their deformations, see \cite{Kang:2022vab,Kang:2023dsa,Kang:2024elv,Kang:2024inx}.} The Schur sector of such a theory is described by a vertex operator algebra \cite{Beem:2013sza}, whose associated variety is conjecturally the Higgs branch \cite{Beem:2017ooy,Kang:2022zsl}. The Schur and Macdonald limits of the superconformal index count protected operators and their allowed spacetime derivatives \cite{Kinney:2005ej,Gadde:2011uv}; the Hall--Littlewood limit isolates a smaller sector whose relation to the Higgs branch requires the qualifications discussed in \cite{Kang:2022zsl}. For Argyres--Douglas theories these indices and vertex algebras have been determined in many cases; see, for example, \cite{Buican:2015ina,Creutzig:2017qyf,Song:2017oew,Song:2016yfd,Kang:2026nge,Kang:2022zsl}.

In our earlier work \cite{Kang:2025zub,Kang:2026nge,KLSnonpoint}, we proposed that the Higgs branch and protected indices arise from a single graded commutative ring equipped with a compatible filtration. Table~\ref{tab:dictionary} summarizes the correspondence; its mathematical content and the results established here are explained below.

\begin{table}[ht]
\centering\small
\renewcommand{\arraystretch}{1.18}
\begin{tabular}{@{\phantom{S}}p{0.34\textwidth}p{0.61\textwidth}@{\phantom{S}}}
	\toprule
	Protected sector data & Algebraic counterpart\\\midrule
	Higgs branch & Reduced critical locus $\Spec (R_P)_\red$\\
	Schur index & Hilbert series of the coordinate ring $J_\infty(R_P)$ of the arc space, with intrinsic and jet degrees added\\
	Macdonald index & Hilbert series of $\gr_FJ_\infty(R_P)$ for the nilradical filtration (Section~\ref{sec:filtration})\\
	Hall--Littlewood limit of the proposal & Hilbert series of $R_P/\sqrt{0}$ (Section~\ref{sec:arcs})\\
	OPE decoupling relations & Defining equations of the proposed operator ring\\
	Distinguished null relation & Prescribed invariant first jet syzygy (syzygy maximization)\\
	Nilpotent Higgsing & Restriction to a Slodowy slice (Theorem~\ref{thm:mainD})\\
	Stress tensor Schur operator & Generator of $\sqrt{0}$ for $(A_1,A_{2n-1})$ (Remark~\ref{rem:stress})\\\bottomrule
\end{tabular}
\caption{The correspondence proposed in \cite{Kang:2025zub,Kang:2026nge,KLSnonpoint,Kang:2026ney} between protected sector data and the algebraic structures studied in this paper.}
\label{tab:dictionary}
\end{table}

In this language, syzygy maximization is the commutative image of the decoupling of a null state in the operator product expansion \cite{Kang:2026nge,Kang:2026ney}. The variety in Theorem~\ref{thm:mainB} is the Higgs branch of the corresponding $\mathcal{D}_p(G)$ theory \cite{Song:2017oew,Kang:2022zsl,Kang:2026ney}, and restriction to a Slodowy slice describes nilpotent Higgsing \cite{Chacaltana:2012zy,Beem:2014rza}.

This paper is a companion to \cite{Kang:2026ney}, which develops the physical construction for the $\mathcal{D}_p(G)$ theories of ADE type, including the exceptional potentials, the index comparisons, and the relation to vertex operator algebras. The ambient potential is called universal in \cite{KLSnonpoint,Kang:2026ney}, because its restrictions to Slodowy slices supply the potentials for nilpotent Higgsings of the same theory; it must therefore be fixed before any slice is chosen. In \cite{KLSnonpoint}, the potentials were obtained member by member, by solving the coefficient equations directly. Here,we prove closed-form formulas for the $(A_1,A_{2n-1})$, $\mathcal{D}_p(SU(N))$, and $\mathcal{D}_p(SO(2N))$ families at all ranks, determine their critical schemes and distinguished minimal syzygies, prove uniqueness, and establish compatibility with nilpotent Higgsing. The coordinate weights, symmetry representations, distinguished bidegrees, and, in type D, the quartic partner are inputs to this algebraic characterization.

The index proposal seeks a geometric origin for the series that count protected operators. The Schur and Macdonald indices are defined as graded traces over a protected sector of the quantum theory. Their evaluation leads to character formulas and symmetric function expansions, or, for Lagrangian theories, to integrals over the gauge group \cite{Kinney:2005ej,Gadde:2011uv}. The proposal of \cite{Kang:2025zub,Kang:2026nge,KLSnonpoint,Kang:2026ney,SUPERPOINT} identifies these series with Hilbert series constructed from a finitely generated graded ring:
\begin{align}
	I_{\mathrm{Schur}}(q)=\HS_{q,q}\big(J_\infty(R_P)\big),\qquad I_{\Mac}(q,T)=\HS_{q,q,T}\big(\gr_FJ_\infty(R_P)\big).
\end{align}
Here, $J_\infty(R_P)$ is the coordinate ring of the arc space of $\Crit(P)$, and $F$ is induced by the nilradical filtration of $R_P$. The two occurrences of $q$ add intrinsic degree and jet order, while $T$ records the associated graded filtration degree; see Section~\ref{sec:arcs}. Thus, a series defined through the protected spectrum is to be recovered from the equations and nilpotent structure of an affine scheme.

This gives a geometric interpretation to the operator count. The reduced scheme records the Higgs branch, nilpotent classes retain additional proposed operator data, and arcs include the derivatives of the generators and their relations. Syzygies describe dependencies among the equations, and hence, the corrections needed in counting graded dimensions. The ring supplies multiplicative and homological information that the series alone does not specify. There is a concrete mathematical precedent for such a relation between arc geometry and special functions: Bruschek, Mourtada, and Schepers interpret the Rogers--Ramanujan identities through arc space Hilbert series \cite{BMS}.

The results of this paper establish the construction and characterization of the geometric objects required by the proposal for the three families above. Within the classes specified in Theorems~\ref{thm:mainA}--\ref{thm:mainC}, the reduced geometry and the prescribed first jet condition determine $P$ up to the stated normalizations and rescalings; in type D, the condition includes the fixed quartic partner $Q_4$. The essential point is that a condition at first jet order determines the entire polynomial, and thus the equations at every jet order. For the matrix families, this follows by identifying the syzygy equation with the twisted period equation on the Weyl group orbit space. Once the graded ring is determined, its nilradical fixes the filtration intrinsically. Theorem~\ref{thm:mainD} extends the construction to Slodowy slices, with the grading specified in Section~\ref{sec:slices}. The geometric series is consequently determined uniformly throughout these families and their nilpotent Higgsings.

The remaining index conjecture is an equality between these explicitly defined Hilbert series and protected traces. The present theorems prove the algebraic characterization of the schemes, rather than this equality at all orders; the proposed identification of $R_P$ with Zhu's $C_2$-algebra \cite{Zhu} also remains a separate question. One specialization already follows from the geometry: equation \eqref{eq:HL} proves that the Hall--Littlewood limit of the geometric series is the Hilbert series of the reduced ring. The same nonreduced scheme thus supplies both the Higgs branch and a definite candidate for the full protected count.

The paper is organized as follows. Section~\ref{sec:prelim} fixes the notation and conventions for jets and minimal syzygies, defines syzygy maximization, and develops the results on twisted periods used throughout. Section~\ref{sec:A1} treats the $(A_1,A_{2n-1})$ family. Sections~\ref{sec:typeA} and~\ref{sec:typeD} prove the closed-form formulas, critical loci, and uniqueness in types A and D. Section~\ref{sec:minimal} computes the minimal syzygies in type D, including the Pfaffian class for $\mathfrak{so}_8$ and its triality action. Section~\ref{sec:slices} treats Slodowy slices, the nilradical filtration, and the arc space. Section~\ref{sec:discussion} discusses consequences and open questions.

\section{Jacobian rings, jets, and twisted periods}\label{sec:prelim}

In this section we fix conventions for jet schemes and minimal syzygies, define syzygy maximization, and prove the results on polynomial twisted periods used in Sections~\ref{sec:typeA}--\ref{sec:minimal}. For the reader's convenience, Table~\ref{tab:notation} collects the notation used in more than one section.

\begin{table}[ht]
\centering
\renewcommand{\arraystretch}{1.12}
\begin{tabular}{@{\phantom{S}}p{0.23\textwidth}p{0.73\textwidth}@{}}
	\toprule
	Symbol & Meaning\\\midrule
	$I_P,\ R_P$ & Jacobian ideal and Jacobian ring\\
	$J_m(R),\ J_\infty(R)$ & Coordinate rings of the jet scheme and arc space\\
	$K,\ \mathfrak{m}$ & First syzygy module of $J_1(R_P)$; ideal of jet coordinates\\
	$\Eu,\ \Fl,\ \mathrm{Com},\ \Kos$ & Euler, flavor, commutator, and Koszul identities\\
	$\mathfrak{g},\ \mathfrak{h},\ h$ & Lie algebra, Cartan subalgebra, Coxeter number\\
	$C_k,\ Q_\mathfrak{g},\ \Pf$ & Trace invariants, distinguished partner in equation \eqref{eq:introQ}, Pfaffian\\
	$H_P,\ \mathcal{B}_x,\ \TP_d(\mathfrak{g})$ & Hessian, twisted period tensor in equation \eqref{eq:B}, twisted periods of degree $d$\\
	$\Sigma_Q$ & First jet identity in equation \eqref{eq:introSigmaQ} attached to $Q$\\
	$\Psi(s)=\sum_d\psi_ds^d,\ \chi(t)$ & $\det(1-sM)^\alpha$ and its coefficients; $\det(t-M)$\\
	$\mathbb O_p$ & Nilpotent orbit closure $\{M\in\mathfrak{g}:M^p=0\}$\\
	$\mathfrak{n}$ & Nilradical of a Jacobian ring\\\bottomrule
\end{tabular}
\caption{Notation and conventions used in this paper.}
\label{tab:notation}
\end{table}

\subsection{Jets and arcs}

Let $S=\mathbb{C}[u_1,\dots,u_r]$ with positive intrinsic degrees $\deg u_a=d_a$, and let $P\in S$ be weighted homogeneous of degree $D$. Write
\begin{equation*}
	F_a=\partial_{u_a}P,\qquad I_P=(F_1,\dots,F_r),\qquad R_P=S/I_P.
\end{equation*}
For a finitely generated algebra $R$, let $J_m(R)$ denote the coordinate ring of the $m$-th jet scheme of $\Spec R$, which parametrizes maps 
\begin{equation*}
    \Spec A[t]/(t^{m+1})\to\Spec R
\end{equation*}
for varying $\mathbb{C}$-algebras $A$; see \cite{Mustata,EM} for general properties. Its affine presentation is obtained by substituting formal series
\begin{equation}\label{eq:arcs}
	u_a(t)=\sum_{j\ge0}u_{a,j}t^j,\qquad F_{a,j}=[t^j]\,F_a(u(t)),
\end{equation}
so that the coordinate ring of the arc space is
\begin{equation*}
	J_\infty(R_P)=\frac{\mathbb{C}[u_{a,j}:1\le a\le r,\ j\ge0]}{(F_{a,j}:1\le a\le r,\ j\ge0)},
\end{equation*}
and $J_m(R_P)$ retains the coordinates and equations with $j\le m$. In particular,
\begin{equation*}
	J_1(R_P)=S_1/(F_{a,0},F_{a,1})_{a=1}^r,\qquad S_1=\mathbb{C}[u_{a,0},u_{a,1}],\qquad
	F_{a,1}=\sum_b\frac{\partial F_a}{\partial u_b}(u_0)\,u_{b,1}.
\end{equation*}
We keep intrinsic and jet degrees separate,
\begin{equation*}
	\deg u_{a,j}=(d_a,j),\qquad \deg F_{a,j}=(D-d_a,j),
\end{equation*}
and their sum is the Schur degree in the physical interpretation. With the coefficient convention in equation \eqref{eq:arcs}, the natural derivation of $J_\infty(R_P)$ sends $u_{a,j}$ to $(j+1)u_{a,j+1}$.

\subsection{Minimal syzygies}

A first syzygy is an identity among the defining relations. Such identities account for the fact that the equations does not need to impose independent constraints. To fix the indexing, take a minimal bigraded free resolution over $S_1$ \cite{Eisenbud},
\begin{equation*}
	\cdots\longrightarrow\mathcal F_2\longrightarrow\mathcal F_1\longrightarrow S_1\longrightarrow J_1(R_P)\longrightarrow0 .
\end{equation*}
The generators of $\mathcal F_1$ are minimal defining relations and those of $\mathcal F_2$ are minimal first syzygies, which are hence counted by the second Betti numbers of the quotient ring. When the $F_{a,j}$ form a minimal generating set, put
\begin{equation*}
	K=\ker\Big(\bigoplus_{a,\,j\in\{0,1\}}S_1e_{a,j}\longrightarrow S_1\Big),\qquad e_{a,j}\longmapsto F_{a,j},
\end{equation*}
and let $\mathfrak{m}\subset S_1$ be the ideal generated by all jet coordinates. Then,
\begin{equation}\label{eq:betti}
	\beta^{S_1}_{2,(a,b)}\big(J_1(R_P)\big)=\dim_\mathbb{C}(K/\mathfrak{m} K)_{(a,b)},\qquad K/\mathfrak{m} K\simeq\Tor^{S_1}_2\big(J_1(R_P),\mathbb{C}\big).
\end{equation}
The quotient by $\mathfrak{m} K$ is essential. A polynomial identity obtained by multiplying a lower degree syzygy by a coordinate contributes no new minimal generator; we call such identities \emph{inherited}. We thus distinguish throughout between the dimension of the polynomial syzygy space $K$ and the minimal Betti number. Invariant calculations are taken with respect to the complexified flavor group; in particular, the $U(1)$-invariant polynomials are the weight-zero polynomials for the corresponding $\mathbb C^*$-action. When a reductive flavor group $G_F$ acts, a singlet Betti number is the dimension of the invariant part of $K/\mathfrak{m} K$; taking invariants is exact on each finite-dimensional graded piece. We denote this invariant dimension by
\begin{equation*}
    \beta^{G_F}_{2,(a,b)}(J_1(R_P))
    :=\dim_\mathbb C\Tor^{S_1}_2(J_1(R_P),\mathbb C)^{G_F}_{(a,b)}.
\end{equation*}
Thus $G_F$ specifies the invariant part, while the resolution is always over the displayed polynomial ring $S_1$. Eliminating a variable changes that ambient ring and can change the Betti numbers.

\subsection{Universal identities and syzygy maximization}\label{sec:universal}

Three kinds of identity are present for every potential with prescribed weights and symmetry. Differentiating Euler's identity gives
\begin{equation}\label{eq:Euler}
	\Eu=\sum_ad_a\,u_{a,0}F_{a,1}-\sum_a(D-d_a)\,u_{a,1}F_{a,0}=0,
\end{equation}
of bidegree $(D,1)$. If a group acts linearly and preserves $P$, each infinitesimal generator $\xi$ gives 
\begin{equation*}
    \sum_a(\xi u_0)_aF_{a,0}=0,
\end{equation*}
together with its first jet; for the $U(1)$ of Section~\ref{sec:A1}, these are the identities $\Fl_0,\Fl_1$, and for the adjoint action on a Lie algebra, they are the commutator identities $\mathrm{Com}$, expressing 
\begin{equation*}
    [M,\nabla P(M)]=0.
\end{equation*}
Finally, any two generators give a Koszul syzygy \cite{Eisenbud}
\begin{equation*}
    F_{a,0}e_{b,1}-F_{b,1}e_{a,0}\in K;
\end{equation*}
for the matrix families, the invariant combination is $\Kos=\ip{F_0}{\varepsilon_1}-\ip{F_1}{\varepsilon_0}$, of bidegree $(2D-2,1)$. Polynomial multiples of these identities are inherited at higher degree.

The distinguished bidegrees are inputs from the operator count of \cite{KLSnonpoint,Kang:2026ney}: 
\begin{align}
\begin{aligned}
    (A_1,A_{2n-1})\ \text{family}:&\qquad (n+3,1)=(D+1,1)\quad &\text{where}\ D=n+2, \\
    \mathcal{D}_p(SU(N))\ \text{or type }A\ \text{family}:&\qquad (p+2,1)=(D+1,1)\quad &\text{where}\ D=p+1, \\
    \mathcal{D}_p(SO(2N))\ \text{or type }D\  \text{family}:&\qquad (p+3,1)=(D+2,1)\quad &\text{where}\ D=p+1. \\
\end{aligned}
\end{align}

\begin{definition}[Syzygy maximization]\label{def:maximization}
Let $G_F$ be the connected flavor group and $(a,1)$ the distinguished bidegree. A potential $P$ of the prescribed form satisfies \emph{syzygy maximization} if
\begin{equation}
	\beta^{G_F}_{2,(a,1)}\big(J_1(R_P)\big)\neq0,
\end{equation}
that is, if the first jet equations admit a flavor singlet minimal syzygy at $(a,1)$.
\end{definition}

For fixed degree, coefficient comparison turns the syzygy equation into a finite linear system. Extra polynomial identities appear when its rank drops, a condition expressed by minors. To count minimal identities one must also quotient by those inherited from lower degree. For $n\ge4$, Theorem~\ref{thm:A1} proves uniqueness under the reduced locus condition, up to the stated normalization and rescaling of $w$; Theorem~\ref{thm:typeA} gives the type A characterization under its nilpotence hypothesis. At $p=2$ in type A and $p=3$ in type D, the distinguished degree coincides with the Koszul degree, so the existence of a numerical class does not by itself impose an additional condition. In type D we require the fixed identity in equation \eqref{eq:introidentity}; Section~\ref{sec:quartic} shows why its form matters.

\begin{remark}\label{rem:count}
In \cite{Kang:2026ney}, the condition is imposed in the stronger form that the Betti number $\beta^{G}_{2,(p+c-1,1)}\big(J_1(R_P)\big)$ equals $2$ for $\mathfrak g=\mathfrak{so}_8$ and $1$ otherwise, and with $\beta^{U(1)}_{2,(n+3,1)}=1$ for the $(A_1,A_{2n-1})$ family in \cite{KLSnonpoint}. Each of these counts implies the nonvanishing in Definition~\ref{def:maximization}, so they are special cases of it. By Theorem~\ref{thm:TPform}, the potentials $P_{p,\mathfrak g}$ attain these values, which equal $\dim\TP_c(\mathfrak g)$; for the $(A_1,A_{2n-1})$ and type A families, the uniqueness theorems show that the two formulations select the same potential. For $\mathfrak{so}_8$, the count $2$ is the stronger condition: the second class is the Pfaffian class of Section~\ref{sec:Pf-class}, which by Lemma~\ref{lem:pairing} arises from a twisted period partner of $P$, which will be described in detail in Section~\ref{sec:triality}.
\end{remark}

\subsection{Polynomial twisted periods}\label{sec:TP}

Let $\mathfrak{g}=\mathfrak{sl}_N(\mathbb{C})$ with $N\ge3$ or $\mathfrak{g}=\mathfrak{so}_{2N}(\mathbb{C})$ with $N\ge4$, and let $h=N$ or $h=2N-2$ be its Coxeter number \cite{Humphreys}. We use the invariant pairings
\begin{equation}\label{eq:pairings}
	\ip AB=\tr(AB)\ \ \text{for}\;\;\mathfrak{sl}_N,\qquad \ip AB=\frac12\tr(AB)\ \ \text{for}\;\;\mathfrak{so}_{2N},
\end{equation}
with traces taken in the defining representation, and the invariants $C_k=\tfrac{1}{k}\tr(M^k)$. For an invariant polynomial $P$ let $\nabla P:\mathfrak{g}\to\mathfrak{g}$ be its gradient with respect to equation \eqref{eq:pairings} and $H_P(M)=d(\nabla P)_M$ its Hessian, a self-adjoint endomorphism of $\mathfrak{g}$.

The diagonal Cartan subalgebra of 
$\mathfrak{sl}_N$ is
\begin{equation*}
    \mathfrak{h}=\{x\in\mathbb{C}^N:\sum_ix_i=0\}.
\end{equation*}
The split Cartan subalgebra of $\mathfrak{so}_{2N}$, with eigenvalues $(x_1,\dots,x_N,-x_1,\dots,-x_N)$ in the defining representation, is $\mathfrak{h}\cong\mathbb{C}^N$. In both cases the pairing restricts to $\sum_ix_iy_i$, the positive roots are $\Phi^+=\{e_i-e_j\}_{i<j}$, respectively $\{e_i\pm e_j\}_{i<j}$, all of squared length two, and one checks directly in coordinates that
\begin{equation}\label{eq:casimir}
	\sum_{\beta\in\Phi^+}(\beta,u)\,\beta=h\,u,\qquad u\in\mathfrak{h} .
\end{equation}
Chevalley restriction identifies the invariant rings \cite{GW,Kostant}:
\begin{equation*}
    \mathbb C[\mathfrak g]^G\xrightarrow{\ \sim\ }\mathbb C[\mathfrak h]^W,
    \qquad P\longmapsto P|_{\mathfrak h},
\end{equation*}
where $W$ is the Weyl group. We write $f=P|_\mathfrak{h}$, $g=\nabla f$, and $\Hess f$ for the Cartan restriction, gradient, and Hessian, respectively. This identification lets us test invariant polynomial identities on the Cartan. The Jacobian ideal $I_P$, however, is generated by derivatives in all directions of $\mathfrak g$.

\begin{lemma}\label{lem:chevalley}
Let $P$ be invariant and $x\in\mathfrak{h}$. Then $\nabla P(x)=\nabla f(x)\in\mathfrak{h}$, the Hessian $H_P(x)$ preserves $\mathfrak{h}$, and $H_P(x)|_\mathfrak{h}=\Hess f(x)$.
\end{lemma}

\begin{proof}
The maximal torus $T$ fixes $x$. For $t\in T$, equivariance thus gives
\begin{equation*}
	\Ad_t\nabla P(x)=\nabla P(x),\qquad \Ad_t\,H_P(x)\,\Ad_t^{-1}=H_P(x).
\end{equation*} The $T$-fixed subspace of $\mathfrak{g}$ is $\mathfrak{h}$, and an operator commuting with $T$ preserves its weight spaces. Since the pairing is nondegenerate on $\mathfrak{h}$, the identities 
\begin{equation*}
    \ip{\nabla P(x)}u=df_x(u),\qquad \ip{H_P(x)u}v=d^2f_x(u,v)
\end{equation*}
for $u,v\in\mathfrak{h}$ give the remaining assertions.
\end{proof}

For $x$ in the regular part $\mathfrak{h}^{\rm reg}=\{(\beta,x)\neq0\ \forall\beta\}$ define the symmetric bilinear map
\begin{equation}\label{eq:B}
	\mathcal{B}_x(u,v)=\sum_{\beta\in\Phi^+}\frac{(\beta,u)(\beta,v)}{(\beta,x)}\,\beta\ \in\mathfrak{h} .
\end{equation}

\begin{definition}\label{def:TP}
An invariant homogeneous polynomial $P$ of degree $q+1\ge2$ is a \emph{(polynomial) twisted period} if
\begin{equation}\label{eq:TP}
	\Hess f(x)\,u=\frac qh\,\mathcal{B}_x\big(u,\nabla f(x)\big)\qquad\text{for all }x\in\mathfrak{h}^{\rm reg},\ u\in\mathfrak{h} .
\end{equation}
We write $\TP_{q+1}(\mathfrak{g})$ for the space of such polynomials.
\end{definition}

\begin{remark}\label{rem:parameter}
The factor $q/h$ is forced. If
\begin{equation*}
    \Hess f(x)u=\nu\,\mathcal{B}_x(u,\nabla f)
\end{equation*}
for a nonconstant homogeneous $f$, then $u=x$ and equation \eqref{eq:casimir} give
\begin{equation*}
    q\nabla f=\nu h\nabla f
\end{equation*}
by Euler's identity. Thus, $\nu=q/h$. Our roots have squared length two, and equations \eqref{eq:B} and \eqref{eq:TP} use the tensor and the positive sign in \cite[(2.1), (2.3)]{FS}. They are the polynomial twisted period equations for Dubrovin's almost dual structure \cite{DubrovinAlmost}. The corresponding Dunkl operators have the minus sign in \cite[(2.10)]{FS}, as displayed below.
\end{remark}

Dunkl operators combine differentiation with the Weyl reflections. They are useful here because equation \eqref{eq:TP} is exactly the condition that all components of the gradient lie in their common kernel. This brings the uniqueness question into the classification of singular polynomials. We first write the equation in coordinates. In type A we use the translation-invariant extension $\tilde f(x)=f(x-\bar x\mathbf 1)$ to $\mathbb{C}^N$, with $\bar x=N^{-1}\sum_ix_i$, so that $\sum_ig_i=0$. Let $s_{ik}$ exchange $x_i$ and $x_k$, and in type D let $\tilde s_{ik}$ send $(x_i,x_k)$ to $(-x_k,-x_i)$. The Dunkl operators \cite{Dunkl} at parameter $\nu$ are
\begin{align*}
	D_i&=\partial_i-\nu\sum_{k\neq i}\frac{1-s_{ik}}{x_i-x_k} &&(\text{for type }A_{N-1}),\\
	D_i&=\partial_i-\nu\sum_{k\neq i}\Big(\frac{1-s_{ik}}{x_i-x_k}+\frac{1-\tilde s_{ik}}{x_i+x_k}\Big) &&(\text{for type }D_N).
\end{align*}
Each numerator vanishes on its reflection hyperplane, so these operators preserve polynomials. A polynomial annihilated by all $D_i$ is called singular \cite{DJO}.

\begin{lemma}\label{lem:TPcoords}
Let $P$ be invariant and homogeneous of degree $q+1\ge2$, with $g=\nabla f$ (extended as above in type A). The following are equivalent:
\begin{enumerate}[label=(\alph*)]
\item $P\in\TP_{q+1}(\mathfrak{g})$;
\item for all $i\neq j$,
\begin{align}
	(x_i-x_j)\,\partial_ig_j&=-\frac qN\,(g_i-g_j) &&(\text{for }\mathfrak{sl}_N), \label{eq:TPA}\\
	(x_i^2-x_j^2)\,\partial_ig_j&=\frac{2q}h\,(x_ig_j-x_jg_i) &&(\text{for }\mathfrak{so}_{2N}); \label{eq:TPD}
\end{align}
\item every component $g_j$ is annihilated by all Dunkl operators $D_i$ at $\nu=q/h$.
\end{enumerate}
\end{lemma}

\begin{proof}
In type A, $e_i$ in the following calculation means its orthogonal projection $e_i-N^{-1}\mathbf 1$ to $\mathfrak h$; pairing with a root gives the same result. For $i\neq j$, the $j$-th coordinate of $\mathcal{B}_x(e_i,v)$ receives contributions only from roots with nonzero $i$ and $j$ components. In type A, this is $-(v_i-v_j)/(x_i-x_j)$; in type D, the root $e_i+e_j$ adds $(v_i+v_j)/(x_i+x_j)$, and the sum is $2(x_iv_j-x_jv_i)/(x_i^2-x_j^2)$. Hence, equation \eqref{eq:TPA} and equation \eqref{eq:TPD} are the off-diagonal entries of equation \eqref{eq:TP}. Let 
\begin{equation*}
    E=\Hess f-\frac qh\mathcal{B}_x(\cdot,g)
\end{equation*}
and suppose its off-diagonal entries vanish. In type A, $E$ extends to $\mathbb{C}^N$ with $E\mathbf 1=0$, because $\Hess\tilde f\,\mathbf 1=0$ and $(\beta,\mathbf 1)=0$; a diagonal matrix annihilating $\mathbf 1$ is zero. In type D, Euler's identity and $\mathcal{B}_x(x,g)=hg$ give 
\begin{equation*}
    Ex=qg-qg=0 ,
\end{equation*}
so each diagonal entry vanishes wherever $x_i\neq0$, hence identically. This proves (a)$\Leftrightarrow$(b).

For (b)$\Leftrightarrow$(c), Weyl invariance of $f$ gives $g(wx)=w\,g(x)$, so $s_{ik}g_j=g_{s_{ik}(j)}$, while $\tilde s_{ij}g_j=-g_i$ and $\tilde s_{ik}g_j=g_j$ for $j\notin\{i,k\}$. For $i\neq j$ only $k=j$ contributes to the reflection terms, and
\begin{equation*}
	D_ig_j=\partial_ig_j-\nu\,\frac{g_j-g_i}{x_i-x_j}\quad(\text{A}),\qquad D_ig_j=\partial_ig_j-2\nu\,\frac{x_ig_j-x_jg_i}{x_i^2-x_j^2}\quad(\text{D}).
\end{equation*}
These vanish exactly when (b) holds with $\nu=q/h$. For the diagonal operators: in type A, $\sum_jg_j=0$ gives $D_ig_i=-\sum_{j\neq i}D_ig_j$. In type D, Euler's identity, the symmetry of the Hessian, and equation \eqref{eq:TPD} give
\begin{equation*}
	x_iD_ig_i=qg_i-\sum_{j\neq i}\Big[x_j\partial_ig_j+\nu\,x_i\Big(\frac{g_i-g_j}{x_i-x_j}+\frac{g_i+g_j}{x_i+x_j}\Big)\Big]=\big(q-2\nu(N-1)\big)g_i ,
\end{equation*}
which vanishes at $\nu=q/h$.
\end{proof}

The uniqueness statements of this paper rest on the following classification.

\begin{corollary}[of Feigin--Silantyev {\cite[Theorem~4.8 and Corollary~4.9]{FS}}]\label{thm:FS}
Let $q\ge1$ with $\gcd(q,h)=1$. Then, $\dim\TP_{q+1}(\mathfrak{g})\le1$.
\end{corollary}

\begin{proof}
By Chevalley restriction and Lemma~\ref{lem:TPcoords}, $\TP_{q+1}(\mathfrak g)$ is the degree-$(q+1)$ part of the space $L_\nu$ of Weyl-invariant polynomial solutions in \cite{FS} at $\nu=q/h$. By \cite[Theorem~4.8 and Corollary~4.9]{FS}, $L_\nu$ consists of the constants unless $\nu=(d_\beta-1)/h+m$ for a fundamental degree $d_\beta$ and an integer $m\ge0$; in that case $\dim L_\nu=2$, except in type $D_N$ with $N$ even and $d_\beta=N$, where $\dim L_\nu=3$. With $\nu=q/h$ the condition reads $q+1\equiv d_\beta\pmod h$. The fundamental degrees are pairwise incongruent modulo $h$, except for the repeated degree $N$ in type $D_N$ with $N$ even. In that case $q\equiv N-1\pmod{2N-2}$, so $N-1$ divides $\gcd(q,h)$, which is impossible since $N-1\ge3$. Hence $L_\nu$ contains at most one nonconstant solution up to scale.
\end{proof}

\begin{remark}\label{rem:FSrange}
In \cite{FS}, the case $\nu\le1$ is proved in full, and the case $\nu>1$ is reduced to it through an isomorphism between the solution spaces at $\nu$ and $\nu-1$, which is sketched there. Corollary~\ref{thm:FS} enters only the uniqueness statements in Theorems~\ref{thm:mainC}(i), \ref{thm:typeA}(4), and \ref{thm:typeD}(4), and the triality invariance in Proposition~\ref{prop:triality}; for $p>h$, these rely on that step. The closed formulas, the identities they satisfy, the reduced critical loci, and the minimal syzygies are proved independently of it. For the examples with $p>h$, namely those in Table~\ref{tab:typeA} and $\mathcal{D}_7(SO(8))$, we verified directly that the solutions of equation \eqref{eq:TPA} or \eqref{eq:TPD} among invariant polynomials of degree $p+1$ form a one-dimensional space; see Table~\ref{tab:FScheck}.
\end{remark}

\begin{table}[H]
\centering
\renewcommand{\arraystretch}{1.12}
\begin{tabular}{@{\phantom{T}}lcccc@{\phantom{T}}}
\toprule
Theory & $(N,p)$ & $\nu=p/h$ & $\dim \mathbb{C}[\mathfrak g]^{G}_{p+1}$ & Solutions\\\midrule
$\mathcal{D}_5(SU(3))$    & $(3,5)$  & $5/3$  & 2 & 1\\
$\mathcal{D}_{11}(SU(3))$ & $(3,11)$ & $11/3$ & 3 & 1\\
$\mathcal{D}_7(SU(4))$    & $(4,7)$  & $7/4$  & 4 & 1\\
$\mathcal{D}_6(SU(5))$    & $(5,6)$  & $6/5$  & 3 & 1\\
$\mathcal{D}_7(SO(8))$    & $(4,7)$  & $7/6$  & 7 & 1\\
$\mathcal{D}_9(SO(10))$   & $(5,9)$  & $9/8$  & 7 & 1\\\bottomrule
\end{tabular}
\caption{Direct check of uniqueness for cases with $p>h$. For each case, equation \eqref{eq:TPA} or \eqref{eq:TPD} is solved as a linear system in the coefficients of a general invariant polynomial of degree $p+1$, written in the basis of monomials in the fundamental invariants. The fourth column is the dimension of that space, and the last column is the dimension of the solution space. The last row is an additional check beyond the examples of Table~\ref{tab:typeA} and Section~\ref{sec:rank4}.}
\label{tab:FScheck}
\end{table}

The exceptional type D space reappears for the quartic partners of $\mathfrak{so}_8$, whose parameter is $3/6=1/2$. Their multiplicity is distinct from the uniqueness of the potential at the coprime parameter $p/h$; Section~\ref{sec:triality} computes the two-dimensional quartic space directly.

\begin{lemma}[Pairing]\label{lem:pairing}
Let $P\in\TP_{p+1}(\mathfrak{g})$ and $Q\in\TP_{q+1}(\mathfrak{g})$. Then,
\begin{equation}\label{eq:pairing}
	q\,H_P\nabla Q-p\,H_Q\nabla P=0\qquad\text{on }\mathfrak{g} .
\end{equation}
Equivalently, $\Sigma_Q=q\ip{\nabla Q}{\varepsilon_1}-p\ip{H_QY}{\varepsilon_0}$ is a syzygy of the first jet equations of $P$, of bidegree $(p+q,1)$.
\end{lemma}

\begin{proof}
The left side of equation \eqref{eq:pairing} is an $\Ad$-equivariant polynomial map $\mathfrak{g}\to\mathfrak{g}$. At $x\in\mathfrak{h}^{\rm reg}$, Lemma~\ref{lem:chevalley} and equation \eqref{eq:TP} give
\begin{equation*}
	q\,\Hess f_P\,\nabla f_Q-p\,\Hess f_Q\,\nabla f_P=\frac{pq}h\Big(\mathcal{B}_x(\nabla f_Q,\nabla f_P)-\mathcal{B}_x(\nabla f_P,\nabla f_Q)\Big)=0
\end{equation*}
by the symmetry of $\mathcal{B}_x$. The map thus vanishes on the dense set of regular semisimple elements, hence identically. The image of $\Sigma_Q$ in $S_1$ is 
\begin{equation*}
    q\ip{\nabla Q}{H_PY}-p\ip{H_QY}{\nabla P}=\ip{qH_P\nabla Q-pH_Q\nabla P}{Y},
\end{equation*}
by self-adjointness of the Hessians.
\end{proof}

\begin{lemma}[Commuting Hessians]\label{lem:commute}
Let $P,Q$ be invariant and homogeneous of degrees $p+1,q+1$, and suppose $qH_P\nabla Q=pH_Q\nabla P$. Then, $(p+q)[H_P,H_Q]=0$.
\end{lemma}

\begin{proof}
Differentiate in a direction $V$:
\begin{equation*}
	q\big((dH_P)[V]\nabla Q+H_PH_QV\big)=p\big((dH_Q)[V]\nabla P+H_QH_PV\big).
\end{equation*}
The operator $V\mapsto(dH_P)[V]\nabla Q$ is self-adjoint, since $\ip W{(dH_P)[V]\nabla Q}=d^3P(W,V,\nabla Q)$ is symmetric in $V,W$; the same holds with $P,Q$ exchanged. Taking adjoints exchanges $H_PH_Q$ and $H_QH_P$. Subtracting the adjoint equation from the original gives $q[H_P,H_Q]=-p[H_P,H_Q]$.
\end{proof}

Recall $c_\mathfrak{g}$ and $Q_\mathfrak{g}$ from equation \eqref{eq:introQ}. With the pairings in equation \eqref{eq:pairings}, the gradients are
\begin{equation}\label{eq:gradQ}
\begin{gathered}
	\nabla C_3=(M^2)_0,\qquad H_{C_3}(V)=(MV+VM)_0\qquad(\text{for }\mathfrak{sl}_N),\\
	\nabla Q_4=M^3-\lambda C_2M,\qquad \lambda=\frac3h\qquad(\text{for }\mathfrak{so}_{2N}),
\end{gathered}
\end{equation}
where $A_0=A-N^{-1}(\tr A)\mathbf 1$ is the traceless part.

\begin{proposition}[The distinguished identity]\label{prop:identity-TP}
Let $P$ be invariant and homogeneous of degree $p+1\ge2$, and let $c=c_\mathfrak{g}$, $Q=Q_\mathfrak{g}$. Then
\begin{equation}
	(c-1)\,H_P\nabla Q=p\,H_Q\nabla P\quad\Longleftrightarrow\quad P\in\TP_{p+1}(\mathfrak{g}).
\end{equation}
In particular, $Q_\mathfrak{g}\in\TP_c(\mathfrak{g})$.
\end{proposition}

\begin{proof}
($\Rightarrow$) By Lemma~\ref{lem:commute}, $[H_P,H_Q]=0$. By Lemma~\ref{lem:chevalley} its restriction to $x\in\mathfrak{h}^{\rm reg}$ is $[H,K]=0$ with $H=\Hess f$ and $K=\Hess f_Q$.

In type A, work on $\mathbb{C}^N$ with the translation-invariant extensions. Then $H\mathbf 1=0$, $\mathbf 1^TH=0$, $Hx=pg$, and $K=2\Pi\diag(x)\Pi$ with $\Pi=I-N^{-1}\mathbf 1\mathbf 1^T$. Hence $H\Pi=\Pi H=H$, and
\begin{equation*}
	[H,\Pi\diag(x)\Pi]=H\diag(x)-\diag(x)H-\frac pN\big(g\mathbf 1^T-\mathbf 1g^T\big).
\end{equation*}
Its $(i,j)$ entry is $(x_j-x_i)H_{ij}-\tfrac pN(g_i-g_j)$, so $[H,K]=0$ gives equation \eqref{eq:TPA} with $q=p$.

In type D,
\begin{equation*}
    f_Q=\frac14\sum_ix_i^4-\frac\lambda4s_2^2\qquad\text{with}\qquad s_2=\sum_ix_i^2.
\end{equation*}
Then,
\begin{equation*}
	K=3\diag(x_i^2)-\lambda s_2I-2\lambda xx^T .
\end{equation*}
Using $Hx=pg$ and $x^TH=pg^T$,
\begin{equation*}
	[H,K]_{ij}=3(x_j^2-x_i^2)H_{ij}-2\lambda p\,(g_ix_j-x_ig_j),
\end{equation*}
so $[H,K]=0$ gives equation \eqref{eq:TPD} with $q=p$, since $2\lambda p/3=2p/h$. In both types Lemma~\ref{lem:TPcoords} shows $P\in\TP_{p+1}(\mathfrak{g})$.

For $P=Q$ the identity is tautological, so $Q\in\TP_c(\mathfrak{g})$. ($\Leftarrow$) now follows from Lemma~\ref{lem:pairing} with $q=c-1$.
\end{proof}

\begin{proposition}[Basic twisted periods]\label{prop:basicTP}
\leavevmode
\begin{enumerate}[label=(\roman*)]
	\item $C_2\in\TP_2(\mathfrak{g})$.
	\item For $\mathfrak{so}_{2N}$, the Pfaffian, normalized by $\Pf|_\mathfrak{h}=x_1\cdots x_N$, lies in $\TP_N(\mathfrak{so}_{2N})$ for every $N\ge4$.
	\item $\TP_3(\mathfrak{sl}_N)=\mathbb{C} C_3$; $\TP_4(\mathfrak{so}_{2N})=\mathbb{C} Q_4$ for $N\ge5$; and $\TP_4(\mathfrak{so}_8)=\mathrm{span}\{Q_4,\Pf\}$.
	\item $C_3=-\tfrac N2\,[s^3]\det(1-sM)^{2/N}$ and $Q_4=-\tfrac h6\,[s^4]\det(1-sM)^{3/h}$.
\end{enumerate}
\end{proposition}

\begin{proof}
\begin{enumerate}[label=(\roman*)]
	\item On $\mathfrak{h}$, $\nabla f=\kappa x$ and $\Hess f=\kappa I$, with $\kappa=1$ for $\mathfrak{sl}_N$ and $\kappa=2$ for $\mathfrak{so}_{2N}$. By equation \eqref{eq:casimir}, $\tfrac1h\mathcal{B}_x(u,\kappa x)=\kappa u$.
	\item With $g_j=\Pf/x_j$ one has $\partial_ig_j=\Pf/(x_ix_j)$ for $i\neq j$, while $x_ig_j-x_jg_i=\Pf\,(x_i^2-x_j^2)/(x_ix_j)$. Since $2(N-1)/h=1$, this is equation \eqref{eq:TPD} with $q=N-1$.
	\item The degree-three invariants of $\mathfrak{sl}_N$ are $\mathbb{C} C_3$. The degree-four invariants of $\mathfrak{so}_{2N}$ are spanned by $C_2^2$ and $C_4$, together with $\Pf$ when $N=4$. Proposition~\ref{prop:identity-TP} and (ii) place $C_3$, $Q_4$, and (for $N=4$) $\Pf$ in the relevant spaces. On the other hand $C_2^2\notin\TP_4$: for $f=s_2^2$ one has $g=4s_2x$, so $x_ig_j-x_jg_i=0$, while $\partial_ig_j=8x_ix_j\neq0$. Since $\TP_4$ is a linear subspace, the claims follow.
	\item Since $\log\det(1-sM)=-\sum_{k\ge1}C_ks^k$, the coefficient of $s^3$ in $\det(1-sM)^{2/N}$ is $-\tfrac2NC_3$, and in type D, where odd traces vanish, the coefficient of $s^4$ in $\det(1-sM)^{\alpha}$ is $-\alpha C_4+\tfrac{\alpha^2}2C_2^2$. At $\alpha=3/h$ this equals $-\tfrac6hQ_4$.
\end{enumerate}
\end{proof}

\begin{remark}\label{rem:tower}
Lemma~\ref{lem:pairing} gives, for a twisted period $P$ of degree $p+1$, a family of singlet syzygies $\Sigma_Q$ indexed by $Q\in\bigoplus_d\TP_d(\mathfrak{g})$. Three choices recover identities already encountered: $Q=C_2$ gives $\kappa\,\Eu$, $Q=P$ gives $p\,\Kos$, and $Q=Q_\mathfrak{g}$ gives the distinguished identity. By Proposition~\ref{prop:basicTP}(iv), the last pairs $P$ with the lowest non-quadratic member of the determinant family to which $P$ belongs.
\end{remark}

\section{The \texorpdfstring{$(A_1,A_{2n-1})$}{(A1,A2n-1)} family}\label{sec:A1}

Consider $S=\mathbb{C}[x,y,z,w]$ with intrinsic degrees
\begin{equation*}
	\deg x=\deg y=\frac n2,\qquad \deg z=1,\qquad \deg w=2,
\end{equation*}
and let $U(1)$ act on $x,y,z,w$ with charges $1,-1,0,0$.
The filtered presentation of \cite{KLSnonpoint} also assigns
\begin{equation}\label{eq:A1-qT}
	\wt_{q,T}(x)=\wt_{q,T}(y)=\big(\frac n2,\frac n2\big),\qquad \wt_{q,T}(z)=(1,1),\qquad \wt_{q,T}(w)=(2,1).
\end{equation}
Physically, $x,y,z$ represent Higgs branch operators, whose two weights coincide, while $w$ represents the Schur operator of the stress tensor multiplet. The Jacobian ideal is homogeneous in intrinsic degree but generally inhomogeneous in the second weight, which then defines a quotient filtration; Remark~\ref{rem:stress} identifies it intrinsically. The Betti numbers below use intrinsic and jet degrees, with this filtration forgotten.

We first derive the allowed form of $P$. Flavor invariance forces $x,y$ to occur through $xy$, and for $n\ge4$ the inequality $2n>n+2$ excludes $(xy)^2$. Thus,
\begin{equation*}
	P=(a_0w+b\,z^2)\,xy+\sum_{\ell=0}^{m+1}c_\ell\,z^{n+2-2\ell}w^\ell,\qquad m=\lfloor n/2\rfloor .
\end{equation*}
Every point of the prescribed reduced critical locus
\begin{equation*}
	w=0,\qquad xy+z^n=0
\end{equation*}
must be critical. The $x$-derivative on this locus gives $b=0$, the $z$-derivative then gives $c_0=0$, and 
\begin{equation*}
    \partial_wP|_{w=0,\,xy=-z^n}=(c_1-a_0)z^n
\end{equation*}
gives $c_1=a_0$. Setting $a_k=c_{k+1}$ for $k\ge1$,
\begin{equation}\label{eq:A1-form}
	P=a_0(xy+z^n)w+\sum_{k=1}^ma_k\,z^{n-2k}w^{k+1}.
\end{equation}
We require
\begin{equation}\label{eq:A1radical}
	\sqrt{I_P}=(w,\ xy+z^n),
\end{equation}
so that the reduced critical variety is the Kleinian surface $\mathbb{C}^2/\mathbb{Z}_n$ \cite{Slodowy}, the Higgs branch of the $(A_1,A_{2n-1})$ theory \cite{KLSnonpoint}.

\begin{theorem}\label{thm:A1}
Let $n\ge4$, and let $P$ have the form in equation \eqref{eq:A1-form} and satisfy equation \eqref{eq:A1radical}. Then
\begin{equation}
	\beta^{U(1)}_{2,(n+3,1)}\big(J_1(R_P)\big)=1
\end{equation}
if and only if, up to an overall nonzero scalar and a nonzero rescaling of $w$,
\begin{equation}\label{eq:Pn}
	P=P_n=(xy+z^n)w+\sum_{k=1}^mC_{n,k}\,z^{n-2k}w^{k+1},\qquad
	C_{n,k}=\frac{2^k\,n!}{(n-2k)!\,[n(n-1)]^k\,k!\,(k+1)!}.
\end{equation}
For every other $P$ of this form the singlet Betti number at $(n+3,1)$ vanishes.
\end{theorem}

The same coefficient formula defines $P_n$ for $n=2,3$; their first jet presentations are treated in Section~\ref{sec:A1-min}. To prove the theorem we first solve for the coefficients admitting an identity beyond those forced by symmetry and homogeneity, then determine the reduced critical locus, and finally show that the additional identity is minimal.

\subsection{The singlet syzygy space}

Write $F_{x,j},F_{y,j},F_{z,j},F_{w,j}$ for the first jet Jacobian generators. A coefficient multiplying $F_{a,j}$ must have the complementary intrinsic degree, jet degree, and charge. Listing the monomials with these properties gives the complete singlet ansatz at $(n+3,1)$ for $n\ge4$:
\begin{equation}\label{eq:A1-ansatz}
\begin{split}
	\Sigma={}&(\lambda_0x_1z_0+\lambda_1x_0z_1)F_{x,0}+(\lambda_2y_1z_0+\lambda_3y_0z_1)F_{y,0} +(\lambda_4w_1+\lambda_5z_0z_1)F_{z,0}\\
	&+(\lambda_6z_0^2z_1+\lambda_7w_0z_1+\lambda_8z_0w_1)F_{w,0} +\lambda_9x_0z_0F_{x,1}+\lambda_{10}y_0z_0F_{y,1}\\
	&+(\lambda_{11}w_0+\lambda_{12}z_0^2)F_{z,1}+(\lambda_{13}z_0^3+\lambda_{14}z_0w_0)F_{w,1}.
\end{split}
\end{equation}
The flavor symmetry gives
\begin{equation*}
	\Fl_0=x_0F_{x,0}-y_0F_{y,0}=0,\qquad \Fl_1=x_1F_{x,0}+x_0F_{x,1}-y_1F_{y,0}-y_0F_{y,1}=0,
\end{equation*}
and together with equation \eqref{eq:Euler} these supply three independent inherited identities at $(n+3,1)$:
\begin{equation}\label{eq:A1inherited}
	z_1\Fl_0,\qquad z_0\Fl_1,\qquad z_0\Eu .
\end{equation}
Here and below we use the name of an identity also for its coefficient vector in the presentation module.

Comparing coefficients in equation \eqref{eq:A1-ansatz} first gives
\begin{equation*}
\begin{gathered}
	\lambda_6=\lambda_{13}=0,\qquad \lambda_8=-n\lambda_{12},\qquad \lambda_9+\lambda_{10}=n\lambda_{12},\\
	\lambda_1+\lambda_3=-\lambda_7,\qquad \lambda_2+\lambda_9+\lambda_{14}=0,\qquad \lambda_0+\lambda_{10}+\lambda_{14}=0 .
\end{gathered}
\end{equation*}
Two free parameters give $z_1\Fl_0$ and $z_0\Fl_1$. After these are separated, the remaining coefficient equations involve only
\begin{equation*}
	(A,B,C,D,E,F)=(\lambda_4,\lambda_5,\lambda_7,\lambda_{11},\lambda_{12},\lambda_{14}).
\end{equation*}
Set $a_{-1}=a_{m+1}=0$. The neutral equations are
\begin{equation}\label{eq:Ok}
\begin{split}
	\mathcal{O}_k={}&(n-2k)\,a_kB+(k+1)\,a_kC+(n-2k+2)(n-2k+1)\,a_{k-1}D\\
	&+(n-2k)(n-2k-1)\,a_kE+(k+1)(n-2k)\,a_kF=0,
\end{split}
\end{equation}
\begin{equation}\label{eq:Ek}
\begin{split}
	\mathcal{E}_k={}&(n-2k)\,a_kA+(k+1)(n-2k)\,a_kD\\
	&-2(k+1)(k+2)\,a_{k+1}E+(k+1)(k+2)\,a_{k+1}F=0.
\end{split}
\end{equation}
For odd $n=2m+1$ both families occur for $0\le k\le m$; for even $n=2m$ the equation $\mathcal{E}_m$ vanishes identically. The vector
\begin{equation}\label{eq:A1-universal}
	\left(0,\ -\frac{n+1}2,\ 0,\ 0,\ \frac12,\ 1\right)
\end{equation}
always solves these equations; it is the neutral part of $\tfrac12z_0\Eu$. Let $M(a)$ be the coefficient matrix of equations \eqref{eq:Ok}--\eqref{eq:Ek}. Then,
\begin{equation*}
	\dim K^{U(1)}_{(n+3,1)}=2+\dim\ker M(a).
\end{equation*}
An additional identity exists precisely when $\dim\ker M(a)\ge2$, equivalently $\rank M(a)\le4$. We first solve this rank condition and then establish minimality.

\subsection{Solving the coefficient equations}\label{sec:A1-solve}

The reduced geometry requires
\begin{equation*}
	a_0a_m\neq0 .
\end{equation*}
If $a_0=0$ the potential loses its $x,y$ dependence, and if $a_m=0$ the entire line $x=y=z=0$ becomes critical; both contradict equation \eqref{eq:A1radical}.

Given an additional neutral solution, subtract a multiple of equation \eqref{eq:A1-universal} to set $F=0$. Then
\begin{equation*}
a_1\neq0,\qquad E\neq0 .
\end{equation*}
Suppose first that $a_1=0$. Then $\mathcal{E}_0=0$ gives $A=-D$. If $E\neq0$, the recurrence $\mathcal{E}_k=0$ forces every higher $a_k$ to vanish, contrary to $a_m\neq0$. If $E=0$, the equation $\mathcal{O}_1=0$ gives $n(n-1)a_0D=0$, hence $A=D=0$; the equations $\mathcal{O}_0=\mathcal{O}_m=0$ then force $B=C=0$, leaving only the zero solution. Thus $a_1\neq0$. Now suppose $E=0$. The equations $\mathcal{E}_0=\mathcal{E}_1=0$ give $A+D=A+2D=0$, and $\mathcal{O}_0=\mathcal{O}_1=0$ give $B=C=0$, again leaving only the zero solution.

Put $u=a_1/a_0$. The first equations give
\begin{equation}\label{eq:EBC}
\begin{gathered}
	E=\frac{n(A+D)}{4u},\qquad C=-nB-n(n-1)E,\qquad B=\frac{n(n-1)D/u-(n^2+3n-6)E}{n+2}.
\end{gathered}
\end{equation}
Since $E\neq0$, equation \eqref{eq:Ek} is a recurrence for the $a_k$. If one intermediate coefficient vanished, so would all subsequent ones, contrary to $a_m\neq0$; hence every $a_k$ is nonzero. The equation $\mathcal{E}_1=0$ gives
\begin{equation*}
	a_2=\frac{(n-2)\,a_1(A+2D)}{12E}.
\end{equation*}
Substituting this and equation \eqref{eq:EBC} into $\mathcal{O}_2$ yields
\begin{equation*}
	\mathcal{O}_2=\frac{(n-2)\,a_1}{3(A+D)}\,(A-D)\big(2A+(n+1)D\big).
\end{equation*}
All divisions are permitted, since $E\neq0$ gives $A+D=4uE/n\neq0$. Thus a nonzero additional solution has $D\neq0$; after setting $D=1$ and $A=-\sigma$ there are two possibilities,
\begin{equation*}
	\sigma=-1\qquad\text{or}\qquad\sigma=\frac{n+1}2 ,
\end{equation*}
and the recurrence becomes
\begin{equation*}
	\frac{a_{k+1}}{a_k}=\frac{2u\,(n-2k)(k+1-\sigma)}{n(k+1)(k+2)(1-\sigma)} .
\end{equation*}
For $\sigma=(n+1)/2$ this reads
\begin{equation*}
	\frac{a_{k+1}}{a_k}=\frac{2u\,(n-2k)(n-2k-1)}{n(n-1)(k+1)(k+2)},
\end{equation*}
so $a_k=a_0u^kC_{n,k}$, and a rescaling of $w$ together with an overall normalization sets $a_0=a_1=1$.

For odd $n=2m+1$, the last equation $\mathcal{E}_m=0$ forces $\sigma=m+1=(n+1)/2$. For even $n=2m$, the alternative $\sigma=-1$ gives
\begin{equation*}
	a_k=a_0\binom mk\Big(\frac um\Big)^k,\qquad P_{\rm deg}=a_0\,w\Big[xy+\Big(z^2+\frac umw\Big)^m\Big].
\end{equation*}
When $m\ge2$ its critical locus contains 
\begin{equation*}
    \left\{x=y=0,\ z^2+\frac umw=0\right\},
\end{equation*}
including points with $w\neq0$. The reduced geometry given in equation \eqref{eq:A1radical} therefore excludes this branch.

\subsection{Hypergeometric form and the explicit syzygy}
\label{sec:hypergeo}

We define
\begin{equation}\label{eq:Hn}
	r=\frac w{z^2},\qquad H_n(r)=\sum_{k=0}^mC_{n,k}\,r^k,\qquad C_{n,0}=1 .
\end{equation}
Then, $P_n=xyw+z^nw\,H_n(r)$, interpreted as a polynomial after expansion. The recurrence is equivalent to
\begin{equation}\label{eq:hypergeometric-ode}
	r(1-\kappa r)\,H_n''+\Big(2+\big(n-\tfrac32\big)\kappa r\Big)H_n'-2H_n=0,\qquad \kappa=\frac8{n(n-1)}.
\end{equation}
Recall the hypergeometric series and rising factorial \cite{AAR}:
\begin{equation*}
    {}_2F_1(a,b;c;t)=\sum_{k\ge0}\frac{(a)_k(b)_k}{(c)_k\,k!}t^k,
    \qquad (a)_0=1,\quad (a)_k=\prod_{j=0}^{k-1}(a+j).
\end{equation*}
Equation \eqref{eq:hypergeometric-ode} and $H_n(0)=1$ then give
\begin{equation}\label{eq:Hn-2F1}
	H_n(r)={}_2F_1\Big(-\frac n2,\ \frac{1-n}2;\ 2;\ \frac{8r}{n(n-1)}\Big).
\end{equation}
One of the first two parameters is a nonpositive integer, so the series terminates at degree $m$. An explicit additional syzygy is
\begin{equation}\label{eq:Sigman}
\begin{split}
	\Sigma_n={}&\Big(-x_1z_0+\frac n2x_0z_1\Big)F_{x,0}+\Big(-y_1z_0+\frac n2y_0z_1\Big)F_{y,0} -\frac{2(n+1)}{n(n-1)}\,w_1F_{z,0}\\
	&-n\,w_0z_1F_{w,0}+\frac4{n(n-1)}\,w_0F_{z,1}+z_0w_0F_{w,1}=0 .
\end{split}
\end{equation}
Substitution cancels the $x,y$-dependent terms, and the remaining expression is equation \eqref{eq:hypergeometric-ode}. Thus the hypergeometric differential equation is exactly the condition that $\Sigma_n$ vanish identically.

\subsection{The reduced critical locus}

\begin{proposition}\label{prop:A1crit}
For every $n\ge2$, $\sqrt{I_{P_n}}=(w,\ xy+z^n)$.
\end{proposition}

\begin{proof}
At $w=0$ the critical equations reduce to $xy+z^n=0$. Suppose $w\neq0$. The $x$ and $y$ derivatives force $x=y=0$. If $z=0$, the highest coefficient makes $\partial_zP_n$ or $\partial_wP_n$ nonzero, according to the parity of $n$. Thus $z\neq0$, and the remaining equations are
\begin{equation*}
	nH_n(r)-2rH_n'(r)=0,\qquad H_n(r)+rH_n'(r)=0 .
\end{equation*}
Since $r\neq0$, they imply 
\begin{equation*}
    H_n(r)=H_n'(r)=0 .
\end{equation*}
We show that this is impossible. Away from the singular points $0,\kappa^{-1}$ of equation \eqref{eq:hypergeometric-ode}, a solution and its derivative cannot vanish simultaneously unless the solution is identically zero. At $r=0$, $H_n(0)=1$. At $r=\kappa^{-1}$, the terminating Chu--Vandermonde identity \cite{AAR}
\begin{equation}\label{eq:vandermonde}
	{}_2F_1(-m,b;c;1)=\frac{(c-b)_m}{(c)_m},
\end{equation}
gives
\begin{equation*}
	H_{2m}(\kappa^{-1})=\frac{(m+\frac32)_m}{(2)_m}\neq0,\qquad H_{2m+1}(\kappa^{-1})=\frac{(m+\frac52)_m}{(2)_m}\neq0 .
\end{equation*}
Hence there is no critical point with $w\neq0$.
\end{proof}

\subsection{Minimality and the low rank cases}\label{sec:A1-min}

For $n\ge4$ the eight first jet Jacobian generators are minimal. A charged coefficient of the required degree cannot express one of $F_{x,j},F_{y,j}$ through the other generators. For the neutral generators, the only possible relation would express $F_{z,j}$ through degree-one coefficients multiplying $F_{w,0},F_{w,1}$, and the $xy$-terms force those coefficients to vanish; the intrinsic and jet degrees exclude the reverse dependence.

The coefficient of $F_{z,1}$ in $\Sigma_n$ contains $4w_0/[n(n-1)]$. Since $w_0$ is a coordinate, an inherited contribution to this monomial would require a syzygy at $(n+1,1)$ with a nonzero constant coefficient of $F_{z,1}$. The only other singlet terms there are $z_0F_{w,1}$ and $z_1F_{w,0}$, and their $x,y$-terms force both coefficients to vanish, leaving no identity. Hence, $[\Sigma_n]\neq0$ in $K/\mathfrak{m} K$. For the coefficients in equation \eqref{eq:Pn} the neutral solution space is spanned by equation \eqref{eq:A1-universal} and the vector of $\Sigma_n$, so the full singlet kernel is four-dimensional, and the three inherited identities in equation \eqref{eq:A1inherited} leave exactly one minimal class. For every other $P$ satisfying equation \eqref{eq:A1radical}, Section~\ref{sec:A1-solve} gives $\dim\ker M(a)=1$, so all singlet syzygies at $(n+3,1)$ are inherited. This completes the proof of Theorem~\ref{thm:A1}.

We now turn to the low rank cases. For $n=3$ the potential is
\begin{equation*}
	P_3=(xy+z^3)w+zw^2 .
\end{equation*}
At $(6,1)$, the ansatz by equation \eqref{eq:A1-ansatz} acquires three terms,
\begin{equation*}
    \lambda_{15}x_1y_0F_{w,0}+\lambda_{16}x_0y_1F_{w,0}+\lambda_{17}x_0y_0F_{w,1}.
\end{equation*}
Coefficient comparison leaves $\lambda_3,\lambda_{10},\lambda_{12},\lambda_{14},\lambda_{17}$ free. Put
\begin{equation*}
	\eta=\lambda_{14}-2\lambda_{12}-2\lambda_{17},\qquad \zeta=\lambda_{17},\qquad \Kos_w=F_{w,0}e_{w,1}-F_{w,1}e_{w,0}.
\end{equation*}
The complete solution is
\begin{equation*}
	\Sigma=\Big(\frac32\eta-\lambda_3\Big)z_1\Fl_0+\Big(\frac32\lambda_{12}-\lambda_{10}\Big)z_0\Fl_1+\lambda_{12}\,z_0\Eu+\eta\,\Sigma_3+\zeta\,\Kos_w .
\end{equation*}
The lower degree equations give $K^{U(1)}_{(5,0)}=\mathbb{C}\Fl_0$ and $K^{U(1)}_{(5,1)}=\mathbb{C}\Fl_1\oplus\mathbb{C}\Eu$. Since $\mathfrak{m}$ is generated by the coordinates, the inherited part at $(6,1)$ is obtained by multiplying lower syzygies by one coordinate. The degree and charge sectors needed for multiplication by $x_j,y_j,w_j$ have zero kernel, so
\begin{equation*}
	(\mathfrak{m} K)^{U(1)}_{(6,1)}=\mathrm{span}_\mathbb{C}\{z_1\Fl_0,\ z_0\Fl_1,\ z_0\Eu\}.
\end{equation*}
The inherited subspace is $\eta=\zeta=0$, and
\begin{equation*}
	\beta^{U(1)}_{2,(6,1)}\big(J_1(R_{P_3})\big)=2 .
\end{equation*}
The second class is the Koszul class $\Kos_w$, whose degree coincides with the distinguished one because $2n=n+3$ at $n=3$.

For $n=2$ one has $P_2=(xy+z^2)w+w^2$, and eliminating $w$ by $\partial_wP_2=0$ gives
\begin{equation*}
	R_{P_2}\simeq\frac{\mathbb{C}[x,y,z]}{\big(x(xy+z^2),\ y(xy+z^2),\ z(xy+z^2)\big)}.
\end{equation*}
The Jacobian generators are minimal as generators of an ideal in the original four-variable ring, even though that presentation has a removable algebra coordinate. Eliminating $w$ and its jet coordinates changes the polynomial ring over which the resolution is taken, so the uniform Betti calculation does not apply.

\begin{table}[H]
\centering
\renewcommand{\arraystretch}{1.25}
\begin{tabular}{@{\phantom{n}}c@{\phantom{nnn}}l@{\phantom{n}}}
	\toprule
	$n$ & $H_n(r)$\\\midrule
	4 & $1+r+\frac1{18}r^2$\\
	5 & $1+r+\frac1{10}r^2$\\
	6 & $1+r+\frac2{15}r^2+\frac1{675}r^3$\\
	7 & $1+r+\frac{10}{63}r^2+\frac5{1323}r^3$\\
	8 & $1+r+\frac5{28}r^2+\frac5{784}r^3+\frac1{43904}r^4$\\\bottomrule
\end{tabular}
\caption{The first hypergeometric polynomials $H_n$ of equation \eqref{eq:Hn}, so that $P_n=xyw+z^nwH_n(w/z^2)$. They agree with the potentials found in \cite{KLSnonpoint}.}
\label{tab:Hn}
\end{table}

\begin{remark}[The stress tensor generates the nilradical]\label{rem:stress}
The equation $\partial_wP_n=0$ reads 
\begin{equation*}
    xy+z^n=-\sum_{k\ge1}(k+1)C_{n,k}z^{n-2k}w^k ,
\end{equation*}
so $xy+z^n\in(w)$ in $R_{P_n}$. Thus, Proposition~\ref{prop:A1crit} gives $\sqrt0=(w)$: the stress tensor operator generates the nilradical. Consequently the weights given by equation \eqref{eq:A1-qT} define exactly the nilradical filtration of Section~\ref{sec:filtration}, in which each factor of $w$ lowers the filtration weight by one. The first few $H_n$ are listed in Table~\ref{tab:Hn}.
\end{remark}

\section{The type A potential}\label{sec:typeA}

Let $\mathfrak{g}=\mathfrak{sl}_N(\mathbb{C})$ with the pairing $\ip AB=\tr(AB)$, and assign degree one to its linear coordinates. The invariant ring is $\mathbb{C}[\mathfrak{g}]^{SL_N}=\mathbb{C}[C_2,\dots,C_N]$. For an invariant homogeneous $P$ of degree $p+1$, write 
\begin{equation*}
    G=\nabla P,\qquad H_M=H_P(M).
\end{equation*}
The first jet equations are the components of
\begin{equation*}
	F_0=G(M_0),\qquad F_1=H_{M_0}M_1,\qquad \deg M_0=(1,0),\quad\deg M_1=(1,1).
\end{equation*}

\begin{theorem}\label{thm:typeA}
Let $N\ge3$, $p\ge2$, and $\gcd(p,N)=1$, and define
\begin{equation}\label{eq:PA}
	P^A_{p,N}(M)=[s^{p+1}]\,\det(1-sM)^{p/N}.
\end{equation}
\begin{enumerate}[label=(\arabic*)]
\item $P^A_{p,N}$ is a nonzero invariant polynomial of degree $p+1$ and a twisted period. Its first jet equations satisfy
\begin{equation}\label{eq:Aidentity}
	\tr(M_0^2F_1)-p\,\tr(M_0M_1F_0)=0,
\end{equation}
which is equivalent to $2H_P\nabla C_3=p\,H_{C_3}\nabla P$.
\item This identity represents the unique minimal singlet at $(p+2,1)$:
\begin{equation}
	\beta^{SL_N}_{2,(p+2,1)}\big(J_1(R_{P^A_{p,N}})\big)=1,
\end{equation}
and $Q\mapsto[\Sigma_Q]$ is an isomorphism
\begin{equation}
	\TP_3(\mathfrak{sl}_N)=\mathbb{C} C_3\longrightarrow(K/\mathfrak{m} K)^{SL_N}_{(p+2,1)}.
\end{equation}
\item $\Crit(P^A_{p,N})_\red=\{M\in\mathfrak{g}: M^p=0\}$.
\item Let $P$ be a nonzero invariant homogeneous polynomial of degree $p+1$ whose reduced critical locus is contained in the nilpotent cone. Then
\begin{equation*}
	\beta^{SL_N}_{2,(p+2,1)}\big(J_1(R_P)\big)\neq 0
\end{equation*}
if and only if $P$ is proportional to $P^A_{p,N}$.
\end{enumerate}
\end{theorem}

\subsection{The invariant first jet condition}\label{sec:Acondition}

Each first jet equation has intrinsic degree $p$, so a syzygy at $(p+2,1)$ has quadratic coefficients: those multiplying $\varepsilon_0$ are bilinear in $M_0,M_1$, and those multiplying $\varepsilon_1$ are quadratic in $M_0$. For $N\ge3$ the invariant trilinear forms on the adjoint representation are spanned by $\tr(ABC)$ and $\tr(BAC)$, so every singlet candidate is a combination of
\begin{equation*}
	\mathcal{X}=\tr(M_0M_1\varepsilon_0),\qquad \mathcal{X}'=\tr(M_1M_0\varepsilon_0),\qquad \mathcal{Y}=\tr(M_0^2\varepsilon_1).
\end{equation*}
Invariance of $P$ gives $[M,G(M)]=0$, and hence the universal commutator identity
\begin{equation*}
	\mathrm{Com}=\mathcal{X}-\mathcal{X}'\in K,\qquad \tr\big([M_0,M_1]F_0\big)=0 .
\end{equation*}
An additional identity exists precisely when the images $\tr(M_0M_1F_0)$ and $\tr(M_0^2F_1)$ of $\mathcal{X}$ and $\mathcal{Y}$ are linearly dependent.

We determine the proportionality by evaluating at a semisimple matrix whose centralizer forces the gradient onto a line. Let
\begin{equation}\label{eq:Mt}
	M_t=\diag\big((N-1)t,-t,\dots,-t\big),\qquad t\neq0 .
\end{equation}
A matrix fixed by the centralizer of $M_t$ is scalar on each eigenspace, and tracelessness leaves only the line $\mathbb{C} M_t$. Hence, equivariance gives
\begin{equation*}
	G(M_t)=b(t)\,M_t ,
\end{equation*}
and $b(t)\neq0$ when $\Crit(P)$ lies in the nilpotent cone. Consequently, $\varphi(M)=\tr(M^2G(M))$ is nonzero, since $\varphi(M_t)=b(t)N(N-1)(N-2)t^3$. Setting $M_1=M_0$ gives $F_1=pF_0$, so the images of $\mathcal{Y}$ and $\mathcal{X}$ become $p\varphi$ and $\varphi$. Hence, every additional identity is a multiple of $\mathcal{Y}-p\mathcal{X}$, whose image is equation \eqref{eq:Aidentity}.

To compare with Proposition~\ref{prop:identity-TP}, note that $\tr(M_0^2F_1)=\ip{H_P\nabla C_3}{M_1}$ and, since $M$ and $G$ commute, 
\begin{equation*}
    \tr(M_0M_1F_0)=\ip{(MG)_0}{M_1}=\tfrac12\ip{H_{C_3}G}{M_1}
\end{equation*}
by equation \eqref{eq:gradQ}. Thus equation \eqref{eq:Aidentity} is equivalent to $2H_P\nabla C_3=pH_{C_3}\nabla P$, and as formal elements
\begin{equation}\label{eq:SigmaC3}
	\Sigma_{C_3}=2(\mathcal{Y}-p\mathcal{X})+p\,\mathrm{Com} .
\end{equation}

\subsection{Construction from the characteristic polynomial}\label{sec:Aconstruction}

Write
\begin{equation*}
	\Delta(s)=\det(1-sM)=1+\delta_2s^2+\dots+\delta_Ns^N,\qquad \Psi(s)=\Delta(s)^\alpha=\sum_{\ell\ge0}\psi_\ell s^\ell,
\end{equation*}
the fractional power being the formal series with $\psi_0=1$. Since $d\log\Delta=-\sum_{j\ge0}s^{j+1}\tr(M^j\,dM)$,
\begin{equation}\label{eq:Agradpsi}
	\nabla\psi_d=-\alpha\Big(\sum_{j=0}^{d-1}\psi_{d-1-j}M^j\Big)_0 .
\end{equation}
Set $B(t)=\sum_{j=0}^p\psi_{p-j}t^j$. The proof rests on two cancellations, both of which require $\alpha=p/N$. The first concerns the trace of $B(M)$. From $\tr(1-sM)^{-1}=N-s\Delta'/\Delta$ one obtains
\begin{equation*}
	\Psi(s)\,\tr(1-sM)^{-1}=N\Psi(s)-\frac s\alpha\Psi'(s),
\end{equation*}
whose coefficient of $s^p$ gives $\tr B(M)=(N-p/\alpha)\psi_p$ for any $\alpha\neq0$. At $\alpha=p/N$ the trace vanishes, the projection in equation \eqref{eq:Agradpsi} is unnecessary, and
\begin{equation}\label{eq:AgradP}
	\nabla\psi_{p+1}=-\alpha B(M),\qquad \nabla\psi_{p+2}=(M\nabla\psi_{p+1})_0 .
\end{equation}

The second cancellation concerns differentiation along the gradient of the cubic Casimir. Let $\partial_v$ denote differentiation along $v(M)=(M^2)_0=\nabla C_3$. Direct differentiation gives
\begin{equation*}
	\partial_v\Psi(s)=\Big(1-\frac{\tr(M^2)}Ns^2\Big)\Psi'(s)+\alpha\tr(M^2)\,s\,\Psi(s),
\end{equation*}
and thus,
\begin{equation*}
	\partial_v\psi_d=(d+1)\psi_{d+1}+\tr(M^2)\Big(\alpha-\frac{d-1}N\Big)\psi_{d-1}.
\end{equation*}
At $d=p+1$ and $\alpha=p/N$ the second term vanishes.

Now let $P=\psi_{p+1}$ with $\alpha=p/N$. Then, 
\begin{equation*}
    \varphi=\tr(M^2\nabla P)=\partial_vP=(p+2)\psi_{p+2} ,
\end{equation*}
so equation \eqref{eq:AgradP} gives 
\begin{equation*}
    \nabla\varphi=(p+2)(M\nabla P)_0 .
\end{equation*}
On the other hand, with $L_M(V)=\tfrac12(MV+VM)_0$,
\begin{equation*}
	\nabla\varphi=2L_M\nabla P+H_Pv .
\end{equation*}
Since $L_M\nabla P=(M\nabla P)_0$, we obtain $H_Pv=p\,L_M\nabla P$, that is, 
\begin{equation*}
    2H_P\nabla C_3=pH_{C_3}\nabla P.
\end{equation*}
Pairing with $M_1$ gives equation \eqref{eq:Aidentity}, and Proposition~\ref{prop:identity-TP} shows that $P^A_{p,N}$ is a twisted period. That $P^A_{p,N}\neq0$ follows from Proposition~\ref{prop:A-crit} below. This proves part~(1) of Theorem~\ref{thm:typeA}.

For computation, $\Delta\Psi'=\alpha\Delta'\Psi$ yields the recurrence
\begin{equation}\label{eq:Arecurrence}
	\psi_\ell=\frac1\ell\sum_{j=2}^{\min(N,\ell)}\big((\alpha+1)j-\ell\big)\,\delta_j\,\psi_{\ell-j},\qquad \psi_0=1,
\end{equation}
and, equivalently,
\begin{equation*}
	P^A_{p,N}=\sum_{\substack{m_2,\dots,m_N\ge0\\ \sum_jjm_j=p+1}}\frac{\alpha^{\underline m}}{m_2!\cdots m_N!}\,\delta_2^{m_2}\cdots\delta_N^{m_N},\qquad m=\sum_jm_j,
\end{equation*}
where $\alpha^{\underline m}=\alpha(\alpha-1)\cdots(\alpha-m+1)$ is the falling factorial, distinct from the rising factorial in equation \eqref{eq:vandermonde}.

\subsection{The reduced critical locus}

\begin{proposition}\label{prop:A-crit}
For coprime $p,N$, $\Crit(P^A_{p,N})_\red=\{M\in\mathfrak{sl}_N:M^p=0\}$.
\end{proposition}

\begin{proof}
Let $\chi(t)=\det(t-M)$ and expand at infinity:
\begin{equation*}
\chi(t)^{p/N}=t^p\sum_{k\ge0}\psi_kt^{-k}=B(t)+\Theta(t),\qquad \Theta(t)=O(t^{-1}).
\end{equation*}
By equation \eqref{eq:AgradP}, 
\begin{equation*}
    \nabla P^A_{p,N}(M)=-\frac pNB(M) .
\end{equation*}
Suppose $M$ is critical. Then $B(M)=0$, so every eigenvalue of $M$ is a root of $B$, whether or not $M$ is semisimple. The polynomial
\begin{equation*}
	L(t)=\chi(t)B'(t)-\frac pN\chi'(t)B(t)=-\chi\Theta'+\frac pN\chi'\Theta
\end{equation*}
has degree at most $N-2$. If $\lambda$ is a root of $\chi$ of multiplicity $m_\lambda$, both terms of $L$ are divisible by $(t-\lambda)^{m_\lambda}$, since $B(\lambda)=0$. Hence $\chi\mid L$, so $L=0$ and
\begin{equation*}
	\frac{B'}B=\frac pN\,\frac{\chi'}\chi,\qquad \operatorname{ord}_\lambda B=\frac pN\,m_\lambda\in\mathbb{Z} .
\end{equation*}
Coprimality forces $N\mid m_\lambda$, so all eigenvalues coincide, and tracelessness makes their common value zero. For nilpotent $M$, $\Delta=1$, so $\psi_k(M)=0$ for $k>0$, $B(M)=M^p$, and the critical equation is exactly $M^p=0$. In particular $P^A_{p,N}\neq0$, since otherwise every traceless matrix would be critical.
\end{proof}

Writing $N=ap+r$ with $0\le r<p$, the reduced critical locus is the closure of the nilpotent orbit with partition $[p^a,r]$, omitting a zero part \cite{CM}; for $p>N$ it is the full nilpotent cone.

\begin{remark}\label{rem:coprime}
Without coprimality, the integrality condition $\tfrac pNm_\lambda\in\mathbb{Z}$ can permit several eigenvalues, and the reduced support of the same determinant expression requires a separate analysis.
\end{remark}

\subsection{Uniqueness}

Let $P$ be as in Theorem~\ref{thm:typeA}(4), and suppose that its Betti number at $(p+2,1)$ is nonzero. The first jet presentation is minimal by the irreducibility argument at the beginning of Section~\ref{sec:Amin}, which applies to every nonzero invariant homogeneous $P$. The nonvanishing Betti number therefore gives a nonzero invariant class in $(K/\mathfrak mK)_{(p+2,1)}$. Since taking $SL_N$-invariants is exact, this class has an invariant representative $\Sigma\in K_{(p+2,1)}$.

Invariance of $P$ gives the differentiated commutator syzygy
\begin{equation*}
	C_1=[M_1,\varepsilon_0]+[M_0,\varepsilon_1],
\end{equation*}
whose components lie in $K$, and hence
\begin{equation*}
	\mathrm{Com}=\tr(M_0C_1)\in\mathfrak mK.
\end{equation*}
Since the class of $\Sigma$ in $K/\mathfrak mK$ is nonzero, $\Sigma$ is not a multiple of $\mathrm{Com}$. By Section~\ref{sec:Acondition}, such a syzygy forces the identity $\mathcal Y-p\mathcal X\in K$, equivalently,
\begin{equation*}
	2H_P\nabla C_3=pH_{C_3}\nabla P.
\end{equation*}
Proposition~\ref{prop:identity-TP} then gives $P\in\TP_{p+1}(\mathfrak{sl}_N)$. By Corollary~\ref{thm:FS}, this space is at most one-dimensional and contains the nonzero polynomial $P^A_{p,N}$. Therefore $P$ is proportional to $P^A_{p,N}$. Conversely, $P^A_{p,N}$ has Betti number one at $(p+2,1)$ by part~(2). This proves part~(4). The argument uses only an invariant syzygy at $(p+2,1)$ that is not a multiple of $\mathrm{Com}$, which is a weaker hypothesis than the nonvanishing of the Betti number.

The residue formula for twisted periods \cite{DubrovinAlmost,FS} is consistent with the determinant formula. Here,$[t^{-1}]$ denotes the coefficient of $t^{-1}$ in the expansion at infinity, so that 
\begin{equation*}
    [t^{-1}]A(t)=-\Res_{t=\infty}A(t)\,dt
\end{equation*}
for the conventional residue of a one-form at infinity. With the integrand expanded with leading term $t^p$,
\begin{equation}\label{eq:Aresidue}
\begin{split}
	P^A_{p,N}(M)&=-\Res_{t=\infty}\det(t-M)^{p/N}\,dt\\
	&=[t^{-1}]\prod_{i=1}^N(t-x_i)^{p/N}=[s^{p+1}]\prod_{i=1}^N(1-sx_i)^{p/N}
\end{split}
\end{equation}
on the Cartan subalgebra.

\subsection{Minimality of the singlet syzygy}\label{sec:Amin}

All first jet generators have intrinsic degree $p$, and the presentation is minimal. A constant linear relation among the components of $F_0$ would define a nonzero invariant subspace of the irreducible adjoint representation annihilated by $G$, forcing $G=0$. A constant relation among the components of $F_1=H_{M_0}M_1$ would make the corresponding combination of components of $G$ constant, hence zero by homogeneity. The two jet degrees cannot mix in a constant relation. Thus there are no syzygies of intrinsic degree $p$.

Since $\mathfrak{m}$ is generated by the linear coordinates, which form two copies of the adjoint representation, an inherited singlet at intrinsic degree $p+2$ is a contraction of an adjoint-valued syzygy of intrinsic degree $p+1$ with $M_0$ or $M_1$. We determine these syzygies.

At jet degree zero the candidates are $[M_0,\varepsilon_0]$ and $L_{M_0}\varepsilon_0$. The image of the first vanishes. The image of the second does not, since at equation \eqref{eq:Mt}
\begin{equation*}
	L_{M_t}G(M_t)=b(t)(N-2)\,t\,M_t\neq0 .
\end{equation*}

At jet degree one the candidates are $[M_1,\varepsilon_0]$, $L_{M_1}\varepsilon_0$, $[M_0,\varepsilon_1]$, $L_{M_0}\varepsilon_1$, with coefficients $a_1,a_2,a_3,a_4$. As an operator on $M_1$, their image is
\begin{equation*}
	-a_1\ad_G+a_3\ad_MH+a_2L_G+a_4L_MH,\qquad L_G(V)=\tfrac12(GV+VG)_0 .
\end{equation*}
Differentiating the equivariance of $G$ gives $H\ad_M=\ad_G$, and transposing, with $H$ self-adjoint and $\ad$ skew-adjoint, gives $\ad_MH=\ad_G$. By Lemma~\ref{lem:commute} applied to the identity of part~(1), $H$ commutes with $H_{C_3}=2L_M$, so $L_MH$ is self-adjoint, as is $L_G$. A sum of a self-adjoint and a skew-adjoint operator vanishes only if both parts do. The skew part gives $a_3=a_1$, which is the differentiated commutator identity
\begin{equation*}
	C_1=[M_1,\varepsilon_0]+[M_0,\varepsilon_1].
\end{equation*}
The self-adjoint part requires $a_2L_G+a_4L_MH=0$. Applying it to $M$ and using $HM=pG$ and $L_MG\neq0$ gives $a_2=-pa_4$, so a nonzero solution would require $L_MH=pL_G$. Differentiating equivariance at $M_t$ gives $H_{M_t}[A,M_t]=b(t)[A,M_t]$, and $A=E_{12}$ with $[E_{12},M_t]=-NtE_{12}$ gives $H_{M_t}E_{12}=b(t)E_{12}$. Since $L_{M_t}E_{12}=\tfrac12(N-2)tE_{12}$, the required identity would read
\begin{equation*}
	\tfrac12\,b(t)(N-2)\,t\,E_{12}=\tfrac p2\,b(t)(N-2)\,t\,E_{12},
\end{equation*}
which is impossible for $p\ge2$, $N\ge3$, and $t\neq0$.

The relevant lower syzygies are then $C_0=[M_0,\varepsilon_0]$ and $C_1$. Their singlet contractions at jet degree one are
\begin{equation*}
	\tr(M_1C_0)=-\mathrm{Com},\qquad \tr(M_0C_1)=\mathrm{Com},
\end{equation*}
which span the same line. The full singlet kernel at $(p+2,1)$ is two-dimensional, spanned by $\mathrm{Com}$ and $\mathcal{Y}-p\mathcal{X}$: the images of $\mathcal{X}$ and $\mathcal{X}'$ coincide, the image of $\mathcal{Y}$ is $p$ times that of $\mathcal{X}$ by part~(1), and the image of $\mathcal{X}$ is nonzero because it restricts to $\varphi\neq0$ at $M_1=M_0$. The quotient by the inherited line is one-dimensional, and by equation \eqref{eq:SigmaC3} it is spanned by $[\Sigma_{C_3}]=2[\mathcal{Y}-p\mathcal{X}]$. Since $\TP_3(\mathfrak{sl}_N)=\mathbb{C} C_3$ by Proposition~\ref{prop:basicTP}, this proves part~(2). Part~(3) is Proposition~\ref{prop:A-crit}.

\begin{remark}\label{rem:A-small}
For $N=2$ the symmetric trilinear form is absent and invariant theory alone fixes $P^A_{p,2}\propto C_2^{(p+1)/2}$ for odd $p$, so the cubic criterion is unnecessary. For $p=2$ the degree-three invariants are $\mathbb{C} C_3$, so $P^A_{2,N}\propto C_3=Q_{\mathfrak{sl}_N}$, part~(4) is automatic, and $\Sigma_{C_3}$ is proportional to the Koszul class $\Kos$.
\end{remark}

\subsection{Examples}

The characteristic coefficients begin
\begin{equation*}
	\delta_2=-C_2,\qquad \delta_3=-C_3,\qquad \delta_4=\tfrac12C_2^2-C_4,\qquad \delta_5=C_2C_3-C_5 .
\end{equation*}
Using equation \eqref{eq:Arecurrence} with a convenient overall normalization gives Table~\ref{tab:typeA}. These agree with \cite{KLSnonpoint}, and the formula applies uniformly, including when the reduced critical locus is a proper orbit closure.

\begin{table}[H]
\centering
\renewcommand{\arraystretch}{1.3}
\begin{tabular}{@{\phantom{T}}ll@{\phantom{T}}}
	\toprule
	Theory & Potential (up to normalization)\\\midrule
	$\mathcal{D}_5(SU(3))$ & $C_3^2+\frac19C_2^3$\\
	$\mathcal{D}_{11}(SU(3))$ & $C_3^4+\frac23C_2^3C_3^2+\frac2{135}C_2^6$\\
	$\mathcal{D}_7(SU(4))$ & $C_4^2-\frac34C_4C_2^2+\frac14C_3^2C_2+\frac{29}{192}C_2^4$\\
	$\mathcal{D}_6(SU(5))$ & $C_5C_2+C_4C_3-\frac{11}{10}C_3C_2^2$\\\bottomrule
\end{tabular}
\caption{The potentials $P^A_{p,N}$ for several $\mathcal{D}_p(SU(N))$ theories.}
\label{tab:typeA}
\end{table}

\section{The type D potential}\label{sec:typeD}

Let $N\ge4$, $\mathfrak{g}=\mathfrak{so}_{2N}(\mathbb{C})$, and
\begin{equation}\label{eq:Dhyp}
h=2N-2,\qquad p\ge2,\qquad \gcd(p,h)=1 ;
\end{equation}
in particular $p$ is odd. We use the defining $2N$-dimensional representation and the pairing $\ip AB=\tfrac12\tr(AB)$ of equation \eqref{eq:pairings}, so that on the split Cartan subalgebra it is the standard pairing on $\mathbb{C}^N$. All gradients and Hessians in this section and the next are taken with respect to it; rescaling the pairing rescales them consistently. With $C_k=\tfrac1k\tr(M^k)$ we have $\nabla C_2=2M$ and $\nabla C_4=2M^3$, and 
\begin{equation*}
    Q_4=\frac12C_4-\frac3{4h}C_2^2
\end{equation*}
is the fixed quartic of equation \eqref{eq:introQ}. Replacing $Q_4$ by another quartic changes the condition imposed on $P$, and Section~\ref{sec:quartic} shows that this choice cannot be dispensed with. The Pfaffian $\Pf$ is normalized by $\Pf|_\mathfrak{h}=x_1\cdots x_N$.

\begin{theorem}\label{thm:typeD}
Under equation \eqref{eq:Dhyp}, put $\alpha=p/h$ and
\begin{equation}\label{eq:PD}
	P^D_{p,N}(M)=[s^{p+1}]\,\det(1-sM)^\alpha .
\end{equation}
\begin{enumerate}[label=(\arabic*)]
\item $P^D_{p,N}$ is a nonzero invariant homogeneous polynomial of degree $p+1$ and a twisted period.
\item With $G=\nabla P^D_{p,N}$ and $H=H_{P^D_{p,N}}$,
\begin{equation}\label{eq:Didentity}
	3H\nabla Q_4-p\,H_{Q_4}G=0 .
\end{equation}
Consequently the first jet equations $F_0=G(M)$, $F_1=H_MY$ satisfy
\begin{equation}\label{eq:Djetidentity}
	3\ip{\nabla Q_4(M)}{F_1}-p\ip{H_{Q_4}(M)Y}{F_0}=0
\end{equation}
at bidegree $(p+3,1)$.
\item Over $S_1=\mathbb{C}[\mathfrak{g}\oplus\mathfrak{g}]$, with coordinate bidegrees $(1,0)$ and $(1,1)$,
\begin{equation}\label{eq:Dbetti}
	\beta^{SO(2N)}_{2,(p+3,1)}\big(J_1(R_{P^D_{p,N}})\big)=\begin{cases}2,&N=4,\\1,&N\ge5,\end{cases}
\end{equation}
and $Q\mapsto[\Sigma_Q]$ is an isomorphism from $\TP_4(\mathfrak{so}_{2N})$ onto this space. For $N\ge5$ it is spanned by the class of equation \eqref{eq:Djetidentity}. For $N=4$ the Pfaffian supplies a second class,
\begin{equation}\label{eq:DPfidentity}
	3\ip{\nabla\Pf(M)}{F_1}-p\ip{H_{\Pf}(M)Y}{F_0}=0,
\end{equation}
and the two classes form the standard doublet of the outer automorphism group $\Sym_3$ of $\mathfrak{so}_8$.
\item Every invariant homogeneous polynomial $P$ of degree $p+1$ satisfying equation \eqref{eq:Didentity}, with $G=\nabla P$ and $H=H_P$, is proportional to $P^D_{p,N}$.
\item $\Crit(P^D_{p,N})_\red=\{M\in\mathfrak{g}:M^p=0\}$.
\end{enumerate}
\end{theorem}

We first construct the potential and prove the quartic identity, then prove uniqueness and determine the critical locus. Section~\ref{sec:minimal} establishes part~(3).

\subsection{The generating function and the quartic identity}

Write
\begin{equation*}
	\Psi(s,M)=\det(1-sM)^\alpha=\sum_{d\ge0}\psi_d(M)\,s^d .
\end{equation*}
The characteristic polynomial is even, so $\psi_d=0$ for odd $d$. An even power of $M$ is orthogonal to $\mathfrak{so}_{2N}$, while an odd power lies in it, so differentiating the determinant gives, for even $d$,
\begin{equation}\label{eq:Dgradpsi}
	\nabla\psi_d=-2\alpha\sum_{\substack{1\le j\le d-1\\ j\ \mathrm{odd}}}\psi_{d-1-j}\,M^j,
\end{equation}
with no projection needed. Put $P=\psi_{p+1}$ and $G=\nabla P$. Then
\begin{equation}\label{eq:Dpsip3}
	\nabla\psi_{p+3}=M^2G-2\alpha PM .
\end{equation}
Let $\partial_{M^3}$ denote differentiation along the vector field $M\mapsto M^3$. Direct differentiation gives
\begin{equation*}
	\partial_{M^3}\Psi=\frac1s\,\partial_s\Psi+\alpha\tr(M^2)\,\Psi,\qquad \partial_{M^3}\psi_d=(d+2)\psi_{d+2}+2\alpha C_2\psi_d .
\end{equation*}
The constant term on the right of the first equation is $2\psi_2+\alpha\tr(M^2)=0$, since $\psi_2=-\alpha C_2$, so the second equation also holds at $d=0$. The derivative of $M\mapsto M^3$ in the direction $V$ is $M^2V+MVM+VM^2$. Since $G$ commutes with $M$, differentiating $\partial_{M^3}P$ and using equation \eqref{eq:Dpsip3} gives
\begin{equation*}
	H M^3+3M^2G=(p+3)\nabla\psi_{p+3}+2\alpha C_2G+4\alpha PM,
\end{equation*}
that is,
\begin{equation}\label{eq:Dmatrixidentity}
	HM^3=pM^2G+2\alpha C_2G-2\alpha(p+1)PM .
\end{equation}
Set $\lambda=3/h$. From equation \eqref{eq:introQ},
\begin{equation}\label{eq:DgradQ4}
	\nabla Q_4=M^3-\lambda C_2M,\qquad H_{Q_4}(V)=M^2V+MVM+VM^2-\lambda C_2V-2\lambda M\ip MV.
\end{equation}
Substituting equation \eqref{eq:Dmatrixidentity} and equation \eqref{eq:DgradQ4}, together with $HM=pG$ and $\ip MG=(p+1)P$, gives
\begin{equation}\label{eq:Dforced}
	3H\nabla Q_4-pH_{Q_4}G=(6\alpha-2p\lambda)\big(C_2G-(p+1)PM\big).
\end{equation}
The right side vanishes because $2p\lambda=6p/h=6\alpha$. This proves equation \eqref{eq:Didentity}, and pairing with $Y$ gives equation \eqref{eq:Djetidentity}.

\begin{remark}\label{rem:D-exponent}
Equations \eqref{eq:Dmatrixidentity} and \eqref{eq:Dforced} hold for every exponent $\alpha$. Whenever $C_2G-(p+1)PM$ is not the zero polynomial, the quartic identity forces $\alpha=p/h$. Section~\ref{sec:trace-singlets} proves the required nonvanishing for the coprime determinant potentials considered here. This qualification is necessary for arbitrary $\alpha$, since the extracted coefficient can vanish at other exponents.
\end{remark}

At the Cartan point $(x,0,\dots,0)$ the potential restricts to
\begin{equation*}
	(-1)^{(p+1)/2}\binom{\alpha}{(p+1)/2}x^{p+1}.
\end{equation*}
Coprimality makes $\alpha$ nonintegral, so this coefficient is nonzero. Invariance and homogeneity follow from the determinant formula, and Proposition~\ref{prop:identity-TP} shows that $P^D_{p,N}$ is a twisted period. This proves parts~(1) and~(2).

\subsection{Uniqueness}

Let $P$ be invariant and homogeneous of degree $p+1$ and satisfy equation \eqref{eq:Didentity}. By Proposition~\ref{prop:identity-TP}, $P\in\TP_{p+1}(\mathfrak{so}_{2N})$. By Corollary~\ref{thm:FS} this space is at most one-dimensional, and it contains $P^D_{p,N}\neq0$. This proves part~(4).

The type D residue formula for twisted periods \cite[Section~5]{FS}, related to the period formulas of \cite{EYY}, explains the exponent. With the series at infinity expanded with leading term $t^p$,
\begin{equation*}
	[t^{-1}]\,t^{-2\alpha}\prod_{i=1}^N(t^2-x_i^2)^\alpha=[t^{-1}]\,t^p\prod_{i=1}^N(1-t^{-2}x_i^2)^\alpha=[s^{p+1}]\prod_{i=1}^N(1-s^2x_i^2)^\alpha,
\end{equation*}
so that, with the conventional residue of a one-form at infinity,
\begin{equation}\label{eq:Dresidue2}
	P^D_{p,N}(M)=-\Res_{t=\infty}\,t^{-p/(N-1)}\det(t-M)^{p/(2N-2)}\,dt .
\end{equation}
The prefactor compensates for the difference between the dimension $2N$ of the defining representation and the Coxeter number $2N-2$.

\subsection{The reduced critical locus}

Fix $M$, put $\chi(t)=\det(t-M)$, and expand
\begin{equation*}
	t^{-2\alpha}\chi(t)^\alpha=t^p\sum_{d\ge0}\psi_d(M)\,t^{-d}=B(t)+\Theta(t),\qquad \Theta(t)=O(t^{-1}),
\end{equation*}
where $B(t)=t^p+\psi_2t^{p-2}+\dots+\psi_{p-1}t$ is odd and monic. Equation \eqref{eq:Dgradpsi} gives
\begin{equation}\label{eq:DGB}
	G(M)=-2\alpha B(M).
\end{equation}
The Laurent series $F=t^{-2\alpha}\chi^\alpha$ satisfies $t\chi F'=\alpha(t\chi'-2\chi)F$. Hence the polynomial
\begin{equation*}
	L=t\chi B'-\alpha(t\chi'-2\chi)B=-t\chi\Theta'+\alpha(t\chi'-2\chi)\Theta
\end{equation*}
has degree at most $2N-1$. If $M$ is critical, then $B(M)=0$, so every eigenvalue is a root of $B$. At a root $a$ of $\chi$ of multiplicity $m_a$, the polynomial $\chi'$ vanishes to order at least $m_a-1$ and $B(a)=0$, so $\chi'B$ vanishes to order at least $m_a$; at $a=0$ the extra factor $t$ compensates similarly. It follows that $\chi\mid L$, and the degree bound forces $L=0$:
\begin{equation*}
	\frac{B'}B=\alpha\,\frac{\chi'}\chi-\frac{2\alpha}t .
\end{equation*}
At a nonzero root $a$ of $\chi$, comparing residues gives $\operatorname{ord}_aB=\tfrac ph\operatorname{ord}_a\chi$, so coprimality forces $h\mid\operatorname{ord}_a\chi$. Since $\chi$ is even, $-a$ is a distinct root of the same multiplicity, which would require $2h=4N-4\le2N$, impossible for $N\ge4$. Therefore $M$ is nilpotent. For nilpotent $M$, $\chi(t)=t^{2N}$, so $B(t)=t^p$, and equation \eqref{eq:DGB} gives $G(M)=0$ exactly when $M^p=0$. This proves part~(5) as an equality of reduced varieties; it does not identify the Jacobian ideal with the ideal of entries of $M^p$.

\subsection{The nilpotent orbit closure}\label{sec:orbit}

At the boundary admissible level $k=-h+h/p$ \cite{KW88,KWb}, the associated variety of the simple affine vertex algebra is \cite{AvEM}
\begin{equation*}
	\mathbb O_p=\{M\in\mathfrak{g}:(\ad M)^{2p}=0\}.
\end{equation*}
For odd $p$ in type D this equals $\{M:M^p=0\}$. To see this, let $d$ be the largest Jordan block of a nilpotent $M$ in the defining representation, recall $\mathfrak{g}\simeq\Lambda^2\mathbb{C}^{2N}$, and decompose under the $\mathfrak{sl}_2$-triple through $M$ using the Clebsch--Gordan rule. If $d\le p$, every $\mathfrak{sl}_2$ highest weight in $\mathfrak{g}$ is at most $2p-2$, so $(\ad M)^{2p}=0$. Conversely, if $d\ge p+2$, the exterior square of that block contains highest weight $2d-4\ge2p$; if $d=p+1$, this even block occurs at least twice in an orthogonal partition, and the tensor product of two copies contains highest weight $2d-2=2p$. In both cases $(\ad M)^{2p}\neq0$.

Write $2N=ap+r$ with $0\le r<p$. The orbit partition is the largest orthogonal partition dominated by $[p^a,r]$, its $D$-collapse \cite{CM}. Since $p$ is odd, it is
\begin{equation*}
	[p^a,r]\ (r\text{ odd}),\qquad [p^a,r-1,1]\ (r>0\text{ even}),\qquad [p^a]\ (r=0),
\end{equation*}
with zero parts omitted. For $p\ge2N-1$ this is the principal partition $[2N-1,1]$.

\subsection{Trace formulas}

The odd traces vanish, and
\begin{equation*}
	\det(1-sM)^\alpha=\exp\Big(-\alpha\sum_{j\ge1}C_{2j}s^{2j}\Big).
\end{equation*}
Thus, up to a nonzero overall normalization,
\begin{align*}
	P^D_{3,N}&\propto C_4-\frac\alpha2C_2^2, & \alpha&=\frac3h,\\
	P^D_{5,N}&\propto C_6-\alpha C_2C_4+\frac{\alpha^2}6C_2^3, & \alpha&=\frac5h,\\
	P^D_{7,N}&\propto C_8-\alpha C_2C_6-\frac\alpha2C_4^2+\frac{\alpha^2}2C_2^2C_4-\frac{\alpha^3}{24}C_2^4, & \alpha&=\frac7h.
\end{align*}
Each expression is subject to coprimality, and traces above the fundamental degrees are reduced by the characteristic polynomial identities. In particular the $\mathcal{D}_3(SO(10))$ potential is $C_4-\frac3{16}C_2^2$. By Proposition~\ref{prop:basicTP}(iv), the first of these expressions is $Q_4$ up to normalization for every $N$, whether or not $\gcd(3,h)=1$.

\subsection{The role of the fixed quartic}\label{sec:quartic}

The uniqueness statement in Theorem~\ref{thm:typeD} uses the identity with the prescribed quartic $Q_4$. The following family shows why the reduced critical locus and the invariant Betti number alone do not suffice. On $\mathfrak{so}_{10}$ take $p=3$ and
\begin{equation}\label{eq:Plambda}
	P_\lambda=\frac12C_4-\frac\lambda4C_2^2,\qquad
	\lambda\in\mathbb{C}^\times\setminus\Big\{1,\frac12,\frac13,\frac14,\frac15\Big\}.
\end{equation}
At $\lambda=3/h=3/8$, this is the determinant potential up to normalization:
\begin{equation*}
	P_{3/8}=Q_4=-\frac43P^D_{3,5}.
\end{equation*}

\begin{proposition}[The numerical condition does not determine the Jacobian ideal]\label{prop:quartic-nonuniqueness}
For the family in equation \eqref{eq:Plambda}, the following hold.
\begin{enumerate}[label=(\roman*)]
	\item The reduced critical locus is independent of $\lambda$:
	\begin{equation*}
		\Crit(P_\lambda)_\red=\{M\in\mathfrak{so}_{10}:M^3=0\}.
	\end{equation*}
	\item The first jet scheme has
	\begin{equation*}
		\beta^{SO(10)}_{2,(6,1)}\big(J_1(R_{P_\lambda})\big)=1,
	\end{equation*}
	and its minimal invariant syzygy space at this bidegree is spanned by the Koszul class.
	\item The Jacobian ideals $I_{P_\lambda}\subset\mathbb{C}[\mathfrak{so}_{10}]$ are pairwise distinct.
	\item The identity
	\begin{equation*}
		3H_{P_\lambda}\nabla Q_4=3H_{Q_4}\nabla P_\lambda
	\end{equation*}
	holds if and only if $\lambda=3/8$.
\end{enumerate}
\end{proposition}

\begin{proof}
The gradient is $G_\lambda=M^3-\lambda C_2M$. If a critical matrix has $C_2\neq0$, its minimal polynomial divides $t(t^2-\lambda C_2)$, which has distinct roots. It is therefore semisimple, with nonzero eigenvalues $\pm a$, where $a^2=\lambda C_2$. If there are $j$ nonzero eigenvalue pairs, counted with multiplicity, then $1\le j\le5$ and
\begin{equation*}
	C_2=ja^2=j\lambda C_2,
\end{equation*}
contrary to the excluded values of $\lambda$. Thus every critical point has $C_2=0$ and satisfies $M^3=0$. Conversely, $M^3=0$ implies $C_2=0$ and $G_\lambda=0$, proving (i).

For (ii), the invariant calculation in Sections~\ref{sec:evals}--\ref{sec:inherited} applies to this family. Its matrix identity is equation \eqref{eq:Dmatrixidentity} with $p=3$ and $\alpha$ replaced by $\lambda$, and the Cartan evaluations become
\begin{equation*}
	b_j=1-j\lambda,\qquad
	\rho_j=-\frac{j\lambda}{1-j\lambda},\qquad
	\tau=-\frac{2\lambda}{1-5\lambda}.
\end{equation*}
The arguments there use this identity and the conditions
\begin{equation*}
	b_1b_2b_5\neq0,\qquad
	\rho_1\neq1,\qquad \rho_1\neq\rho_2,\qquad \tau\neq0,
\end{equation*}
all of which hold for the stated values of $\lambda$. Thus, the invariant kernel at $(6,1)$ has dimension three, and its inherited subspace is
\begin{equation*}
	\operatorname{span}\{\mathrm{Com},\,C_2\Eu\},
\end{equation*}
of dimension two. The Koszul syzygy has bidegree $(6,1)$ and coefficient $G_\lambda=M^3-\lambda C_2M$ against $\varepsilon_1$. The inherited identities have no $M^3$ component in that coefficient, so the Koszul class is nonzero and spans the minimal quotient.

For (iii), each Jacobian ideal is generated in degree three. Its degree-three part is a copy of the adjoint representation spanned by the gradient components. Equality of $I_{P_\lambda}$ and $I_{P_\mu}$ would identify these submodules, so Schur's lemma would make $G_\lambda$ and $G_\mu$ proportional. The covariants $M^3$ and $C_2M$ are linearly independent, and comparison of their coefficients gives $\lambda=\mu$.

Finally, direct differentiation gives
\begin{equation}\label{eq:fixed-quartic-residual}
	3H_{P_\lambda}\nabla Q_4-3H_{Q_4}\nabla P_\lambda
	=6\Big(\lambda-\frac38\Big)\big(C_2M^3-2C_4M\big).
\end{equation}
The vector polynomial on the right is nonzero: at the Cartan point $(1,2,0,0,0)$ its first component is $5-17=-12$. This proves (iv).
\end{proof}

At $\lambda=3/8$, the distinguished twisted period identity is proportional to the Koszul identity, since $P_{3/8}=Q_4$. For every other allowed $\lambda$, the Koszul class still survives, but the identity with the fixed $Q_4$ fails. Thus the numerical Betti condition does not impose the twisted period equation; the prescribed form of the syzygy is essential.

A change of invariant generators does not remove this distinction. Writing
\begin{equation*}
	C_4'=\frac12C_4-\frac a4C_2^2
\end{equation*}
rewrites the same potential as
\begin{equation*}
	P_\lambda=C_4'+\frac{a-\lambda}4C_2^2,
\end{equation*}
without changing its Jacobian ideal. Choosing $a=\lambda$ separately for each potential makes that potential its own quartic generator. The resulting identity
\begin{equation*}
	3H_{P_\lambda}\nabla C_4'-3H_{C_4'}\nabla P_\lambda=0
\end{equation*}
is then tautological. It imposes no condition on $\lambda$, whereas the identity with the prescribed $Q_4$ selects $\lambda=3/8$.

For higher $p$, the same fixed quartic pairs $P^D_{p,N}$ with the degree-four member of its determinant family, as in Proposition~\ref{prop:basicTP}(iv). The physical selection of this partner through vacuum null relations is discussed in \cite{Kang:2026ney}.

\section{Minimal syzygies and triality in type D}\label{sec:minimal}

The identities of Section~\ref{sec:typeD} must be checked in the minimal quotient. We determine all invariant identities at $(p+3,1)$ and then identify those obtained from lower degree syzygies. Write $M=M_0$, $Y=M_1$, and
\begin{equation*}
	S_1=\mathbb{C}[\mathfrak{g}\oplus\mathfrak{g}],\qquad F_0=G(M),\qquad F_1=H_MY .
\end{equation*}
Let $\varepsilon_0,\varepsilon_1$ denote the formal adjoint generators of the presentation module, with relation map
\begin{equation*}
	S_1\otimes(\mathfrak{g}\oplus\mathfrak{g})\longrightarrow S_1,\qquad (A,B)\longmapsto\ip A{F_0}+\ip B{F_1}.
\end{equation*}
Let $K$ be its kernel, with the summands placed at $(p,0)$ and $(p,1)$, and let $\mathfrak{m}=(M,Y)$ be the ideal generated by all linear coordinates. Singlets in this section are invariants of $SO(2N,\mathbb{C})$.

The presentation is minimal. A constant relation among the components of $G$ gives a kernel in the irreducible adjoint representation and forces $G=0$. A constant relation among the components of $HY$ makes a combination of components of $G$ constant, which is impossible by homogeneity. The two sets cannot mix because their jet degrees differ. Hence
\begin{equation*}
\beta^{SO(2N)}_{2,(p+3,1)}\big(J_1(R_P)\big)=\dim\,(K/\mathfrak{m} K)^{SO(2N)}_{(p+3,1)} .
\end{equation*}

\subsection{Cartan evaluations}\label{sec:evals}

The following evaluations distinguish possible identities and exclude lower degree syzygies. They are computed from the determinant formula before any assumption about a candidate relation. Set $\ell=(p-1)/2$, and let $x^{(j)}$ have $j$ coordinates equal to one and the rest zero. For $g=\nabla f$, where $f=P^D_{p,N}|_\mathfrak{h}$, the generating function gives
\begin{equation*}
	g_i(x^{(j)})=\begin{cases}b_j,&i\le j,\\0,&i>j,\end{cases}\qquad b_j=-2\alpha(-1)^\ell\binom{j\alpha-1}{\ell}.
\end{equation*}
At a zero-coordinate direction $i>j$,
\begin{equation}\label{eq:evals-rho}
	H_{ii}(x^{(j)})=\rho_jb_j,\qquad \rho_j=\frac{j\alpha}{j\alpha-\ell}.
\end{equation}
At $x^{(N)}=\mathbf 1$ all off-diagonal entries coincide,
\begin{equation*}
	H_{ik}(\mathbf 1)=\tau b_N\quad(i\neq k),\qquad \tau=\frac{\alpha(p-1)}{N\alpha-1},
\end{equation*}
and Euler homogeneity then gives
\begin{equation}\label{eq:evals-euler}
	H\mathbf 1=pb_N\mathbf 1,\qquad Hv=b_N(p-N\tau)\,v\quad\text{if }\textstyle\sum_iv_i=0 .
\end{equation}
These formulas follow by expanding $(1-s^2)^{j\alpha-1}$, $(1-s^2)^{j\alpha}$, and $(1-s^2)^{N\alpha-2}$. They also establish the nonvanishing statements we need:
\begin{equation*}
	b_1b_2b_N\neq0,\qquad \rho_1\neq1,\qquad \rho_1\neq\rho_2,\qquad \tau\neq0 .
\end{equation*}
Indeed $\alpha$ has reduced denominator $h$, and $1,2,N<h$, so $j\alpha$ is nonintegral for these $j$; hence the binomial coefficients and the displayed denominators do not vanish. Moreover $\rho_1-1=\ell/(\alpha-\ell)$ and $\rho_1-\rho_2=\alpha\ell/[(\alpha-\ell)(2\alpha-\ell)]$ are nonzero since $p\ge3$.

\subsection{The trace singlets}\label{sec:trace-singlets}

The trace contractions give six candidates at $(p+3,1)$:
\begin{equation*}
\begin{aligned}
	\sT_1&=\ip{M^3}{\varepsilon_1}, & \sT_2&=C_2\ip M{\varepsilon_1}, & \sT_3&=\ip{M^2Y+YM^2}{\varepsilon_0},\\
	\sT_4&=\ip{MYM}{\varepsilon_0}, & \sT_5&=C_2\ip Y{\varepsilon_0}, & \sT_6&=\ip MY\ip M{\varepsilon_0}.
\end{aligned}
\end{equation*}
For $N\ge5$ this list is complete. Each expression contracts four adjoint tensors, or eight defining representation indices, and orthogonal invariant theory expresses such contractions through pairings, giving the single and double traces above \cite{Weyl,GW}. A contraction with the alternating volume tensor requires $2N$ distinct indices and cannot occur when $2N>8$. At $N=4$ it contributes the Pfaffian terms, and the trace list is then the part even under an orthogonal transformation of determinant $-1$; we treat the odd part in Section~\ref{sec:Pf-class}.

Equivariance, homogeneity, and equation \eqref{eq:Dmatrixidentity} give three elements of $K$:
\begin{align*}
	\mathrm{Com}=\sT_3-2\sT_4,\qquad
	C_2\Eu=\sT_2-p\sT_5,\qquad
	\Sigma_{\rm tr}=\sT_1-p\sT_4-2\alpha\sT_5+2\alpha\sT_6 .
\end{align*}
They exhaust the trace-sector kernel. To see this, consider the images of $\sT_4,\sT_5,\sT_6$. After pairing with $Y$, their coefficient vectors are
\begin{equation}\label{eq:T456}
	M^2G,\qquad C_2G,\qquad \ip MG\,M .
\end{equation}
Suppose 
\begin{equation*}
    a_1M^2G+a_2C_2G+a_3\ip MGM=0 .
\end{equation*}
Evaluation at $x^{(1)}$ and $x^{(N)}$ gives $a_1+a_2+a_3=0$ and $a_1+Na_2+Na_3=0$, so $a_1=0$ and $a_3=-a_2$. Differentiating the second Cartan component of $a_2(s_2g-(x\cdot g)x)=0$ at $x^{(1)}$ in the $x_2$-direction and using equation \eqref{eq:evals-rho} gives $a_2b_1(\rho_1-1)=0$, so $a_2=a_3=0$. The three images are independent, and the trace kernel has dimension $6-3=3$, with basis $\mathrm{Com},\ C_2\Eu,\ \Sigma_{\rm tr}$. In particular, 
\begin{equation*}
    C_2G-(p+1)PM=C_2G-\ip MGM\neq 0 .
\end{equation*}

\subsection{Inherited trace singlets}\label{sec:inherited}

Because $\mathfrak{m}$ is generated by linear coordinates, every inherited singlet at $(p+3,1)$ comes from an adjoint-valued syzygy of intrinsic degree $p+2$, contracted with $M$ or $Y$. Taking invariant parts is exact on the finite-dimensional graded pieces, since the group is reductive.

At jet degree zero the complete trace candidate is
\begin{equation*}
	a_1(M^2\varepsilon_0+\varepsilon_0M^2)+a_2\,M\varepsilon_0M+a_3\,C_2\varepsilon_0+a_4\,M\ip M{\varepsilon_0}.
\end{equation*}
Since $MGM=M^2G$, its image is $(2a_1+a_2)M^2G+a_3C_2G+a_4\ip MGM$, so independence of equation \eqref{eq:T456} shows that the kernel is spanned by $[M,[M,\varepsilon_0]]$.

At jet degree one write the coefficient operator of $\varepsilon_1$ as
\begin{equation*}
	A_M(V)=a_1(M^2V+VM^2)+a_2\,MVM+a_3\,C_2V+a_4\,M\ip MV ,
\end{equation*}
and let the remaining candidates, with coefficients $r_1,\dots,r_6$, be
\begin{equation*}
\begin{aligned}
	\sU_1&=MY\varepsilon_0+\varepsilon_0YM, & \sU_2&=M\varepsilon_0Y+Y\varepsilon_0M, & \sU_3&=YM\varepsilon_0+\varepsilon_0MY,\\
	\sU_4&=\ip MY\varepsilon_0, & \sU_5&=\ip M{\varepsilon_0}Y, & \sU_6&=\ip Y{\varepsilon_0}M .
\end{aligned}
\end{equation*}
These exhaust the same fourth-order trace tensors. Set
\begin{equation*}
	\mu_1=2a_1+a_2,\qquad \mu_2=a_3,\qquad \mu_3=2(r_1+r_2+r_3),\qquad \mu_4=r_6+pa_4 .
\end{equation*}
Restricted to the Cartan subalgebra, the syzygy condition becomes the operator equation
\begin{equation*}
	\big(\mu_1\diag(x_i^2)+\mu_2s_2I\big)H+\mu_3\diag(x_ig_i)+r_4\,gx^T+r_5\,(x\cdot g)I+\mu_4\,xg^T=0,
\end{equation*}
where the term $pa_4$ uses $Hx=pg$. This gives necessary linear conditions on the coefficients; we return afterwards to the full Lie algebra.

Evaluating on zero-coordinate directions at $x^{(1)}$ and $x^{(2)}$ gives $\mu_2\rho_1+r_5=0$ and $\mu_2\rho_2+r_5=0$, so $\mu_2=r_5=0$. The first diagonal entry at $x^{(1)}$, and the off-diagonal and diagonal minus off-diagonal entries at $x^{(N)}$, give
\begin{equation*}
	p\mu_1+\mu_3+r_4+\mu_4=0,\qquad \tau\mu_1+r_4+\mu_4=0,\qquad (p-N\tau)\mu_1+\mu_3=0 .
\end{equation*}
These imply $(N-1)\tau\mu_1=0$, hence $\mu_1=\mu_3=0$ and $r_4+\mu_4=0$. The remaining equation is $r_4(gx^T-xg^T)=0$. The gradient is not everywhere parallel to $x$, since differentiating $s_2g-(x\cdot g)x$ at $x^{(1)}$ in a zero-coordinate direction gives the nonzero factor $b_1(\rho_1-1)$. Therefore
\begin{equation*}
	\mu_1=\mu_2=\mu_3=\mu_4=r_4=r_5=0 .
\end{equation*}

It remains to impose the operator identity on the root spaces, which the Cartan calculation does not see. Equivariance and self-adjointness give $H\ad_M=\ad_G=\ad_MH$. After the preceding substitutions, the full operator equation reduces to
\begin{equation*}
	(a_1-r_1)\,\ad_M\ad_G=0 .
\end{equation*}
At $x^{(1)}$ this operator is nonzero on the root space of $e_1-e_2$, since $b_1\neq0$; thus $r_1=a_1$. Then, the adjoint-valued kernel is spanned by
\begin{equation*}
\begin{gathered}
	Z_1=[M,[M,\varepsilon_1]+[Y,\varepsilon_0]],\qquad Z_2=[Y,[M,\varepsilon_0]],\\
	Z_3=M\big(\ip M{\varepsilon_1}-p\ip Y{\varepsilon_0}\big),
\end{gathered}
\end{equation*}
which come from the commutator and Euler identities. Their scalar contractions are
\begin{equation*}
\begin{gathered}
	\ip Y{[M,[M,\varepsilon_0]]}=\mathrm{Com},\qquad \ip M{Z_1}=0,\\
	\ip M{Z_2}=-\mathrm{Com},\qquad \ip M{Z_3}=C_2\Eu .
\end{gathered}
\end{equation*}
Hence, the inherited trace singlet subspace is exactly $\mathrm{span}\{\mathrm{Com},C_2\Eu\}$, leaving one minimal class. With $\lambda=3/h$, the distinguished identity is
\begin{equation}\label{eq:SigmaQ4decomp}
	\Sigma_{Q_4}=3\ip{\nabla Q_4}{\varepsilon_1}-p\ip{H_{Q_4}Y}{\varepsilon_0}=3\Sigma_{\rm tr}-p\mathrm{Com}-3\lambda\,C_2\Eu ,
\end{equation}
so $[\Sigma_{Q_4}]=3[\Sigma_{\rm tr}]\neq0$. For $N\ge5$ there are no further singlet candidates, and Theorem~\ref{thm:typeD}(3) follows for $N\ge5$, since $\TP_4(\mathfrak{so}_{2N})=\mathbb{C} Q_4$ by Proposition~\ref{prop:basicTP}.

\subsection{The Pfaffian class for \texorpdfstring{$\mathfrak{so}_8$}{so8}}\label{sec:Pf-class}

Let $N=4$. Conjugation by an orthogonal matrix of determinant $-1$ fixes $P$ and every trace expression but sends $\Pf$ to $-\Pf$. The singlet calculation thus splits into even and odd parts, and the preceding calculation gives one minimal class in the even part.

Since $\Pf\in\TP_4(\mathfrak{so}_8)$ by Proposition~\ref{prop:basicTP}(ii) and $P\in\TP_{p+1}(\mathfrak{so}_8)$, Lemma~\ref{lem:pairing} with $q=3$ gives
\begin{equation}\label{eq:Pf-identity}
	3H_P\nabla\Pf-p\,H_{\Pf}\nabla P=0 .
\end{equation}
There are two odd singlet candidates,
\begin{equation*}
	\sV_1=\ip{\nabla\Pf}{\varepsilon_1},\qquad \sV_2=\ip{H_{\Pf}Y}{\varepsilon_0}.
\end{equation*}
The list is complete because the only odd fourth-order invariant tensor is contraction with the eight-dimensional volume form, the polarization of $\Pf$. At $x=\mathbf 1$ the Cartan Hessian $K_{\Pf}$ of $x_1x_2x_3x_4$ has zero diagonal and unit off-diagonal entries, so $K_{\Pf}g=3b_4\mathbf 1\neq0$. The image of $\sV_2$ is therefore nonzero, and equation \eqref{eq:Pf-identity} shows that the odd kernel is one-dimensional, spanned by $\Sigma_{\Pf}=3\sV_1-p\sV_2$.

To check minimality, note that at intrinsic degree $p+2$ and jet degree zero the only odd adjoint candidate is $H_{\Pf}\varepsilon_0$, whose image is nonzero by the same evaluation. At jet degree one the two odd adjoint candidates give the operator equation
\begin{equation*}
	a_1\,K_{\Pf}H(Y)+a_2\,(dK_{\Pf})[Y]G=0 .
\end{equation*}
At $x=\mathbf 1$, $K_{\Pf}$ acts by $3$ on the radial direction and by $-1$ on $\sum_iY_i=0$, while $(dK_{\Pf})[Y]G=2b_4K_{\Pf}Y$. Using equation \eqref{eq:evals-euler}, the two subspaces give
\begin{equation*}
	pa_1+2a_2=0,\qquad (p-4\tau)a_1+2a_2=0 .
\end{equation*}
Since $\tau\neq0$, $a_1=a_2=0$. There is no lower odd adjoint syzygy from which an inherited odd singlet could arise, so the Pfaffian class is minimal. Combining the two parities proves equation \eqref{eq:Dbetti}.

The calculations above distinguish the dimension of the invariant kernel from the number of minimal generators. Their results are summarized in Table~\ref{tab:typeDsyzygies} below; the two sectors together give two minimal classes for $\mathfrak{so}_8$. Section~\ref{sec:triality} identifies the representation of the triality group carried by these two classes.
\begin{table}[H]
    \centering
    \begin{tabular}{c|cccc}
        \toprule
        Sector at bidegree $(p+3,1)$ & Candidate space & Syzygies & Inherited syzygies & Minimal classes \\\midrule
        Trace sector, $N\ge4$ & $6$ & $3$ & $2$ & $1$ \\
        Pfaffian sector, $N=4$ & $2$ & $1$ & $0$ & $1$ \\\bottomrule
    \end{tabular}
    \caption{Invariant syzygies of the first jet equations of $P^D_{p,N}$ at bidegree $(p+3,1)$. For each sector, the columns give the number of candidate invariant identities, the dimension of the space of invariant syzygies among them, the dimension of its inherited subspace, and the number of minimal classes. For $N=4$ the trace sector is the part even under $O(8)$, and the Pfaffian sector is the odd part.}
\label{tab:typeDsyzygies}
\end{table}

\subsection{Triality and the twisted period form of the minimal syzygies}\label{sec:triality}

The fundamental invariant degrees of $\mathfrak{so}_{2N}$ are $2,4,\dots,2N-2$ and $N$, the last belonging to the Pfaffian. For even $N$ the degree $N$ occurs twice, and only at $N=4$ does this repetition occur at the first primitive degree above the quadratic:
\begin{equation*}
	\mathbb{C}[\mathfrak{so}_8]^{SO(8)}=\mathbb{C}[C_2,Q_4,\Pf,C_6],\qquad \deg(C_2,Q_4,\Pf,C_6)=(2,4,4,6).
\end{equation*}
Triality permutes the three eight-dimensional representations $8_v,8_s,8_c$ of $\mathrm{Spin}(8)$ \cite{Esole:2017qeh,MW}, and its action on quartic invariants can be read off from their weights. On the Cartan put $s_2=\sum_ix_i^2$ and $s_4=\sum_ix_i^4$, and recall $\Pf|_\mathfrak{h}=x_1x_2x_3x_4$. The vector weights are $\{\pm e_i\}$, and the spinor weights are $\tfrac12(\pm e_1\pm e_2\pm e_3\pm e_4)$ with the two sign parities \cite{Esole:2017qeh}. Summing fourth powers gives
\begin{equation*}
\begin{aligned}
	\tr_{8_v}(M^4)&=2s_4,\\
	\tr_{8_s}(M^4)&=\tfrac32s_2^2-s_4+12\Pf,\\
	\tr_{8_c}(M^4)&=\tfrac32s_2^2-s_4-12\Pf .
\end{aligned}
\end{equation*}
Subtracting the common invariant $s_2^2$ defines
\begin{equation*}
	q_v=8Q_4,\qquad q_s=-4Q_4+12\Pf,\qquad q_c=-4Q_4-12\Pf,\qquad q_v+q_s+q_c=0 .
\end{equation*}
Triality permutes these three quartics, which span $\TP_4(\mathfrak{so}_8)=\mathrm{span}\{Q_4,\Pf\}$. Hence $\TP_4(\mathfrak{so}_8)$ is the standard irreducible doublet of $\Sym_3$, and $\mathbb{C}[\mathfrak{so}_8]^{SO(8)}_4=\mathbb{C} C_2^2\oplus\TP_4(\mathfrak{so}_8)$ is the decomposition into the trivial and standard isotypic components.\footnote{For the standard action by Dynkin diagram automorphisms, the subalgebras of $\mathfrak{so}_8$ fixed by a transposition and by a $3$-cycle are $\mathfrak{so}_7$ and $\mathfrak{g}_2$, respectively. In an elliptic fibration, the monodromy permuting the components of a Kodaira $I_0^*$ fiber is a subgroup of $\Sym_3$. Trivial monodromy, monodromy of order two, and monodromy of order three or six give the gauge algebras $\mathfrak{so}_8$, $\mathfrak{so}_7$, and $\mathfrak{g}_2$, respectively, corresponding to the $\mathrm{Spin}(8)$-, $\mathrm{Spin}(7)$-, and $G_2$-models \cite{Esole:2017qeh,Esole:2018mqb}.}

\begin{remark}\label{rem:not-quotient}
The assignment $Q\mapsto\Sigma_Q$ is defined on the subspace $\TP_4(\mathfrak{so}_8)$, not on the quotient of the quartic invariants by $\mathbb{C} C_2^2$. Indeed $\Sigma_{C_2^2}=12\,C_2\Eu+8p(\sT_5-\sT_6)$, and $\sT_5-\sT_6\notin K$ by the independence of equation \eqref{eq:T456}. The twisted period condition is what singles out the standard isotypic component.
\end{remark}

\begin{proposition}\label{prop:triality}
Let $p\ge2$ with $\gcd(p,6)=1$. The potential $P^D_{p,4}$ is invariant under triality, and the map
\begin{equation}\label{eq:triality-map}
	\TP_4(\mathfrak{so}_8)\longrightarrow(K/\mathfrak{m} K)^{SO(8)}_{(p+3,1)},\qquad Q\longmapsto[\Sigma_Q],
\end{equation}
is an isomorphism of $\Sym_3$-representations.
\end{proposition}

\begin{proof}
We first show that the map given by equation \eqref{eq:triality-map} is a linear isomorphism; this part does not use the triality invariance of $P$. The map is well defined by Lemma~\ref{lem:pairing}. Conjugation by an orthogonal matrix of determinant $-1$ fixes $P$, so it acts on $(K/\mathfrak{m}K)^{SO(8)}_{(p+3,1)}$ as an involution. Under this involution $[\Sigma_{Q_4}]$ is even, since $Q_4$ is a trace invariant, while $[\Sigma_{\Pf}]$ is odd, since $\Pf\mapsto-\Pf$. Both classes are nonzero by Sections~\ref{sec:inherited} and~\ref{sec:Pf-class}, and nonzero vectors in different eigenspaces of an involution are linearly independent. Since $Q_4$ and $\Pf$ span $\TP_4(\mathfrak{so}_8)$ by Proposition~\ref{prop:basicTP}(iii), the map is injective, and the target has dimension two by the same sections, so the map is an isomorphism.

Next, we show that triality fixes $P^D_{p,4}$. The outer automorphism group acts on $\mathbb{C}[\mathfrak{so}_8]^{SO(8)}$, since inner automorphisms act trivially on invariant polynomials, and each class has a representative preserving $\mathfrak{h}$ and acting on it by an orthogonal map that permutes the roots. The tensor $\mathcal{B}_x$ of equation \eqref{eq:B} is unchanged under $\beta\mapsto-\beta$, so it depends only on $\mathfrak{h}$, the pairing, and the root system. Such a map therefore preserves the twisted period equation \eqref{eq:TP}, and $\Sym_3$ acts on each $\TP_d(\mathfrak{so}_8)$. By Corollary~\ref{thm:FS}, $\TP_{p+1}(\mathfrak{so}_8)=\mathbb{C} P^D_{p,4}$ is one-dimensional, so $\Sym_3$ acts on it by a character, either the trivial character or the sign. A transposition is represented on the Cartan by $x_4\mapsto-x_4$, which exchanges the two spinor nodes of the Dynkin diagram, and it fixes the determinant polynomial because the latter depends only on the $x_i^2$. The character is therefore trivial, and the entire triality group fixes $P^D_{p,4}$.

Finally, the pairing is a fixed multiple of the Killing form and is preserved by triality. Since triality fixes $P$, it preserves the gradient and Hessian of $P$, the first jet ideal, the syzygy module $K$, and $\mathfrak{m} K$, and $\Sigma_{\phi^*Q}=\phi^*\Sigma_Q$ for every $\phi$ in the triality group. Hence the isomorphism given by equation \eqref{eq:triality-map} is $\Sym_3$-equivariant.
\end{proof}

The triality action gives a second reason for the survival of the Pfaffian class in the minimal quotient: the inherited subspace is triality-stable and cannot remove one direction of the irreducible doublet. There are two singlets of the connected flavor group, forming one triality doublet. If one takes invariants under $O(8)$, the Pfaffian class is odd and is removed. Throughout we take $SO(8)$-invariants and retain the outer automorphism action on the resulting space. The calculation uses the adjoint representation only, so it is independent of the global form of the flavor group.

\begin{remark}\label{rem:F4}
Since $P^D_{p,4}$ is triality-invariant, its Cartan restriction is invariant under $W(D_4)\rtimes\Sym_3=W(F_4)$, so it can be written in $W(F_4)$-invariant generators of degrees $2,6,8,12$ \cite{Humphreys}. Such a presentation makes triality invariance manifest, but the Jacobian ring, its syzygy module, and its minimal resolution remain those of $\mathbb{C}[\mathfrak{so}_8]$; writing $P$ in triality-invariant generators does not quotient them by triality. In particular the doublet $\TP_4(\mathfrak{so}_8)$ contains no triality invariant.
\end{remark}

The results of Sections~\ref{sec:Amin} and~\ref{sec:inherited}--\ref{sec:triality} combine into the following statement.

\begin{theorem}\label{thm:TPform}
Let $(\mathfrak{g},h)$ and $p$ be as in equations \eqref{eq:introgh}--\eqref{eq:introP}, $c=c_\mathfrak{g}$, and $G=SL_N$ or $SO(2N)$. For $P=P_{p,\mathfrak{g}}$ the map
\begin{equation}
	\TP_c(\mathfrak{g})\longrightarrow(K/\mathfrak{m} K)^G_{(p+c-1,1)},\qquad Q\longmapsto[\Sigma_Q],
\end{equation}
is an isomorphism. Consequently $\beta^G_{2,(p+c-1,1)}(J_1(R_{P_{p,\mathfrak{g}}}))=\dim\TP_c(\mathfrak{g})$, which is $2$ for $\mathfrak{so}_8$ and $1$ otherwise.
\end{theorem}

\begin{proof}
The map is well defined by Lemma~\ref{lem:pairing}. For $\mathfrak{sl}_N$ this is Theorem~\ref{thm:typeA}(2). For $\mathfrak{so}_{2N}$ with $N\ge5$, both sides are one-dimensional and $[\Sigma_{Q_4}]\neq0$ by Section~\ref{sec:inherited}. For $\mathfrak{so}_8$ it is Proposition~\ref{prop:triality}. The dimensions are given by Proposition~\ref{prop:basicTP}(iii).
\end{proof}

Theorem~\ref{thm:TPform} also explains the exceptional behavior of $\mathfrak{so}_8$. The distinguished bidegree is $(p+c-1,1)$, and the minimal singlet classes there correspond to the twisted periods of degree $c$. For $N\ge5$ these are spanned by $Q_4$. For $N=4$, the parameter $3/h=1/2$ is the case $q\equiv N-1\pmod{2N-2}$ in which \cite{FS} find a two-dimensional solution space, and the second solution is the Pfaffian.

\subsection{The Pfaffian partner in every rank}

\begin{proposition}\label{prop:Pf-all-N}
Let $N\ge4$ and $P=P^D_{p,N}$. Then
\begin{equation}
	\Sigma_{\Pf}=(N-1)\ip{\nabla\Pf}{\varepsilon_1}-p\ip{H_{\Pf}Y}{\varepsilon_0}
\end{equation}
is an $SO(2N)$-invariant, $O(2N)$-odd syzygy of bidegree $(p+N-1,1)$.
\end{proposition}

\begin{proof}
By Proposition~\ref{prop:basicTP}(ii), $\Pf\in\TP_N(\mathfrak{so}_{2N})$, and Lemma~\ref{lem:pairing} applies with $q=N-1$.
\end{proof}

Thus, the Pfaffian identity exists in every rank, and its bidegree coincides with the distinguished bidegree $(p+3,1)$ precisely when $N=4$. We do not address whether $\Sigma_{\Pf}$ is minimal for $N\ge5$.

\subsection{Explicit rank-four examples}\label{sec:rank4}

For comparison with the potentials recorded in the atlas \cite{atlas}, set
\begin{equation*}
	I_{2k}=\sum_{i=1}^4x_i^{2k}=\tfrac12\tr(M^{2k}),\qquad \Pf|_\mathfrak{h}=x_1x_2x_3x_4 ,
\end{equation*}
so that $I_2=C_2$ and $I_4=2C_4$. Define
\begin{align*}
	P_6&=25I_2^3-90I_2I_4+72I_6,\\
	P_8&=127I_2^4-468I_2^2I_4+108I_4^2+5184\Pf^2+288I_2I_6 .
\end{align*}
The determinant formulas give exactly
\begin{equation*}
	P^D_{5,4}=-\frac5{1296}P_6,\qquad P^D_{7,4}=\frac7{31104}P_8 ,
\end{equation*}
the second equality using Newton's identity
\begin{equation*}
	I_8=\tfrac43I_2I_6+\tfrac12I_4^2-I_2^2I_4+\tfrac16I_2^4-4\Pf^2 .
\end{equation*}
These are equalities of potentials up to scale, and therefore of their Jacobian ideals. Theorem~\ref{thm:typeD} supplies two minimal singlets at $(8,1)$ and $(10,1)$, respectively, at Schur degrees $9$ and $11$. The two atlas records cited in \cite{atlas} report a singlet kernel of dimension four and an inherited subspace of dimension two, in agreement with the minimal quotient computed here. The additional ranks $5360$ and $5362$ in the $\mathcal{D}_5(SO(8))$ record refer to a zero-Cartan-weight block evaluated at the determinant potential and at a specified rational probe. They are not dimensions of the singlet space, and this comparison by itself is not a proof of the generic rank. Our proof establishes the minimal classes independently of those computations. A Hilbert series truncated below Schur degree $9$ or $11$, respectively, cannot detect these classes.

\section{Slodowy slices and the nilradical filtration}\label{sec:slices}

We treat $\mathfrak{g}=\mathfrak{sl}_N$ and $\mathfrak{g}=\mathfrak{so}_{2N}$ together, writing $P_{p,\mathfrak{g}}=P^A_{p,N}$ or $P^D_{p,N}$ with the corresponding coprimality assumption. Restriction to a Slodowy slice gives further critical schemes and, in the physical interpretation, describes nilpotent Higgsing \cite{Chacaltana:2012zy,Beem:2014rza}. Higgsed $\mathcal{D}_p^b(G)$ theories of types B, C, and D, together with their holographic duals, are studied in \cite{Couzens:2023kyf}.

The ambient Jacobian rings carry a natural Poisson structure. Identify $\mathfrak{g}$ with $\mathfrak{g}^*$ by the invariant pairing and give $\mathbb{C}[\mathfrak{g}]$ its Lie--Poisson bracket. For invariant $P$, equivariance of the gradient makes $I_P$ stable under the adjoint group. The Hamiltonian derivations of linear functions are infinitesimal adjoint actions, so they preserve $I_P$, and since linear functions generate $\mathbb{C}[\mathfrak{g}]$ the Leibniz rule gives
\begin{equation}\label{eq:poisson}
	\{\mathbb{C}[\mathfrak{g}],I_P\}\subseteq I_P .
\end{equation}
Thus the bracket descends to $R_P$, and the nonreduced critical schemes of this paper are Poisson.

By the Jacobson--Morozov theorem \cite{CM}, a nilpotent $e$ lies in an $\mathfrak{sl}_2$-triple $(e,h_e,f_e)$; we write $h_e$ to distinguish it from the Coxeter number. The Slodowy slice \cite{Slodowy} is
\begin{equation*}
	\mathcal{S}_e=e+\ker(\ad f_e).
\end{equation*}
Choose an $h_e$-eigenbasis $X_1,\dots,X_d$ of $\ker(\ad f_e)$ and write
\begin{equation*}
	M(y)=e+\sum_{i=1}^dy_iX_i,\qquad [h_e,X_i]=-\kappa_iX_i .
\end{equation*}

\subsection{The grading and the restricted critical scheme}

The Kazhdan action \cite{GG} is
\begin{equation*}
	\varrho_t(x)=t^2\Ad(t^{-h_e})\,x .
\end{equation*}
It preserves $\mathcal{S}_e$, since $[h_e,e]=2e$, and
\begin{equation*}
	\varrho_t(M(y))=M\big(t^{\kappa_1+2}y_1,\dots,t^{\kappa_d+2}y_d\big).
\end{equation*}
For nonzero $e$, ordinary dilations do not preserve the affine slice through $e$, and the Kazhdan action supplies a positive grading. We divide its weights by two, so that the potential retains degree $p+1$, and call the result the half-Kazhdan degree:
\begin{equation*}
	q_i=\deg y_i=1+\frac{\kappa_i}2 .
\end{equation*}
The $\kappa_i$ are nonnegative, so the action contracts the slice to $e$. Define
\begin{equation*}
	P_e(y)=P_{p,\mathfrak{g}}(M(y)),\qquad R_e=\mathbb{C}[\mathcal{S}_e]/I_{P_e}.
\end{equation*}
Invariance and ordinary homogeneity of $P_{p,\mathfrak{g}}$ give
\begin{equation*}
\begin{aligned}
    P_e\big(t^{\kappa_1+2}y_1,\dots,t^{\kappa_d+2}y_d\big)
    &=P_{p,\mathfrak{g}}\big(t^2\Ad(t^{-h_e})M(y)\big)\\
    &=t^{2(p+1)}P_e(y),
\end{aligned}
\end{equation*}
so $P_e$ has half-Kazhdan degree $p+1$.

The next proposition identifies the ideals after restriction. The chain rule gives one inclusion, and transversality supplies the converse; we recall the transversality argument and then prove the radical statement using the smooth action map. See \cite{Slodowy,LNS} for the classical slice geometry.

\begin{proposition}\label{prop:slice}
Let $P$ be an adjoint-invariant polynomial on a complex semisimple Lie algebra $\mathfrak{g}$, and let $\iota:\mathcal{S}_e\hookrightarrow\mathfrak{g}$. Then
\begin{equation}\label{eq:slice-ideal}
	I_{P|_{\mathcal{S}_e}}=\iota^*I_P\cdot\mathbb{C}[\mathcal{S}_e],\qquad\text{so}\qquad \Crit(P|_{\mathcal{S}_e})=\Crit(P)\times_\mathfrak{g}\mathcal{S}_e
\end{equation}
as closed schemes. Moreover,
\begin{equation}\label{eq:slice-radical}
	\sqrt{I_{P|_{\mathcal{S}_e}}}=\iota^*\big(\sqrt{I_P}\big)\cdot\mathbb{C}[\mathcal{S}_e],
\end{equation}
so the scheme-theoretic intersection of the reduced critical variety with $\mathcal{S}_e$ is reduced.
\end{proposition}

\begin{proof}
First, at every $x\in\mathcal{S}_e$,
\begin{equation}\label{eq:transversal}
	\mathfrak{g}=[\mathfrak{g},x]+\ker(\ad f_e).
\end{equation}
Consider the family of linear maps $\Upsilon_x:\mathfrak{g}\oplus\ker(\ad f_e)\to\mathfrak{g}$, $(A,V)\mapsto[A,x]+V$. Failure of surjectivity is defined by maximal minors and is a closed condition. Let $\gamma_t=\Ad(t^{-h_e})$. It preserves $\ker(\ad f_e)$, and
\begin{equation*}
	\Upsilon_{\varrho_t(x)}\big(\gamma_tA,\,t^2\gamma_tV\big)=t^2\gamma_t\,\Upsilon_x(A,V).
\end{equation*}
For $t\neq0$ the changes of source and target are invertible, so the rank is unchanged. A nonempty failure locus would then be closed and invariant under the contracting action, and would contain $e$. But the $\mathfrak{sl}_2$-decomposition gives $\mathfrak{g}=[\mathfrak{g},e]\oplus\ker(\ad f_e)$, so there is no failure locus.

The inclusion $I_{P|_{\mathcal{S}_e}}\subseteq\iota^*I_P\cdot\mathbb{C}[\mathcal{S}_e]$ is the chain rule. For the converse, view $\Upsilon$ as a surjective map of trivial vector bundles over the slice; locally it admits a regular right inverse. Every constant ambient tangent vector $V$ can thus be written locally as $V=[A(x),x]+\sum_ia_i(x)X_i$ with regular coefficients. Invariance of $P$ gives $dP_x([A(x),x])=0$, so $dP_x(V)=\sum_ia_i(x)\,\partial_{y_i}(P|_{\mathcal{S}_e})$. Every restricted ambient derivative thus lies locally in the slice Jacobian ideal, and equality of the ideal sheaves gives equation \eqref{eq:slice-ideal}.

For the radical statement, let $G$ be the connected adjoint group. By equation \eqref{eq:transversal}, the action map $a:G\times\mathcal{S}_e\to\mathfrak{g}$, $(g,x)\mapsto\Ad(g)x$, has surjective differential everywhere and is smooth. If $\mathcal{Z}=\Crit(P)_\red$, its inverse image is $G\times(\mathcal{Z}\times_\mathfrak{g}\mathcal{S}_e)$, since $\mathcal{Z}$ is invariant. This scheme is reduced, being the base change of a smooth morphism to the reduced scheme $\mathcal{Z}$, and reducedness descends along the faithfully flat projection to $\mathcal{Z}\times_\mathfrak{g}\mathcal{S}_e$. Thus $\iota^*(\sqrt{I_P})\cdot\mathbb{C}[\mathcal{S}_e]$ is radical, and by equation \eqref{eq:slice-ideal} its vanishing set is the critical set of $P|_{\mathcal{S}_e}$, proving equation \eqref{eq:slice-radical}.
\end{proof}

For either matrix potential write $\mathbb O_p=\{M\in\mathfrak{g}:M^p=0\}$; for $\mathfrak{so}_{2N}$ this is the variety of Section~\ref{sec:orbit}. The proposition gives
\begin{equation}\label{eq:slice-HB}
	\Spec(R_e)_\red=\mathbb O_p\cap\mathcal{S}_e
\end{equation}
with the reduced scheme structure. For the interpretation as nilpotent Higgsing we require $e\in\mathbb O_p$, so that the Higgs vev lies in the original Higgs branch; the algebraic statement holds for every slice. Equation \eqref{eq:slice-HB} parallels the description of associated varieties of Drinfeld--Sokolov reductions as intersections with Slodowy slices \cite{ArakawaW}.

\begin{corollary}\label{cor:nilradical}
Restriction induces a surjection $\pi_e:R_P\to R_{P|_{\mathcal{S}_e}}$. If $\mathfrak{n}=\sqrt0\subset R_P$ and $\mathfrak{n}_e=\sqrt0\subset R_{P|_{\mathcal{S}_e}}$, then
\begin{equation}\label{eq:nil-restrict}
	\mathfrak{n}_e=\pi_e(\mathfrak{n}).
\end{equation}
In particular $\pi_e(\mathfrak{n}^k)=\mathfrak{n}_e^k$ for every $k\ge0$, and restricting generators of the ambient nilradical gives generators of the slice nilradical.
\end{corollary}

\begin{proof}
The surjection follows from equation \eqref{eq:slice-ideal}, since $\mathcal{S}_e$ is a closed affine subspace. Equation \eqref{eq:slice-radical}, followed by quotienting by $I_{P|_{\mathcal{S}_e}}$, gives equation \eqref{eq:nil-restrict}, and taking products gives the statement for every power.
\end{proof}

The critical scheme identity also determines the restriction at every jet order. Since the jet functor preserves fiber products \cite{EM}, for $0\le m\le\infty$,
\begin{equation}\label{eq:jets-fiber}
	J_m\big(R_{P|_{\mathcal{S}_e}}\big)\simeq J_m(R_P)\otimes_{J_m(\mathbb{C}[\mathfrak{g}])}J_m(\mathbb{C}[\mathcal{S}_e]):
\end{equation}
a jet into the fiber product is a pair of jets with the same image in $\mathfrak{g}$. For the affine slice one can equivalently substitute $M(t)=e+\sum_iy_i(t)X_i$ into every coefficient equation. This is an equality of coordinate rings; it does not assert that minimal resolutions are obtained by base change.

\subsection{The intrinsic filtration}\label{sec:filtration}

Following \cite{Kang:2026ney}, the additional filtration is defined by powers of the full nilradical. Let $R=\bigoplus_dR_d$ be one of the positively graded Jacobian rings above, and let $\mathfrak{n}=\sqrt0$. Set
\begin{equation}\label{eq:nil-filtration}
F_sR_d=\big(\mathfrak{n}^{\max\{\lceil d-s\rceil,0\}}\big)_d,\qquad \mathfrak{n}^0=R .
\end{equation}
The filtration indices lie in the degree lattice, $\mathbb{Z}$ or $\tfrac12\mathbb{Z}$ in our examples, and the ceiling makes the exponent an integer. A degree-$d$ element of $\mathfrak{n}^k$ lies in $F_{d-k}R_d$; its weight is lower if it also lies in a higher power.

The filtration is increasing and multiplicative. If $a=\max\{\lceil d_1-s\rceil,0\}$ and $b=\max\{\lceil d_2-t\rceil,0\}$, then $a+b\ge\max\{\lceil d_1+d_2-s-t\rceil,0\}$, so $F_sR_{d_1}\cdot F_tR_{d_2}\subseteq F_{s+t}R_{d_1+d_2}$. It is exhaustive because $F_dR_d=R_d$, and separated because the nilradical of a Noetherian ring is nilpotent. A homogeneous nilpotent element of degree $d$ lies in $F_{d-1}R_d$, and a product of $k$ such factors lies in filtration degree at most $d-k$.

The associated graded ring records the successive nilpotent layers:
\begin{equation*}
    \bigl(\operatorname{gr}_F R\bigr)_{d,d-k}\cong\bigl(\mathfrak n^k/\mathfrak n^{k+1}\bigr)_d .
\end{equation*}
Equivalently, it is $\operatorname{gr}_{\mathfrak n}R$, with a class of intrinsic degree $d$ and nilradical order $k$ assigned filtration degree $d-k$. For $X=\operatorname{Spec}R$, its underlying algebra is therefore the coordinate ring of the normal cone $C_{X_{\mathrm{red}}}X$ \cite[Tag~062Z, Definition~31.20.5]{stacks-project}. The filtration expresses the nilpotent structure of the critical scheme through these successive layers.

The filtration depends only on the graded ring and its nilradical: it is the smallest increasing multiplicative filtration satisfying $R_d\subseteq F_dR$ and $\mathfrak{n}_d\subseteq F_{d-1}R$. Indeed, a degree-$d$ element of $\mathfrak{n}^k$ is a sum of products of $k$ homogeneous nilpotent factors with homogeneous coefficients; each factor lowers the weight by one, so multiplicativity places the sum in $F_{d-k}R$, and these are exactly the inclusions in equation \eqref{eq:nil-filtration}. A further decrease in weight would change the filtration and require additional input, as discussed in \cite{Kang:2026ney}.

For a concrete presentation, let $u_a$ be homogeneous algebra generators of $R=S/I_P$ of degrees $d_a$, and let $f_\nu$ represent homogeneous generators of $\mathfrak{n}$ of degrees $e_\nu$. Introduce auxiliary variables $v_\nu$:
\begin{equation}\label{eq:aux-presentation}
\begin{gathered}
	R\simeq\frac{\mathbb{C}[u_a,v_\nu]}{\big(I_P,\ v_\nu-f_\nu(u)\big)},\\
	\wt_{q,T}(u_a)=(d_a,d_a),\qquad \wt_{q,T}(v_\nu)=(e_\nu,e_\nu-1).
\end{gathered}
\end{equation}
The second weights define a quotient filtration, which is equation \eqref{eq:nil-filtration}. A monomial with $k$ auxiliary factors has weight $d-k$ and represents an element of $\mathfrak{n}^k_d$; conversely every element of $\mathfrak{n}^k_d$ is a sum of products of at least $k$ nilradical generators with homogeneous coefficients, and these are exactly the representatives of weight at most $d-k$. This also shows independence of the chosen generating set.

If $\mathfrak{n}$ is generated by the fundamental invariants $C_\nu$, equation \eqref{eq:aux-presentation} reduces to the auxiliary Casimir presentation with weights $(d_\nu,d_\nu-1)$. This occurs in particular when the reduced critical locus is the full nilpotent cone, whose ideal is generated by the fundamental invariants \cite{Kostant}. For smaller orbit closures the defining ideal can also require covariants, as in the classes studied in \cite{JS}. Their images must then be included among the nilradical generators, so the Casimir variables alone do not in general define equation \eqref{eq:nil-filtration}; the type D examples in \cite{Kang:2026ney} illustrate this distinction. For the $(A_1,A_{2n-1})$ family $\mathfrak{n}=(w)$, and equation \eqref{eq:aux-presentation} reduces to the weights in equation \eqref{eq:A1-qT} (Remark~\ref{rem:stress}).

On a slice, Corollary~\ref{cor:nilradical} determines $\mathfrak{n}_e$ and all its powers by restriction. We decompose these ideals into homogeneous components for the half-Kazhdan grading and then apply equation \eqref{eq:nil-filtration}. The order matters: an ambient covariant can restrict to terms of different half-Kazhdan degrees, so its ambient degree is not automatically its slice weight.

\subsection{Jets, arcs, and the index proposal}\label{sec:arcs}

Use the presentation in equation \eqref{eq:aux-presentation} and expand both $u_a(t)$ and $v_\nu(t)$ in jets, with
\begin{equation*}
	\deg u_{a,j}=(d_a,\,j,\,d_a),\qquad \deg v_{\nu,j}=(e_\nu,\,j,\,e_\nu-1),\qquad j\ge0,
\end{equation*}
recording intrinsic degree, jet degree, and filtration weight. Quotienting by all coefficients of the Jacobian and graph equations gives the induced multiplicative filtration on $J_\infty(R)$; jet derivatives preserve filtration degree. Equivalently, it is the filtration generated by the jets of every element of $F_sR$ with weight at most $s$, so it is independent of the presentation. The resulting object is the associated graded of the filtered coordinate ring of the arc space,
\begin{equation*}
	\gr_FJ_\infty(R)=\bigoplus_sF_sJ_\infty(R)/F_{<s}J_\infty(R).
\end{equation*}
The order is to form the arc space and then take the associated graded; interchanging these operations would require a separate argument. Define
\begin{equation*}
	\HS_{q,q,T}\big(\gr_FJ_\infty(R)\big)=\sum_{d,j,s}\dim_\mathbb{C}\big(\gr_FJ_\infty(R)\big)_{d,j,s}\,q^{d+j}T^s .
\end{equation*}
The physical proposal \cite{KLSnonpoint,Kang:2026ney} is
\begin{equation}\label{eq:indexproposal}
	I_{\Mac}(q,T)=\HS_{q,q,T}\big(\gr_FJ_\infty(R)\big)
\end{equation}
for the appropriate ambient or slice ring. The present results specify the rings and their filtration; they do not prove this equality at all orders.

One limit follows directly from the construction. Every term has $s\le d$ and $j\ge0$. After $T=t/q$ its weight is $q^{d+j-s}t^s$, which survives as $q\to0$ only for $j=0$ and $s=d$. On that diagonal $(\gr_FR)_{d,d}=R_d/\mathfrak{n}_d=(R/\mathfrak{n})_d$, so
\begin{equation}\label{eq:HL}
	\lim_{q\to0}\HS_{q,q,T}\big(\gr_FJ_\infty(R)\big)\Big|_{T=t/q}=\HS_t(R/\mathfrak{n}).
\end{equation}
This is the Hilbert series of the reduced ring, that is, of the Higgs branch for the rings constructed here. The equality is algebraic; its interpretation as the Hall--Littlewood index uses the physical identification in \cite{Kang:2022zsl,Kang:2026ney}.

\section{Discussion}\label{sec:discussion}

The closed-form formulas determine Jacobian rings supported on Kleinian surfaces and nilpotent orbit closures, and they characterize these rings through relations among their first jet equations. The reduced varieties are familiar; the results specify the nilpotent scheme structures and the syzygies that select them. For the matrix families the Jacobian ideals are Poisson ideals by equation \eqref{eq:poisson}, so the nonreduced critical schemes carry the ambient Lie--Poisson structure. These are the rings entering the protected operator construction, where passing to the reduced ring loses part of the proposed operator data.

\subsection{Critical schemes and first jet syzygies}

Theorems~\ref{thm:A1} and~\ref{thm:typeA} give uniqueness within the stated homogeneous invariant families: the prescribed reduced critical locus and one distinguished minimal first jet syzygy determine every coefficient of the potential, up to the indicated rescalings. The degree of the potential can be arbitrarily large, while the condition involves only first jets. The proofs explain why this is possible. The first jet equations contain the Hessian, so an invariant syzygy imposes a differential equation on the potential: in the four-variable family it is the hypergeometric equation \eqref{eq:hypergeometric-ode}, and in the matrix families it is the twisted period equation \eqref{eq:TP}.

We emphasize the hypotheses needed for uniqueness. For the $(A_1,A_{2n-1})$ and type A families, the numerical condition and the reduced geometry suffice. In type D, the uniqueness theorem uses the identity with the fixed quartic $Q_4$, and the family in equation \eqref{eq:Plambda} shows that this cannot be weakened: it has a fixed reduced critical locus and a fixed distinguished singlet Betti number, while its Jacobian ideals vary with $\lambda$. Thus the degree and multiplicity of the distinguished syzygy must be supplemented by its form. By Proposition~\ref{prop:basicTP}(iv), the partner $Q_4$ is the degree-four member of the same determinant family, so this form is determined by the family itself.

The minimality arguments compute actual classes in $\Tor^{S_1}_2(J_1(R_P),\mathbb{C})\simeq K/\mathfrak{m} K$: the invariant identities counted here survive after all polynomial multiples of lower degree syzygies are removed. They give explicit equivariant pieces of the minimal resolutions of first jet schemes of nonreduced critical loci. The full resolutions and their higher-jet analogues would give further information; the single Betti number computed here does not determine the full Hilbert series.

A related question arises already at jet order zero. If $\mathfrak{n}$ is the nilradical, the quotients $\mathfrak{n}^k/\mathfrak{n}^{k+1}$ are modules over the reduced ring $R_P/\mathfrak{n}$ and describe the successive nilpotent thickenings of the reduced critical locus. Determining these modules, with their symmetry representations, would describe the scheme structures selected by the first jet condition in more detail, and by Section~\ref{sec:filtration} it would also determine the graded pieces of the filtration entering the Macdonald refinement.

\subsection{Twisted periods}

The Cartan restrictions of the matrix potentials are the classical polynomial twisted periods \cite{DubrovinAlmost,EYY,FS}, and with the conventional residue at infinity
\begin{align*}
	P^A_{p,N}(M)&=-\Res_{t=\infty}\det(t-M)^{p/N}\,dt,\\
	P^D_{p,N}(M)&=-\Res_{t=\infty}\,t^{-p/(N-1)}\det(t-M)^{p/(2N-2)}\,dt,
\end{align*}
both integrands being expanded with leading term $t^p$. The residue expressions themselves belong to the earlier theory. Here,we characterize their invariant extensions by first jet identities on the full Lie algebra, determine their reduced critical loci, and prove which distinguished identities survive in the minimal resolution. Lemma~\ref{lem:pairing} gives the connection between these questions and twisted periods. For a twisted period $P$ of degree $p+1$ it produces a linear map
\begin{equation*}
	\Sigma:\ \bigoplus_d\TP_d(\mathfrak{g})\longrightarrow\bigoplus_dK^G_{(p+d-1,1)},\qquad Q\longmapsto\Sigma_Q ,
\end{equation*}
whose image contains the Euler identity ($Q=C_2$), the Koszul identity ($Q=P$), and the distinguished identity ($Q=Q_\mathfrak{g}$). Theorem~\ref{thm:TPform} shows that at $d=c_\mathfrak{g}$ this map is an isomorphism onto the minimal singlet classes. Two questions follow. Which of the remaining classes $[\Sigma_Q]$ are minimal, in particular for the Pfaffian partner of Proposition~\ref{prop:Pf-all-N} when $N\ge5$ and for the determinant partners $P_{q,\mathfrak{g}}$ with $c_\mathfrak{g}-1<q<p$? And does every invariant first jet syzygy of $P_{p,\mathfrak{g}}$ at jet degree one arise from $\Sigma$ modulo inherited ones?

For $\mathfrak{so}_8$, Proposition~\ref{prop:triality} identifies the primitive quartic doublet with the two-dimensional space of minimal singlets, and triality acts on this part of the minimal resolution as the standard representation of $\Sym_3$. The potential itself is unique for the fixed identity in the coprime range, while the syzygy multiplicity is the multiplicity of twisted periods of degree four. Both are governed by the classification of \cite{FS}, at the parameters $p/h$ and $3/h$ respectively, and only the second meets the exceptional case.

\subsection{Slodowy slices and the nilradical filtration}

The slice result applies to every adjoint-invariant polynomial on a complex semisimple Lie algebra. The transversality is classical \cite{Slodowy,LNS}; Proposition~\ref{prop:slice} uses it to identify the restricted Jacobian ideal and its radical. Differentiating after restriction produces the same ideal as restricting the ambient derivatives, and intersecting the reduced critical scheme with the slice creates no new nilpotent elements. Corollary~\ref{cor:nilradical} further gives $\pi_e(\mathfrak{n}^k)=\mathfrak{n}_e^k$, so the nilpotent structure and all its powers are obtained by restriction. The filtration used in the physical construction is intrinsic once the grading is specified: it is defined from these powers, independently of any choice of invariant or covariant generators. On a slice the appropriate grading is the half-Kazhdan grading, and the ambient degrees of individual covariants does not need to remain their degrees after restriction.

Jets preserve the fiber product, as recorded in equation \eqref{eq:jets-fiber}, so the same construction supplies equations at every finite jet order and on the arc space. This does not imply that a minimal ambient resolution restricts to a minimal resolution on the slice. Computing the resulting changes in syzygies would relate the commutative algebra of these critical schemes to the transverse geometry of the nilpotent orbits; in the physical interpretation the same operation describes nilpotent Higgsing and changes both the grading and the residual flavor symmetry.

\subsection{Protected operator algebras}

The physical origin remains essential to the choice of the algebraic conditions: the symmetry representations and degrees of the generators and the distinguished syzygies come from protected operator data. The formulas proved here determine the resulting commutative rings for entire families, including theories without a known Lagrangian description, and their restrictions to Slodowy slices give the rings proposed for the corresponding Higgsed theories. The exceptional constructions and index comparisons are developed in \cite{Kang:2026ney}.

To identify a Jacobian ring with the relevant vertex algebra quotient, one must match its ideal with the commutative images of the quantum null relations. Let $\mathcal V$ be the vertex operator algebra and $R_{\mathcal V}=\mathcal V/C_2(\mathcal V)$ its $C_2$-algebra \cite{Zhu,ArakawaC2}, where $C_2(\mathcal V)$ is the span of $a_{(-2)}b$ for $a,b\in\mathcal V$. Li's filtration \cite{Li} gives a canonical surjection of differential algebras \cite{AL}
\begin{equation*}
	J_\infty(R_{\mathcal V})\twoheadrightarrow\gr^{\rm Li}\mathcal V,
\end{equation*}
whose kernel can contain further differential relations. Thus even an isomorphism $R_P\simeq R_{\mathcal V}$ leaves a question about the full Schur character. The Macdonald refinement additionally requires the nilradical filtration to reproduce the filtration supplied by the four-dimensional charges \cite{Song:2016yfd}. The Hall--Littlewood specialization, as in equation \eqref{eq:HL}, follows algebraically from the nilradical filtration and recovers the reduced ring, while the Schur specialization forgets the filtration; neither fixes its intermediate weights. The construction of Section~\ref{sec:filtration} makes a precise choice: every homogeneous nilpotent class lies at least one step below its intrinsic degree, and products receive the further decreases forced by multiplicativity. Establishing the full filtered correspondence would explain why this intrinsic choice gives the protected operator count. Finally, in view of Theorem~\ref{thm:TPform}, it would be interesting to identify the vacuum null relation responsible for syzygy maximization directly with the pairing between $P_{p,\mathfrak{g}}$ and the lowest member $P_{c_\mathfrak{g}-1,\mathfrak{g}}$ of its family.

\section*{Acknowledgements}

The author is thankful for discussions and collaborations with Craig Lawrie and Jaewon Song. M.~J.~K.\ is supported by the DOE through DE-SC0010813 and the Start-up Research Grant for new faculty provided by Texas A\&M University.

\bibliography{refs}{}

\providecommand{\href}[2]{#2}\begingroup\raggedright\begin{thebibliography}{10}

\bibitem{AAR}
G.~E. Andrews, R.~Askey, and R.~Roy, {\em {Special Functions}}, vol.~71 of {\em Encyclopedia of Mathematics and its Applications}.
\newblock Cambridge University Press, Cambridge, 1999.

\bibitem{SUPERPOINT}
D.~Annamalai, H.~Hur, T.~Hussel, M.~J. Kang, C.~Lawrie, and J.~Song, ``{The Index from a Super-Point},'' {\it in preparation}, 2026.

\bibitem{ArakawaC2}
T.~Arakawa, ``{A remark on the $C_2$-cofiniteness condition on vertex algebras},'' {\em Math. Z.} {\bfseries 270} (2012) 559--575, \href{http://arxiv.org/abs/1004.1492}{{\ttfamily arXiv:1004.1492 [math.QA]}}.

\bibitem{ArakawaW}
T.~Arakawa, ``{Associated varieties of modules over Kac--Moody algebras and $C_2$-cofiniteness of $W$-algebras},'' {\em Int. Math. Res. Not. IMRN} {\bfseries 2015} no.~22, (2015) 11605--11666, \href{http://arxiv.org/abs/1004.1554}{{\ttfamily arXiv:1004.1554 [math.QA]}}.

\bibitem{AL}
T.~Arakawa and A.~R. Linshaw, \href{http://dx.doi.org/10.1007/978-3-030-23531-4_1}{``{Singular support of a vertex algebra and the arc space of its associated scheme},''} in {\em {Representations and Nilpotent Orbits of Lie Algebraic Systems}}, vol.~330 of {\em Progress in Mathematics}, pp.~1--17.
\newblock Birkh{\"a}user/Springer, Cham, 2019.
\newblock \href{http://arxiv.org/abs/1804.01287}{{\ttfamily arXiv:1804.01287 [math.RT]}}.

\bibitem{AvEM}
T.~Arakawa, J.~van Ekeren, and A.~Moreau, ``{On the nilpotent orbits arising from admissible affine vertex algebras},'' \href{http://dx.doi.org/10.1112/plms.12467}{{\em Proc. Lond. Math. Soc. (3)} {\bfseries 125} no.~3, (2022) 681--699}, \href{http://arxiv.org/abs/2010.08429}{{\ttfamily arXiv:2010.08429 [math.RT]}}.

\bibitem{Argyres:1995jj}
P.~C. Argyres and M.~R. Douglas, ``{New phenomena in SU(3) supersymmetric gauge theory},'' \href{http://dx.doi.org/10.1016/0550-3213(95)00281-V}{{\em Nucl. Phys. B} {\bfseries 448} (1995) 93--126}, \href{http://arxiv.org/abs/hep-th/9505062}{{\ttfamily arXiv:hep-th/9505062}}.

\bibitem{Argyres:1995xn}
P.~C. Argyres, M.~R. Plesser, N.~Seiberg, and E.~Witten, ``{New N=2 superconformal field theories in four-dimensions},'' \href{http://dx.doi.org/10.1016/0550-3213(95)00671-0}{{\em Nucl. Phys. B} {\bfseries 461} (1996) 71--84}, \href{http://arxiv.org/abs/hep-th/9511154}{{\ttfamily arXiv:hep-th/9511154}}.

\bibitem{AGZV}
V.~I. Arnold, S.~M. Gusein-Zade, and A.~N. Varchenko, {\em {Singularities of Differentiable Maps, Vol.~I}}, vol.~82 of {\em Monographs in Mathematics}.
\newblock Birkh{\"a}user, Boston, 1985.

\bibitem{Beem:2013sza}
C.~Beem, M.~Lemos, P.~Liendo, W.~Peelaers, L.~Rastelli, and B.~C. van Rees, ``{Infinite Chiral Symmetry in Four Dimensions},'' \href{http://dx.doi.org/10.1007/s00220-014-2272-x}{{\em Commun. Math. Phys.} {\bfseries 336} no.~3, (2015) 1359--1433}, \href{http://arxiv.org/abs/1312.5344}{{\ttfamily arXiv:1312.5344 [hep-th]}}.

\bibitem{Beem:2014rza}
C.~Beem, W.~Peelaers, L.~Rastelli, and B.~C. van Rees, ``{Chiral algebras of class S},'' \href{http://dx.doi.org/10.1007/JHEP05(2015)020}{{\em JHEP} {\bfseries 05} (2015) 020}, \href{http://arxiv.org/abs/1408.6522}{{\ttfamily arXiv:1408.6522 [hep-th]}}.

\bibitem{Beem:2017ooy}
C.~Beem and L.~Rastelli, ``{Vertex operator algebras, Higgs branches, and modular differential equations},'' \href{http://dx.doi.org/10.1007/JHEP08(2018)114}{{\em JHEP} {\bfseries 08} (2018) 114}, \href{http://arxiv.org/abs/1707.07679}{{\ttfamily arXiv:1707.07679 [hep-th]}}.

\bibitem{BMS}
C.~Bruschek, H.~Mourtada, and J.~Schepers, ``{Arc spaces and the Rogers--Ramanujan identities},'' \href{http://dx.doi.org/10.1007/s11139-012-9401-y}{{\em Ramanujan J.} {\bfseries 30} no.~1, (2013) 9--38}, \href{http://arxiv.org/abs/1101.4950}{{\ttfamily arXiv:1101.4950 [math.AG]}}.

\bibitem{Buican:2015ina}
M.~Buican and T.~Nishinaka, ``{On the superconformal index of Argyres{\textendash}Douglas theories},'' \href{http://dx.doi.org/10.1088/1751-8113/49/1/015401}{{\em J. Phys. A} {\bfseries 49} no.~1, (2016) 015401}, \href{http://arxiv.org/abs/1505.05884}{{\ttfamily arXiv:1505.05884 [hep-th]}}.

\bibitem{Cecotti:2012jx}
S.~Cecotti and M.~Del~Zotto, ``{Infinitely many N=2 SCFT with ADE flavor symmetry},'' \href{http://dx.doi.org/10.1007/JHEP01(2013)191}{{\em JHEP} {\bfseries 01} (2013) 191}, \href{http://arxiv.org/abs/1210.2886}{{\ttfamily arXiv:1210.2886 [hep-th]}}.

\bibitem{Cecotti:2013lda}
S.~Cecotti, M.~Del~Zotto, and S.~Giacomelli, ``{More on the N=2 superconformal systems of type $D_p(G)$},'' \href{http://dx.doi.org/10.1007/JHEP04(2013)153}{{\em JHEP} {\bfseries 04} (2013) 153}, \href{http://arxiv.org/abs/1303.3149}{{\ttfamily arXiv:1303.3149 [hep-th]}}.

\bibitem{Cecotti:2010fi}
S.~Cecotti, A.~Neitzke, and C.~Vafa, ``{R-Twisting and 4d/2d Correspondences},'' \href{http://arxiv.org/abs/1006.3435}{{\ttfamily arXiv:1006.3435 [hep-th]}}.

\bibitem{Chacaltana:2012zy}
O.~Chacaltana, J.~Distler, and Y.~Tachikawa, ``{Nilpotent orbits and codimension-two defects of 6d N=(2,0) theories},'' \href{http://dx.doi.org/10.1142/S0217751X1340006X}{{\em Int. J. Mod. Phys. A} {\bfseries 28} (2013) 1340006}, \href{http://arxiv.org/abs/1203.2930}{{\ttfamily arXiv:1203.2930 [hep-th]}}.

\bibitem{CM}
D.~H. Collingwood and W.~M. McGovern, {\em {Nilpotent Orbits in Semisimple Lie Algebras}}.
\newblock Van Nostrand Reinhold, New York, 1993.

\bibitem{Couzens:2023kyf}
C.~Couzens, M.~J. Kang, C.~Lawrie, and Y.~Lee, ``{Holographic duals of Higgsed $\mathcal{D}_p^b(BCD)$},'' \href{http://arxiv.org/abs/2312.12503}{{\ttfamily arXiv:2312.12503 [hep-th]}}.

\bibitem{Creutzig:2017qyf}
T.~Creutzig, ``{W-algebras for Argyres-Douglas theories},'' \href{http://arxiv.org/abs/1701.05926}{{\ttfamily arXiv:1701.05926 [hep-th]}}.

\bibitem{Dubrovin:1993nt}
B.~Dubrovin, ``{Differential geometry of the space of orbits of a Coxeter group},'' \href{http://arxiv.org/abs/hep-th/9303152}{{\ttfamily arXiv:hep-th/9303152}}.

\bibitem{Dubrovin:1994hc}
B.~Dubrovin, ``{Geometry of 2-D topological field theories},'' \href{http://dx.doi.org/10.1007/BFb0094793}{{\em Lect. Notes Math.} {\bfseries 1620} (1996) 120--348}, \href{http://arxiv.org/abs/hep-th/9407018}{{\ttfamily arXiv:hep-th/9407018}}.

\bibitem{DubrovinAlmost}
B.~Dubrovin, \href{http://dx.doi.org/10.1090/trans2/212/05}{``{On almost duality for Frobenius manifolds},''} in {\em Geometry, Topology, and Mathematical Physics}, vol.~212 of {\em Amer. Math. Soc. Transl. Ser. 2}, pp.~75--132.
\newblock American Mathematical Society, Providence, RI, 2004.
\newblock \href{http://arxiv.org/abs/math/0307374}{{\ttfamily arXiv:math/0307374 [math.DG]}}.

\bibitem{Dunkl}
C.~F. Dunkl, ``{Differential-difference operators associated to reflection groups},'' {\em Trans. Amer. Math. Soc.} {\bfseries 311} (1989) 167--183.

\bibitem{DJO}
C.~F. Dunkl, M.~F.~E. de~Jeu, and E.~M. Opdam, ``{Singular polynomials for finite reflection groups},'' {\em Trans. Amer. Math. Soc.} {\bfseries 346} (1994) 237--256.

\bibitem{EYY}
T.~Eguchi, Y.~Yamada, and S.-K. Yang, ``{Topological field theories and the period integrals},'' \href{http://dx.doi.org/10.1142/S0217732393001379}{{\em Mod. Phys. Lett. A} {\bfseries 8} no.~17, (1993) 1627--1637}, \href{http://arxiv.org/abs/hep-th/9304121}{{\ttfamily arXiv:hep-th/9304121}}.

\bibitem{EM}
L.~Ein and M.~Musta{\c{t}}{\u{a}}, ``{Jet schemes and singularities},'' in {\em {Algebraic Geometry--Seattle 2005, Part 2}}, vol.~80 of {\em Proc. Sympos. Pure Math.}, pp.~505--546.
\newblock Amer. Math. Soc., Providence, 2009.
\newblock \href{http://arxiv.org/abs/math/0612862}{{\ttfamily arXiv:math/0612862}}.

\bibitem{Eisenbud}
D.~Eisenbud, {\em {The Geometry of Syzygies}}, vol.~229 of {\em Graduate Texts in Mathematics}.
\newblock Springer, New York, 2005.

\bibitem{Esole:2017qeh}
M.~Esole, R.~Jagadeesan, and M.~J. Kang, ``{The Geometry of G$_2$, Spin(7), and Spin(8)-models},'' \href{http://arxiv.org/abs/1709.04913}{{\ttfamily arXiv:1709.04913 [hep-th]}}.

\bibitem{Esole:2018mqb}
M.~Esole and M.~J. Kang, ``{The Geometry of the SU(2)$\times$ G$_2$-model},'' \href{http://dx.doi.org/10.1007/JHEP02(2019)091}{{\em JHEP} {\bfseries 02} (2019) 091}, \href{http://arxiv.org/abs/1805.03214}{{\ttfamily arXiv:1805.03214 [hep-th]}}.

\bibitem{FS}
M.~Feigin and A.~Silantyev, ``{Singular polynomials from orbit spaces},'' \href{http://dx.doi.org/10.1112/S0010437X1200036X}{{\em Compos. Math.} {\bfseries 148} no.~6, (2012) 1867--1879}, \href{http://arxiv.org/abs/1110.1946}{{\ttfamily arXiv:1110.1946 [math.RT]}}.

\bibitem{Gadde:2011uv}
A.~Gadde, L.~Rastelli, S.~S. Razamat, and W.~Yan, ``{Gauge Theories and Macdonald Polynomials},'' \href{http://dx.doi.org/10.1007/s00220-012-1607-8}{{\em Commun. Math. Phys.} {\bfseries 319} (2013) 147--193}, \href{http://arxiv.org/abs/1110.3740}{{\ttfamily arXiv:1110.3740 [hep-th]}}.

\bibitem{GG}
W.~L. Gan and V.~Ginzburg, ``{Quantization of Slodowy slices},'' {\em Int. Math. Res. Not.} {\bfseries 2002} no.~5, (2002) 243--255, \href{http://arxiv.org/abs/math/0105225}{{\ttfamily arXiv:math/0105225}}.

\bibitem{GW}
R.~Goodman and N.~R. Wallach, {\em {Symmetry, Representations, and Invariants}}, vol.~255 of {\em Graduate Texts in Mathematics}.
\newblock Springer, Dordrecht, 2009.

\bibitem{Humphreys}
J.~E. Humphreys, {\em {Reflection Groups and Coxeter Groups}}, vol.~29 of {\em Cambridge Studies in Advanced Mathematics}.
\newblock Cambridge University Press, Cambridge, 1990.

\bibitem{JS}
B.~Johnson and E.~Sommers, ``{Equations for some nilpotent varieties},'' \href{http://dx.doi.org/10.1093/imrn/rny132}{{\em Int. Math. Res. Not. IMRN} {\bfseries 2020} no.~14, (2020) 4433--4464}, \href{http://arxiv.org/abs/1706.04820}{{\ttfamily arXiv:1706.04820 [math.RT]}}.

\bibitem{KW88}
V.~G. Kac and M.~Wakimoto, ``{Modular invariant representations of infinite-dimensional Lie algebras and superalgebras},'' {\em Proc. Natl. Acad. Sci. USA} {\bfseries 85} (1988) 4956--4960.

\bibitem{KWb}
V.~G. Kac and M.~Wakimoto, ``{A remark on boundary level admissible representations},'' {\em C. R. Math. Acad. Sci. Paris} {\bfseries 355} (2017) 128--132, \href{http://arxiv.org/abs/1612.07423}{{\ttfamily arXiv:1612.07423 [math.RT]}}.

\bibitem{Kang:2022zsl}
M.~J. Kang, C.~Lawrie, K.-H. Lee, M.~Sacchi, and J.~Song, ``{Higgs branch, Coulomb branch, and Hall-Littlewood index},'' \href{http://dx.doi.org/10.1103/PhysRevD.106.106021}{{\em Phys. Rev. D} {\bfseries 106} no.~10, (2022) 106021}, \href{http://arxiv.org/abs/2207.05764}{{\ttfamily arXiv:2207.05764 [hep-th]}}.

\bibitem{Kang:2021ccs}
M.~J. Kang, C.~Lawrie, K.-H. Lee, and J.~Song, ``{Infinitely many 4D N=1 SCFTs with a=c},'' \href{http://dx.doi.org/10.1103/PhysRevD.105.126006}{{\em Phys. Rev. D} {\bfseries 105} no.~12, (2022) 126006}, \href{http://arxiv.org/abs/2111.12092}{{\ttfamily arXiv:2111.12092 [hep-th]}}.

\bibitem{Kang:2023dsa}
M.~J. Kang, C.~Lawrie, K.-H. Lee, and J.~Song, ``{Emergent N=4 Supersymmetry from N=1},'' \href{http://dx.doi.org/10.1103/PhysRevLett.130.231601}{{\em Phys. Rev. Lett.} {\bfseries 130} no.~23, (2023) 231601}, \href{http://arxiv.org/abs/2302.06622}{{\ttfamily arXiv:2302.06622 [hep-th]}}.

\bibitem{Kang:2022vab}
M.~J. Kang, C.~Lawrie, K.-H. Lee, and J.~Song, ``{Operator spectroscopy for 4D SCFTs with a=c},'' \href{http://dx.doi.org/10.1103/PhysRevD.107.066018}{{\em Phys. Rev. D} {\bfseries 107} no.~6, (2023) 066018}, \href{http://arxiv.org/abs/2210.06497}{{\ttfamily arXiv:2210.06497 [hep-th]}}.

\bibitem{Kang:2024elv}
M.~J. Kang, C.~Lawrie, K.-H. Lee, and J.~Song, ``{Generalized symmetry constraints on deformed 4D SCFTs},'' \href{http://dx.doi.org/10.1103/PhysRevD.111.086028}{{\em Phys. Rev. D} {\bfseries 111} no.~8, (2025) 086028}, \href{http://arxiv.org/abs/2408.14532}{{\ttfamily arXiv:2408.14532 [hep-th]}}.

\bibitem{Kang:2024inx}
M.~J. Kang, C.~Lawrie, K.-H. Lee, and J.~Song, ``{Landscape of 4d N=1 SCFTs with a=c},'' \href{http://dx.doi.org/10.1103/PhysRevD.111.086002}{{\em Phys. Rev. D} {\bfseries 111} no.~8, (2025) 086002}, \href{http://arxiv.org/abs/2412.17895}{{\ttfamily arXiv:2412.17895 [hep-th]}}.

\bibitem{Kang:2021lic}
M.~J. Kang, C.~Lawrie, and J.~Song, ``{Infinitely many 4D N=2 SCFTs with a=c and beyond},'' \href{http://dx.doi.org/10.1103/PhysRevD.104.105005}{{\em Phys. Rev. D} {\bfseries 104} no.~10, (2021) 105005}, \href{http://arxiv.org/abs/2106.12579}{{\ttfamily arXiv:2106.12579 [hep-th]}}.

\bibitem{Kang:2025zub}
M.~J. Kang, C.~Lawrie, and J.~Song, ``{Index from a point},'' \href{http://arxiv.org/abs/2507.12510}{{\ttfamily arXiv:2507.12510 [hep-th]}}.

\bibitem{Kang:2026nge}
M.~J. Kang, C.~Lawrie, and J.~Song, ``{Algebro-geometric bootstrapping from OPE decoupling},'' \href{http://arxiv.org/abs/2602.04946}{{\ttfamily arXiv:2602.04946 [hep-th]}}.

\bibitem{atlas}
M.~J. Kang, C.~Lawrie, and J.~Song, ``{ArcsToSCFT Theory Atlas},'' \url{https://jws611.github.io/ArcsToSCFT-atlas/}, 2026.
\newblock Records for \href{https://jws611.github.io/ArcsToSCFT-atlas/data/theories/D4-p5-Y0.json}{$D_5(SO(8))$} and \href{https://jws611.github.io/ArcsToSCFT-atlas/data/theories/D4-p7-Y0.json}{$D_7(SO(8))$}, $Y=0$, consulted November 2, 2026.

\bibitem{Kang:2026ney}
M.~J. Kang, C.~Lawrie, and J.~Song, ``{One ring to rule them all: syzygy maximization for Argyres-Douglas},'' \href{http://arxiv.org/abs/2610.33151}{{\ttfamily arXiv:2610.33151 [hep-th]}}.

\bibitem{KLSnonpoint}
M.~J. Kang, C.~Lawrie, and J.~Song, ``{The index from a nonpoint and RG flows},'' \href{http://arxiv.org/abs/2609.13374}{{\ttfamily arXiv:2609.13374 [hep-th]}}.

\bibitem{Kinney:2005ej}
J.~Kinney, J.~M. Maldacena, S.~Minwalla, and S.~Raju, ``{An Index for 4 dimensional super conformal theories},'' \href{http://dx.doi.org/10.1007/s00220-007-0258-7}{{\em Commun. Math. Phys.} {\bfseries 275} (2007) 209--254}, \href{http://arxiv.org/abs/hep-th/0510251}{{\ttfamily arXiv:hep-th/0510251}}.

\bibitem{Kostant}
B.~Kostant, ``{Lie group representations on polynomial rings},'' \href{http://dx.doi.org/10.2307/2373130}{{\em Amer. J. Math.} {\bfseries 85} (1963) 327--404}.

\bibitem{LNS}
M.~Lehn, Y.~Namikawa, and C.~Sorger, ``{Slodowy slices and universal Poisson deformations},'' \href{http://dx.doi.org/10.1112/S0010437X11005550}{{\em Compos. Math.} {\bfseries 148} (2012) 121--144}, \href{http://arxiv.org/abs/1002.4107}{{\ttfamily arXiv:1002.4107 [math.AG]}}.

\bibitem{Li}
H.~Li, ``{Abelianizing vertex algebras},'' {\em Commun. Math. Phys.} {\bfseries 259} (2005) 391--411, \href{http://arxiv.org/abs/math/0409140}{{\ttfamily arXiv:math/0409140}}.

\bibitem{MW}
M.~Mikosz and A.~Weber, ``{Triality, characteristic classes, $D_4$ and $G_2$ singularities},'' \href{http://dx.doi.org/10.1007/s40062-014-0094-3}{{\em J. Homotopy Relat. Struct.} {\bfseries 10} (2015) 995--1011}, \href{http://arxiv.org/abs/1311.0507}{{\ttfamily arXiv:1311.0507 [math.AG]}}.

\bibitem{Mustata}
M.~Musta{\c{t}}{\u{a}}, ``{Jet schemes of locally complete intersection canonical singularities},'' \href{http://dx.doi.org/10.1007/s002220100152}{{\em Invent. Math.} {\bfseries 145} (2001) 397--424}, \href{http://arxiv.org/abs/math/0008002}{{\ttfamily arXiv:math/0008002}}. With an appendix by D. Eisenbud and E. Frenkel.

\bibitem{Saito}
K.~Saito, ``{On a linear structure of the quotient variety by a finite reflexion group},'' {\em Publ. Res. Inst. Math. Sci.} {\bfseries 29} (1993) 535--579.

\bibitem{Slodowy}
P.~Slodowy, \href{http://dx.doi.org/10.1007/BFb0090294}{{\em {Simple Singularities and Simple Algebraic Groups}}}, vol.~815 of {\em Lecture Notes in Mathematics}.
\newblock Springer, Berlin, 1980.

\bibitem{Song:2016yfd}
J.~Song, ``{Macdonald Index and Chiral Algebra},'' \href{http://dx.doi.org/10.1007/JHEP08(2017)044}{{\em JHEP} {\bfseries 08} (2017) 044}, \href{http://arxiv.org/abs/1612.08956}{{\ttfamily arXiv:1612.08956 [hep-th]}}.

\bibitem{Song:2017oew}
J.~Song, D.~Xie, and W.~Yan, ``{Vertex operator algebras of Argyres-Douglas theories from M5-branes},'' \href{http://dx.doi.org/10.1007/JHEP12(2017)123}{{\em JHEP} {\bfseries 12} (2017) 123}, \href{http://arxiv.org/abs/1706.01607}{{\ttfamily arXiv:1706.01607 [hep-th]}}.

\bibitem{stacks-project}
{The Stacks Project Authors}, ``{The Stacks Project},'' \url{https://stacks.math.columbia.edu}, 2026.

\bibitem{Weyl}
H.~Weyl, {\em {The Classical Groups: Their Invariants and Representations}}.
\newblock Princeton University Press, Princeton, 1939.
\newblock Second edition, 1946.

\bibitem{Xie:2012hs}
D.~Xie, ``{General Argyres-Douglas Theory},'' \href{http://dx.doi.org/10.1007/JHEP01(2013)100}{{\em JHEP} {\bfseries 01} (2013) 100}, \href{http://arxiv.org/abs/1204.2270}{{\ttfamily arXiv:1204.2270 [hep-th]}}.

\bibitem{Zhu}
Y.~Zhu, ``{Modular invariance of characters of vertex operator algebras},'' {\em J. Amer. Math. Soc.} {\bfseries 9} (1996) 237--302.

\end{thebibliography}\endgroup
\bibliographystyle{sortedbutpretty}
\end{document}